\documentclass[12pt]{article}
\usepackage{amssymb,latexsym,amsmath,amsthm,amsfonts,txfonts,bm}

\newtheorem{theorem}{Theorem}[subsection]
\newtheorem{proposition}[theorem]{Proposition}
\newtheorem{corollary}[theorem]{Corollary}
\newtheorem{conjecture}[theorem]{Conjecture}
\newtheorem{lemma}[theorem]{Lemma}

\theoremstyle{remark}
\newtheorem{remark}[theorem]{Remark}
\newtheorem{example}[theorem]{Example}
\newtheorem{construction}[theorem]{Construction}
\theoremstyle{definition}
\newtheorem{definition}[theorem]{Definition}

\numberwithin{equation}{subsection}

\DeclareMathOperator{\Aut}{Aut}

\DeclareMathOperator{\Frob}{Frob}

\DeclareMathOperator{\Gal}{Gal}
\DeclareMathOperator{\Hom}{Hom}

\DeclareMathOperator{\coker}{coker}
\DeclareMathOperator{\im}{im}
\DeclareMathOperator{\id}{id}

\newcommand{\A}{\mathbb A}

\newcommand{\C}{\mathbb C}
\newcommand{\Fqbar}{{\overline \F_q}}
\newcommand{\F}{\mathbb F}
\newcommand{\Gm}{{\mathbb G}_m}

\newcommand{\Spec}{{\rm Spec}}
\newcommand{\Sch}{{\rm Sch}}
\newcommand{\Zeta}{{\rm Zeta}}

\renewcommand{\star}{*}

\newcommand{\M}{\mathcal M}
\newcommand{\tM}{\tilde{\M}}
\newcommand{\Mp}{\M^\prime}
\newcommand{\Unp}{{\mathcal U}_n^\prime}
\newcommand{\N}{\mathbb N}
\newcommand{\Nn}{{\mathcal N}_n}

\newcommand{\Oz}{\mathcal O}

\newcommand{\Q}{\mathbb Q}
\newcommand{\R}{\mathbb R}
\newcommand{\X}{\mathcal X}
\newcommand{\Z}{\mathbb Z}
\newcommand{\U}{\mathbb U}
\newcommand{\Un}{{\mathcal U}_n}

\newcommand{\bfP}{\mathbb P}

\newcommand{\epolc}{e}
\newcommand{\epolnc}{\overline e}
\newcommand{\epoln}{\overline E}
\newcommand{\epol}{E}

\newcommand{\hpolnc}{\overline h}
\newcommand{\hpoln}{\overline H}
\newcommand{\h}{\mathcal H}

\newcommand{\mubf}{{\bbmu}}
\newcommand{\p}{\mathcal P}

\DeclareMathOperator{\Exp}{Exp}
\DeclareMathOperator{\Log}{Log}
\DeclareMathOperator{\tr}{Tr}

{.tfm}
{.tfm}
\def\beq{\begin{eqnarray}}
\def\eeq{\end{eqnarray}}
\def\bes{\begin{eqnarray*}}
\def\ees{\end{eqnarray*}}

\def\C{\mathbb{C}}

\def\R{\mathbb{R}}
\def\F{\mathbb{F}}        
\def\Q{\mathbb{Q}}

\def\gl{{\mathfrak g\mathfrak l}}
\def\tdn{\tilde{d}_n}

\def\Z{\mathbb{Z}}
\def\N{\mathbb{N}}
\def\K{\mathbb{K}}

\newcommand{\bino}[2]{\mbox{ $#1 \choose #2$}}
\newcommand{\nc}{\newcommand}
\nc{\bbmu}{{\bm \mu}}
\nc{\op}[1]{\mathop{\mathchoice{\mbox{\rm #1}}{\mbox{\rm #1}}
{\mbox{\rm \scriptsize #1}}{\mbox{\rm \tiny #1}}}\nolimits}
\nc{\al}{\alpha}
\newcommand{\be}{\beta}
\nc{\ep}{\varepsilon}
\nc{\ga}{\gamma}
\nc{\Ga}{\Gamma}
\nc{\La}{\Lambda}
\nc{\si}{\sigma}
\nc{\Sig}{{\Gamma}}
\nc{\Om}{\Omega}
\nc{\om}{\omega}
\nc{\SL}[1]{{{\rm SL(}#1{\rm )}}}
\nc{\GL}[1]{{{\rm GL(}#1{\rm )}}}
\nc{\PGL}[1]{{{\rm PGL(}#1{\rm )}}}

\nc{\cpt}{{\op{cpt}}}
\nc{\Dol}{{\op{Dol}}}
\nc{\DR}{{\op{DR}}}
\nc{\Triv}{\op{Triv}}
\nc{\Hod}{{\op{Hod}}}
\nc{\Est}{E_{\op{st}}}
\nc{\Hst}{H_{\op{st}}}
\nc{\Left}[1]{\hbox{$\left#1\vbox to
   10.5pt{}\right.\nulldelimiterspace=0pt \mathsurround=0pt$}}
\nc{\Right}[1]{\hbox{$\left.\vbox to
   10.5pt{}\right#1\nulldelimiterspace=0pt \mathsurround=0pt$}}
\nc{\LEFT}[1]{\hbox{$\left#1\vbox to
   15.5pt{}\right.\nulldelimiterspace=0pt \mathsurround=0pt$}}
\nc{\RIGHT}[1]{\hbox{$\left.\vbox to
   15.5pt{}\right#1\nulldelimiterspace=0pt \mathsurround=0pt$}}
\nc{\bee}{{\bf E}}
\nc{\bphi}{{\bf \Phi}}

\begin{document} 

\title{Mixed Hodge polynomials of character varieties } 
\author{ Tam\'as Hausel
\\ {\it Mathematical Institute}
\\ {\it 24-29 St. Giles'}
\\ {\it Oxford, OX1 3LB, UK}
 \\ {\tt hausel@maths.ox.ac.uk} \and Fernando Rodriguez-Villegas \\{\it Department of Mathematics}\\ {\it University Station C1200}
 \\{\it Austin, Texas, 78712, USA} \\ {\tt
   villegas@math.utexas.edu}   \and \\ {with an appendix by Nicholas M. Katz} }  
 \maketitle
   \begin{abstract} We calculate the $E$-polynomials of certain twisted
$\GL{n,\C}$-character varieties  $\M_n$ of Riemann surfaces by counting
points over finite fields using the character table of the finite group of Lie-type $\GL{n,\F_q}$
and a theorem proved in the appendix by N. Katz.  We
deduce from this calculation several geometric results, for example,  the value of the topological Euler characteristic of the associated
$\PGL{n,\C}$-character variety.
The calculation also leads to several conjectures about the cohomology of  $\M_n$: an
explicit conjecture for its mixed Hodge polynomial; a conjectured
curious Hard Lefschetz theorem and a conjecture relating the pure part
to absolutely indecomposable representations of a certain quiver. 
We prove these conjectures for $n=2$.
\end{abstract}\newpage
\tableofcontents
\newpage

\section{Introduction}
\stepcounter{subsection}
Let $g\geq 0$ and $n>0$ be integers. Let $\zeta_n\in\C$ be a primitive $n$-th root of unity. Abbreviating $[A,B]=ABA^{-1}B^{-1}$ and denoting the identity matrix $I_n\in\GL{n,\C}$ we define 
\beq
 \M_n:=\{\,\, A_1,B_1,\dots, A_g,B_g \in \GL{n,\C} \,\,|\,\,
[A_1,B_1] \dots [A_g,B_g] = \zeta_n I_n \}/\!/\GL{n,\C} \label{twisted}  \eeq an affine GIT quotient by the conjugation action of $\GL{n,\C}$. It is a twisted character variety of a genus $g$ closed  
 Riemann surface $\Sigma$; its points can be thought of twisted
 homomorphisms of $\pi_1(\Sigma) \to \GL{n,\C}$ modulo conjugation. It  is a non-singular affine variety of dimension $d_n:=n^2(2g-2)+2$ by  Theorem~\ref{smooth}.

  One of the goals of this paper is to find the Poincar\'e polynomial $P(\M_n;t)=\sum_i b_i(\M_n) t^i$ which encodes the Betti numbers $b_i(\M_n)$ of $\M_n$. They were calculated for $n=2$ by Hitchin \cite{hitchin} and for $n=3$ by Gothen \cite{gothen}. To be precise, Hitchin and Gothen work with a certain moduli space of Higgs bundles on $\Sigma$, which is known to be diffeomorphic to $\M_n$ by non-Abelian Hodge theory \cite{hitchin, simpson}. On the other hand, the Poincar\'e polynomial of the ${\rm U}(n)$-character variety $\Nn^d$ of $\Sigma$, where $\GL{n,\C}$ is replaced by ${\rm U}(n)$ in the above definition and $\zeta_n=\exp(\frac{d}{n}{2\pi i})$, were obtained by Harder-Narasimhan \cite{harder-narasimhan}, using the Weil conjectures proved by Deligne \cite{De-Weil I}, and by Atiyah-Bott \cite{atiyah-bott}  using gauge theory. An explicit
 closed formula for the Poincar\'e polynomial of the ${\rm U}(n)$-character varieties was given by Zagier \cite{zagier}.

 Other motivations to study the cohomology of character varieties  are discussed in \cite{hausel3}, which also announces many of the results of this paper. Character varieties  appear in the Geometric Langlands program of Beilinson-Drinfeld \cite{beilinson-drinfeld}. Recently, many new ideas relating physics, in particular mirror symmetry,  to the Geometric Langlands program have been discussed by Kapustin-Witten in \cite{kapustin-witten}.   One can expect \cite{hausel3, hausel4.5} that the results of this paper will have  analogues for $\SL{n,\C}$ character varieties reflecting the expected relationship between certain Hodge numbers
 of $\PGL{n,\C}$ and $\SL{n,\C}$ character varieties dictated by mirror symmetry considerations.

 In this paper we  
  study  the mixed Hodge polynomials $H(\M_n;q,t)$ and uncover a surprising amount  of structure governing them. 
 The mixed Hodge polynomial is a common 
 deformation of the Poincar\'e polynomial $P(\M_n;t)=H(\M_n;1,t)$ 
 and the so-called $E$-polynomial $E(\M_n;q)=q^{d_n}H(\M_n;1/q,-1)$ and is defined using Deligne's construction of mixed Hodge structures on the cohomology of a complex algebraic variety \cite{De-Hodge II,De-Hodge III} (see Subsection~\ref{mixedhodge}). 

 We explicitly calculate the $E$-polynomial of $\M_n$ in terms of a generating function
 using arithmetic algebraic geometry. One key result used in this calculation 
 is Theorem~\ref{katz}.3  of Katz in the
 appendix, which basically says that if the number of points of a variety over every  finite field $\F_q$ is a polynomial
 in $q$ then this polynomial agrees with the $E$-polynomial of the variety.
 
  Another ingredient is a well-known character
 formula, Theorem~\ref{char-sum-z}, which counts the number of solutions of certain equations in a finite group.  Similar counting formulas go back to the birth of character theory of finite groups  by Frobenius \cite{frobenius} in 1896. Combining these and Corollary~\ref{polynomialcount} we get \beq E(\M_n;q)= \sum_{\chi \in {\rm Irr(\GL{n,\F_q})}} 
\frac{| \GL{n,\F_q} |^{2g-2}}{\chi(I_n)^{2g-1}} \chi(\zeta_n I_n) ,\label{charform}. \eeq The character table
of $\GL{n,\F_q}$ was determined by Green \cite{green} in 1955. Using Green's results, the
evaluation of the formula (\ref{charform}) is carried out in \S~\ref{epolynomial}.  The calculation
makes non-trivial use of the inclusion-exclusion principle for the poset of finite set partitions.  The end result is  an expression for  the $E$-polynomials in terms of  an explicit generating function in Theorem~\ref{main}. An important consequence of our Theorem~\ref{main} is that the number of points of the variety $\M_n$ over the finite field $\F_q$ is a polynomial in $q$. For example, for $n=2$ we prove in Corollary~\ref{e2} that   
\beq E(\M_2;q)/(q-1)^{2g}=(q^2-1)^{2g-2}+q^{2g-2}(q^2-1)^{2g-2}-\frac{1}{2} q^{2g-2} (q-1)^{2g-2}-\frac{1}{2}q^{2g-2}(q+1)^{2g-2}.\label{e2q}\eeq 

An interesting topological outcome  of  our calculation is  the precise value  of the  Euler characteristic of our character varieties. The variety $\M_n$ is cohomologically a product of $(\C^\times)^{2g}$ and the $\PGL{n,\C}$-character variety $\tM_n:=\M_n/\!/(\C^\times)^{2g}$, which is  defined as the quotient of $\M_n$ by the natural action of the torus $(\C^\times)^{2g}$ on 
$\M_n$. Therefore, the Euler characteristic of $\M_n$ is $0$ due to the fact that the Euler characteristic of the torus $(\C^\times)^{2g}$ is $0$. However the Euler characteristic of $\tM_n$ is more interesting (see \S \ref{eulercharacteristic}):
\begin{corollary} \label{euler-char} Let $g>1$. The Euler characteristic of the $\PGL{n,\C}$-character variety $\tM_n$ is $\mu(n) n^{2g-3}$, where $\mu$ is the  M\"obius function.  
\end{corollary}
\noindent The last result of the first part of this paper is a formula for the number of points on the untwisted $\GL{n}$ character variety  (which is defined by replacing $\zeta_n I$  by $I$ in the definition \eqref{twisted} ) over a finite field $\F_q$. Our explicit generating function formula in Theorem~\ref{untwisted} could be interesting to compare
with recent work of Liebeck-Shalev  \cite{liebeck-shalev} studying
asymptotics of the same quantities. 


The second part of this paper concerns  the mixed Hodge
 polynomial of $\M_n$. In Conjecture~\ref{main-conj} we give a formula for it as a natural $t$-deformation of our calculation of the $E$-polynomial of $\M_n$.
 Here we only give our conjecture in the case $g=1$. In
 this case we know {\em a priori} that our character variety is $\M_1=(\C^\times)^2$, the $2$-torus (see Theorem~\ref{g=1-char-var}). Therefore our conjecture becomes a purely combinatorial statement.   \begin{conjecture} The following combinatorial identity holds:
  \begin{equation*}
\sum_\lambda
\prod
\frac{\left(z^{2a+1}-w^{2l+1}\right)^2}
{(z^{2a+2}-w^{2l})(z^{2a}-w^{2l+2})}\,T^{|\lambda|}=\exp\left( \sum_{k\geq 1}
 \frac{(z^k-w^k)^2}{(z^{2k}-1)(1-w^{2k})(1-T^k)}\frac{T^k}k\right) ,
 \end{equation*}
where the sum on the left hand side is over all partitions $\lambda$, and the product
is over all boxes in the Ferrers diagram of $\lambda$, and $a$ and $l$ are its arm and leg-length, 
as defined in \S\ref{partitions}.
 \end{conjecture} 
 This yet unproven identity is reminiscent of  the Macdonald
 identities and the Weyl-Kac character formula; it is conceivable that it has a representation theory
 interpretation. For example, the corresponding formula in the $g=0$ 
 case will be proved as Theorem~\ref{g=0iden} using a result of Garsia-Haiman \cite{garsia-haiman} obtained from the study of Macdonald polynomials. This presently mysterious link between mixed Hodge polynomials of character varieties and Macdonald polynomials is further
 developed in \cite{hausel5, hausel-letellier-villegas}. In particular, the main conjecture of \cite{hausel-letellier-villegas} says
 that the mixed Hodge polynomials of character varieties of Riemann surfaces with semisimple conjugacy classes at the punctures are 
 governed by Macdonald polynomials in a simple way, its structure resembling  a topological quantum field theory.

 We study many implications of our main Conjecture~\ref{main-conj} and prove several
 consistency results. Because we have an explicit description  of the cohomology ring of
 $\M_2$ given in \cite{hausel-thaddeus-generators,hausel-thaddeus-relations} we are able to determine the mixed Hodge polynomial in the $n=2$ case and confirm all our conjectures.  

 Rather than giving a full description of our conjectures for general $n$ here, we present instead the corresponding theorems in the $n=2$
 case.

 \begin{theorem} \label{mhp2} The mixed Hodge polynomials of $\tM_2$ and $\M_2$ are   given by \beq  \label{mhp2-fmla} H(\tM_2;q,t)=\frac{H(\M_2;q,t)}{(qt+1)^{2g}}&=& \frac{(q^2t^3+1)^{2g}}{(q^2t^2-1)(q^2t^4-1)}+
\frac{q^{2g-2}t^{4g-4}(q^2t+1)^{2g}}{(q^2-1)(q^2t^2-1)}-\nonumber \\
&&-\frac{1}{2}\frac{q^{2g-2}t^{4g-4}(qt+1)^{2g}}{(qt^2-1)(q-1)}-\frac{1}{2}
\frac{ q^{2g-2}t^{4g-4}(qt-1)^{2g}}{(q+1)(qt^2+1)}\end{eqnarray} 
\end{theorem}
By setting $t=-1$ in this formula we recover the $E$-polynomial
in \eqref{e2q}. Thus, by purely cohomological calculations on $\M_2$ we derive formula \eqref{e2q}, which reflects the structure
of irreducible characters of $\GL{2,\F_q}$. For example, the four terms
above correspond to the four types of irreducible characters of $\GL{2,\F_q}$.  Looking at the other specialization $q=1$ gives a pleasant formula for the Poincar\'e polynomial  
$P(\M_2,t)=H(\M_2;1,t)$, which agrees with Hitchin's calculation \cite{hitchin}.  

Changing $q$ by $1/qt^2$ in the right hand side of (\ref{mhp2-fmla})  interchanges the  first two terms and fixes the other two. This implies the following 
\begin{corollary} \label{curious2}The mixed Hodge polynomial of $\tM_2$ satisfies the following {\em curious Poincar\'e duality}: \beq \label{curiouspoin-fmla}H(\tM_2;1/qt^2,t)=(qt)^{-\dim \tM_2} H(\tM_2;q,t)\eeq   
\end{corollary}
In fact, we can give a geometrical interpretation of this combinatorial observation. First, $H^2(\tM_2)$ is one dimensional generated by a class $\alpha$.  Define the Lefschetz map $L:H^i(\tM_2)\to H^{i+2}(\tM_2)$ by
$x\mapsto \alpha\cup x$. As it respects mixed Hodge structures and $\alpha$ has weight $4$ it defines a map on the graded pieces of the weight filtration 
$L: Gr^W_{l} H^i(\tM_2)\to Gr^W_{l+4}H^{i+2}(\tM_2)$. In \S\ref{curioushard}  we prove the following {\em curious Hard Lefschetz} 
\begin{theorem} \label{lefschetz} The Lefschetz map  \label{lefschetz-intr} $$L^{l}:Gr^W_{6g-6-2l} H^{i-l}(\tM_2)\to Gr^W_{6g-6+2l} H^{i+l}(\tM_2)$$ is an isomorphism.  
\end{theorem} 
\noindent The agreement of the dimensions of these two isomorphic vector spaces is equivalent to (\ref{curiouspoin-fmla}).

Interestingly, this theorem implies (see Remark~\ref{curioushardremark}) a theorem of \cite{hausel2} 
that the Lefschetz map $L^k:H^{\tilde{d}_2/2-k}(\tM_2)\to H^{\tilde{d}_2/2+k}(\tM_2)$  is injective; where $\tilde{d}_2=\dim \tM_2=6g-6$. As it is explained 
in \cite{hausel2}  this weak version of Hard Lefschetz applied to toric
hyperk\"ahler varieties yields
new inequalities for the $h$-numbers of matroids. See also
\cite{hausel-sturmfels} for the original argument on toric hyperk\"ahler 
varieties. 
Theorem~\ref{lefschetz-intr}  can also be thought of as an analogue
of the Faber 
conjecture \cite{faber} on the cohomology of the moduli space of curves,
which is another non-compact variety whose cohomology ring is conjectured 
to satisfy a certain form of the Hard Lefschetz theorem.  

For any smooth variety $X$ there is an important subring of $H^*(X)$, namely the so-called pure ring $PH^*(X)\cong \oplus_{k} W_{k} H^k(X)$. We denote by $PP(X;t)$ the 
Poincar\'e polynomial of the pure ring. We can obtain $PP(X;t)$ from $H(X;q,t)$
by taking the monomials which are powers of $qt^2$.  In the case of $\M_2$ the pure ring is generated by a single class $\beta\in H^4(\M_2)$ and with one relation $\beta^g=0$, the so-called Newstead relation. Consequently $PP(\M_2;t)=1+t^4+\dots+t^{4g-4}.$ This implies the following:

\begin{theorem} \label{pp2} Let $A_n(q)$ be the number of absolutely indecomposable $g$-tuples of $n$ by $n$  matrices over the 
finite field $\F_q$ modulo conjugation. Then for $n=2$ we have
$$PP(\M_2;\sqrt{q})=q^{d_n/2} A_2(1/q)$$
\end{theorem}
In \S\ref{purepart} we conjecture the same for any $n$. The function $A_n(q)$  is an instance of the $A$-polynomial defined
by Kac \cite{kac} for any quiver. The quiver here is $S_{\! g}$,  $g$ loops on one vertex. Kac showed that the $A$ function is always a polynomial and conjectured it has non-negative coefficients. When the dimension vector is  indivisible  this has been proved by Crawley-Boevey and Van den Bergh \cite{crawley-boevey-etal} by giving a cohomological interpretation of the $A$-polynomial. For  $S_{\! g}$ the
result of \cite{crawley-boevey-etal} only applies in the   $n=1$ case, since all other dimension vectors are divisible. Our Theorem~\ref{pp2} shows a cohomological interpretation for $A_2(q)$. For general $n$  Theorem~\ref{purity} together with our main Conjecture~\ref{main-conj} will then give a conjectural cohomological interpretation of  $A_n(q)$, implying Kac's conjecture for $S_{\! g}$.

Theorem~\ref{pp2} implies that the middle dimensional cohomology 
$H^{6g-6}(\tM_2)$ has a trivial pure part. It follows that  the middle dimensional compactly supported cohomology also has trivial pure part. This implies the following theorem, which will be proved
in Corollary~\ref{corollarym2intersect}.
\begin{theorem} \label{sen} The intersection form on middle dimensional compactly supported cohomology $H^{6g-6}_{c}(\tM_2)$ is trivial or equivalently the forgetful map $H^{6g-6}_{c}(\tM_2)\to H^{6g-6}(\tM_2)$ is $0$. \end{theorem}
\noindent This was the main result of \cite{hausel1} and was interpreted there as the vanishing of "topological $L^2$ cohomology" for $\tM_2$. It is surprising that we can deduce this result only from
the knowledge \cite{hausel-thaddeus-generators, hausel-thaddeus-relations} of the  structure  of the ordinary cohomology ring $H^*(\tM_2)$ and the study of its mixed Hodge structure. In fact, we only need to know that the famous Newstead relation $\beta^g=0$ holds in $H^*(\tM_2)$. See \cite{hausel4.5} for a more detailed discussion on the background and various ramifications of Theorem~\ref{sen}.

 The structure of the paper is as follows. In \S~\ref{preliminaries} we collect various
 facts which we will need later. In \S~\ref{mixedhodge} we define
 and list properties of the mixed Hodge polynomials obtained from Deligne's mixed Hodge structure. In \S~\ref{character} we define and prove the basic properties of the character varieties we study. 
 In \S~\ref{count} we derive a classical character formula for the number of solutions of a certain equation over finite groups. 
 In \S~\ref{partitions} we collect the definitions and notations for partitions which will be used throughout the paper. In \S~\ref{formal-inf-prod}
we introduce a formalism to handle various formal infinite products. 
Then in \S~\ref{epolynomial} we calculate the $E$-polynomial of our variety. In \S~\ref{mixedhodgen} we formulate our main conjecture on
the mixed Hodge polynomial of our character varieties and derive various consequences, several of which we can test for consistency. 
We also relate our conjectured mixed Hodge polynomial to  Kac's $A$-polynomial in \S~\ref{purepart}. Finally in \S~\ref{mixedhodge2} we prove all our conjectures in the $n=2$ case. 

\begin{paragraph} {\bf Acknowledgment} The first author was supported by 
NSF grants DMS-0305505 and DMS-0604775   an Alfred Sloan Fellowship and a Royal Society
University Research Fellowship. The second author was supported by an NSF grant DMS-0200605. We would like to thank Dan Freed for
bringing us together and for useful discussions.  
We would like to thank Nick Katz for his invaluable help and writing the appendix. We would also like to thank Daniel Allcock, William 
Crawley-Boevey, Mark Haiman, Gergely Harcos, Lisa Jeffrey, Sean Keel, Eckhard Meinrenken,  Martin Olsson, Tony Pantev, Alexander Postnikov, Nicholas Proudfoot, Graeme Segal, Michael Thaddeus and the referees for various comments and discussions.  Various calculations in the paper were assisted by the symbolic algebra packages  Maple, Macaulay 2 and PARI-GP. 
\end{paragraph}

\section{Preliminaries}
\label{preliminaries}

\subsection{Mixed Hodge polynomials}
\label{mixedhodge}
Motivated by the (then still
unproven) Weil Conjectures and Grothendieck's "yoga of weights", which drew 
cohomological conclusions about complex varieties from the truth of those
conjectures, Deligne in \cite{De-Hodge II, De-Hodge III} proved the existence of Mixed Hodge structures on the cohomology
of a complex algebraic variety.
\begin{proposition}[Deligne \cite{De-Hodge II, De-Hodge III}] Let X be a complex algebraic variety. For each $j$ there is an increasing weight
filtration
$${0} = W_{-1}  \subseteq W_0   \subseteq \dots \subseteq 
W_{2j} = H^j(X,\Q)$$
and a decreasing Hodge filtration
$$H^j(X,\C) = F^0  \supseteq  F^1  \supseteq \dots   \supseteq F^m \supseteq    F^{m+1} = {0}$$
such that the filtration induced by F on the complexification of the graded pieces $Gr_l^W:=W_{l}/W_{l-1}$ of the weight filtration
endows every graded piece with a pure Hodge structure of weight $l$, or equivalently for every $0\leq p \leq l$ we have \beq\label{conjugates} Gr^{W^\C}_l=F^p Gr^{W^\C}_l \oplus \overline{F^{l-p+1}Gr^{W^\C}_l}.\eeq

\end{proposition}

We  now list properties of this mixed Hodge structure, which we will need in this paper. From now on we use the notation $H^*(X)$  for $H^*(X,\Q)$.
\begin{theorem} \label{mixedhodgeproperties} \begin{enumerate} 
\item The map $f^*:H^*(Y)\to H^*(X) $, induced by an algebraic map $f:X\to Y$, strictly preserves mixed Hodge structures.
\item \label{galois} A field automorphism $\sigma:\C\to \C$ induces an isomorphism 
$H^*(X)\cong H^*(X^\sigma)$, which preserves the
mixed Hodge structure.
\item \label{kunneth} The K\"unneth isomorphism $$H^*(X\times Y)\cong H^*(X)\otimes H^*(Y)$$ is compatible with mixed Hodge structures.
\item \label{preservescup} The cup product $$H^k(X)\times H^{l}(X)\to H^{k+l}(X)$$ is compatible with mixed Hodge structures.
\item If $X$ is smooth $W_{j-1}H^j(X)$ is trivial.
\item If $X$ is smooth the {\em pure part} $PH^*(X):=\oplus_k W_kH^k(X)\subset H^*(X)$ is a subring.
\item If $X$ is smooth and $i:X\to Y$ is a smooth compactification of $X$, then 
$Im(i^*)=PH^*(X).$
\end{enumerate}
\end{theorem}

Using Deligne's \cite[8.3.8]{De-Hodge III} construction  of mixed Hodge structure on relative cohomology one can define \cite{danilov-khovanski} (for a general discussion of this cf. Note 11 on page 141 of \cite{fulton}) a well-behaved mixed Hodge structure on compactly supported cohomology $H^*_{c}(X):=H^*_c(X,\Q)$. Its basic properties are as follows (for proofs see \cite{peters-steenbrink}):
\begin{theorem}\label{mixedhodgec}\begin{enumerate} \item \label{forget} The forgetful map $$H_c^k(X)\to H^k(X)$$ is compatible with  mixed Hodge structures.
\item For a smooth connected $X$ we have 
Poincar\'e duality \beq\label{poincareduality}H^{k}(X)\times H^{2d-k}_c(X) \to H^{2d}_c(X)\cong \Q(-d)\eeq
is compatible with mixed Hodge structures, where $\Q(-d)$ is
the pure mixed Hodge structure on $\Q$ with weight $2d$ and Hodge
filtration $F^d=\Q$ and $F^{d+1}=0$.
\item In particular, for a smooth $X$, $W^{j+1}H_c^j(X)\cong H_c^j(X)$ .
\end{enumerate}
\end{theorem}

\begin{definition} Define the {\em Mixed Hodge numbers} by  
$$h^{p,q;j}(X):=\dim_\C \left(Gr^F_p Gr^W_{p+q}H^j(X)^\C\right),$$  and the {\em compactly supported Mixed Hodge numbers}  by
$$h_c^{p,q;j}(X):=\dim_\C \left(Gr^F_p  Gr^W_{p+q}H_c^j(X)^\C\right).$$ Form the {\em Mixed Hodge polynomial}: $$H(X;x,y,t):=\sum h^{p,q;j}(X) x^p y^q t^j,$$ the {\em compactly supported Mixed Hodge polynomial}: $$H_c(X;x,y,t):=\sum h_c^{p,q;j}(X) x^p y^q t^j,$$ and the {\em $E$-polynomial} of $X$:
$$E(X;x,y):=H_c(X;x,y,-1).$$
\end{definition}

\eqref{poincareduality} implies the following
\begin{corollary}\label{compactly} For a smooth connected $X$ of dimension $d$ we have $$H_c(X;x,y,t)=(xyt^2)^d H(X;1/x;1/y;1/t).$$
\end{corollary}
\begin{remark} \label{eulerchar} By definition $E(X;1,1)=H_c(X;1,1,-1)$ is the Euler
characteristic of $X$. 
\end{remark}

\begin{remark}  For our varieties $\M_n$ we will find in Corollary~\ref{hodgetriv}  that only  Hodge type $(p,p)$ can be non-trivial in the mixed Hodge structure, in other words $h^{p,q;j}=0$ unless $p=q$. Hence, $H(\M_n;x,y,t)$ only depends on $xy$ and $t$.  To simplify our notation we will denote by \beq \label{hqt}H(\M_n;q,t):=H(\M_n;\sqrt{q},\sqrt{q},t)\eeq and $$E(\M_n;q):=E(\M_n;\sqrt{q},\sqrt{q}).$$
\end{remark}

It is in fact the $E$-polynomial which could sometimes be calculated using arithmetic algebraic geometry. Here we explain a theorem of Katz (for details see the appendix). The  setup is the following. Let $X$ be a variety over $\C$. By a {\em spreading out} of $X$ we mean a separated scheme $\mathcal X$ over 
a finitely generated $\Z$-algebra with an embedding $\varphi:R\hookrightarrow \C$, such that the extension of scalars ${\mathcal X}_\varphi\cong X$. 
We say that $X$ has {\em polynomial count}\footnote{A similar property for smooth and proper schemes was studied in
\cite{bogaart-edixhoven}.} if there is a polynomial  
$P_{X}(t)\in \Z[t]$ and a spreading out $\mathcal X$  such that for every homomorphism $\phi: R\to \F_q$
to a finite field, the number of $\F_q$-points of the scheme 
${\mathcal X}_{\phi}$ is $$\#{\mathcal X}_{\phi}(\F_q)=P_{X}(q).$$ 
Then we have the following (cf Theorem~\ref{katz}.3)
 \begin{theorem}[Katz] Let $X$ be a variety
over $\C$. Assume $X$ has  polynomial count with count polynomial $P_X(t)\in \Z[t]$, then the $E$-polynomial of $X$ is 
given by: $$E(X;x,y)=P_X(xy).$$
\end{theorem}
\begin{remark} Informally this means that if we can count the number of solutions of 
the equations defining our variety over $\F_q$, and this number turns out to be  some universal polynomial evaluated at $q$, then this polynomial determines the $E$-polynomial of the variety. 
\end{remark}

In fact
it is enough for this to be true for all finite fields of all but finitely many
characteristics.  
We illustrate this in a  simple example. \begin{example} 
Fix a non-zero integer $m\in \Z$ and let
 $\mathcal X$ be the scheme over $\Z$ determined by the equation
 \begin{equation}
   \label{pol-count-example}
xy=m.    
 \end{equation}
 The extension of scalars ${\mathcal X}_\phi$ of ${\mathcal X}$ determined by a ring
 homomorphism $\phi: \Z \longrightarrow \F_q$ is given by the same
 equation \eqref{pol-count-example} now viewed over $\F_q$. It is easy to count solutions to \eqref{pol-count-example}. Let $p$ be
the characteristic of $\F_q$ (so that $q$ is a power of $p$). Then 
\begin{equation}
 \label{pol-count-example-1}
\#{\mathcal X}_\phi(\F_q)=
\left\{ \begin{array}{ll}
2q-1 & p \mid m\\
q-1 & \mbox{otherwise}
\end{array}
\right.
\end{equation}
Therefore ${\mathcal X}/\Z$ is fiberwise polynomial-count but not strongly
polynomial-count (for precise definitions see the appendix). This is not a contradiction to Theorem~\ref{strongly} of the
appendix; if we extend scalars to $\Z[\frac1m]$ then we eliminate the
primes dividing $m$ and find that in all cases $\#{\mathcal X}_\phi(\F_q)=q-1$, hence ${\mathcal X}$ has polynomial count.
In fact, $\mathcal{X}/\Z[\frac1m]$ is just isomorphic to $\Gm/\Z[\frac1m]$ .
\end{example}
\begin{example} To illustrate Katz's theorem further, 
we consider
the variety $X=\C^\times$. First we determine its mixed Hodge polynomial (cf. proof of Theorem 9.1.1 in \cite{De-Hodge III} ). The only question is to decide the Hodge numbers on
the one-dimensional $H^1(\C^*)$. Because $h^{0,1;1}=h^{1,0;1}$ 
and $h^{2,0;1}=h^{0,2;1}$ 
we must have $h^{1,1;1}(\C^\times)=1$ and the mixed Hodge polynomial is \beq \label{hcstar}H(\C^\times;x,y,t)=1+xyt.\eeq
Consequently, the compactly supported mixed Hodge polynomial is $$H_c(\C^\times;x,y,t)=t+xyt^2,$$ by Corollary~\ref{compactly} . Therefore 
the $E$-polynomial is $$E(\C^\times;x,y)=xy-1.$$ We can obtain the variety $\C^\times$  by extension of scalars $\Z\subset \C$ from the 
group scheme ${\mathbb G}_m$ over $\Z$ . The counting polynomial of this scheme is the polynomial $P_{\C^\times}(q)=q-1=\#{\mathbb G_m}(\F_q)$, which is consistent with Katz's theorem above. 
\end{example}

\subsection{Character varieties}
\label{character}
Here we define the character varieties and list their basic properties. 

Let $g\geq 0$, $n>0$ be integers. Let $\K$ be an algebraically closed field with $\zeta_n\in \K $ a primitive $n$-th root of unity. The
existence of such $\zeta_n$ is  equivalent to the condition \beq \label{char}{\rm char}(\K)\nmid n\eeq which we henceforth assume. Examples to bear in mind are $\K=\C$ and the algebraic closure of a finite field $\K=\overline{\F_q}$, where $q=p^r$ is a prime power, with $p\nmid n$. 

Denote by $I_n\in \GL{n,\K}$ the identity matrix,  $[A,B]:=ABA^{-1}B^{-1}\in \SL{n,\K}$, the commutator. 
The group $\GL{n,\K}$ acts by conjugation on $\GL{n,\K}^{2g}$: \bes \sigma: \GL{n,\K}\times \GL{n,\K}^{2g} & \to & \GL{n,\K}^{2g}\\ \nonumber (h,(A_1,B_1,\dots,A_g,B_g))&\mapsto& (h^{-1}A_1h,h^{-1}B_1h,\dots,h^{-1}A_gh,h^{-1}B_gh,).\ees As the center of $\GL{n,\K}$ acts trivially, this action induces an action \beq \label{actiongln} \bar{\sigma}:\PGL{n,\K}\times \GL{n,\K}^{2g} \to  \GL{n,\K}^{2g}\eeq  of $\PGL{n,\K}$.   Let $\mu: \GL{n,\K}^{2g}\to \SL{n,\K}$  be
given by
 $$\mu(A_1,B_1,\dots,A_g,B_g):= [A_1,B_1] \dots [A_g,B_g].$$  
 We define  \beq \label{undef} \Un:=\mu^{-1}(\zeta_nI_n).\eeq 
 Clearly the $\PGL{n,\K}$-action \eqref{actiongln} will leave the affine variety $\Un$ invariant. Thus we have a $\PGL{n,\K}$ action on $\Un$: \beq \label{action}\bar{\sigma}: \PGL{n,\K}\times \Un\to \Un.\eeq The  categorical quotient  \beq \label{categorical} \pi_n: \Un\to \M_n\eeq 
  exists by \cite[Theorem1.1]{mumford-etal} (cf. also \cite[\S 3]{newstead})  in the sense of geometric invariant theory \cite{mumford-etal}.   Explicitly we have  $$\M_n={\rm Spec}(\K[\Un]^{\PGL{n,\K}})$$  and $\pi_n$ is induced by the obvious embedding of $\K[\Un]^{\PGL{n,\K}}\subset \K[\Un]$. We call $\M_n$ a {\em twisted  $\GL{n,\K}$-character variety} of a closed Riemann surface of genus $g$.  We will use the notation $\M_n/\K$ for the variety $\M_n$, when we want to emphasize the ground field $\K$.

\begin{example}\label{g=0-char-var}
When $g=0$  $\M_n$ is clearly empty, unless $n=1$ when it is a point.\end{example}

\begin{example} When $n=1$, $\SL{1,\K}$ and $\PGL{1,\K}$ are  trivial, and so $\M_1=\GL{1,\K}^{2g}\cong (\K^\times)^{2g}$ is a torus.  
The mixed Hodge polynomial of $\M_n/\C$ then is  \beq \label{htorus} H(\M_1/\C;x,y,t)=(1+xyt)^{2g}.\eeq
by Theorem~\ref{mixedhodgeproperties}.\ref{kunneth} and 
\eqref{hcstar}.
\end{example}

\begin{remark} \label{scheme} It will be important for us to have a spreading out ${\cal X}_n/R$ of the variety 
$\M_n/\C$ over a finitely generated $\Z$-algebra $R$.   Clearly $\Un$ can be defined to be an affine scheme  over $R:=\Z[\zeta_n]$. 
Using Seshadri's extension
of geometric invariant theory quotients for schemes \cite{seshadri},
we can take the categorical quotient by the conjugation action of the  reductive group scheme $\PGL{n,R}$. Explicitly, ${\cal X}_n=\rm{Spec}\left(R[\Un]^{\PGL{n,R}}\right)$. 
As the embedding $\phi: R\to \C$ is a flat morphism, \cite[Lemma 2]{seshadri} implies that 
$$R[\Un]^{\PGL{n,R}}\otimes_R \C=\C[\Un]^{\PGL{n,\C}},$$ thus $\X_n$ is the required 
spreading out of $\M_n/\C$. Roughly speaking the scheme ${\cal X}_n/R$ will be our bridge
between the varieties $\M_n/\C$ and $\M_n/\overline{\F_q}$. 
\end{remark} 

We have the following immediate 
\begin{corollary} The mixed Hodge polynomial $H(\M_n/\C;x,y,t)$ 
does not depend on the choice of the primitive $n$-th root of unity $\zeta_n\in \C$. 
\end{corollary}
\begin{proof}  As $\M_n/\C$ can  be obtained by base change from ${\mathcal X}_n/R$ with $\phi:R\to \C$, we see that
the Galois conjugate $\M_n^\sigma$ for any field automorphism $\sigma:\C\to \C$ can be obtained from the same scheme ${\mathcal X}_n/R$  by extension of scalars $\sigma\phi:R\to \C$.  Now  Theorem~\ref{mixedhodgeproperties}.\ref{galois} implies the corollary. 
\end{proof}

\begin{theorem} \label{smooth} The variety $\M_n$ is non-singular.  \end{theorem}
\begin{proof}  Because of Example~\ref{g=0-char-var} we can assume $g>0$ for the rest of the proof.

We first prove that the affine subvariety $\Un\subset (\GL{n,\K})^{2g}$ is non-singular. 
 By  definition it is enough to show
  that at  a solution $s=(A_1,B_1,\dots, A_g,B_g)$ of the
 equation   \beq \label{abeqn} [A_1,B_1] \dots [A_g,B_g] = \zeta_n I_n\eeq the derivative of $\mu$ on the tangent spaces  $$d\mu_s:T_s(\GL{n,\K}^{2g})\to T_{\zeta_n I_n} \SL{n,\K}$$
 is
 surjective. So take $(X_1,Y_1,\dots,X_g,Y_g)\in
 T_s(\GL{n,\K}^{2g})\cong {\mathfrak g\mathfrak l}(n,\K)^{2g}$.  Then
 differentiate $\mu$ to get: \bes d\mu_s(X_1,Y_1,\dots,X_g,Y_g)=&&
 \sum_{i=1}^g [A_1,B_1]\dots[A_{i-1},B_{i-1}]
 X_iB_iA_i^{-1}B_i^{-1}[A_{i+1},B_{i+1}]\dots [A_g,B_g] \\ & + &
 \sum_{i=1}^g [A_1,B_1]\dots[A_{i-1},B_{i-1}]
 A_iY_iA_i^{-1}B_i^{-1}[A_{i+1},B_{i+1}]\dots [A_g,B_g] \\&- &
 \sum_{i=1}^g [A_1,B_1]\dots[A_{i-1},B_{i-1}]
 A_iB_iA_i^{-1}X_iA_i^{-1}B_i^{-1}[A_{i+1},B_{i+1}]\dots [A_g,B_g]
 \\& -& \sum_{i=1}^g [A_1,B_1]\dots[A_{i-1},B_{i-1}] A_iB_iA_i^{-1}B_i^{-1}Y_i
 B_i^{-1}[A_{i+1},B_{i+1}]\dots [A_g,B_g],  \ees where we used the product rule for matrix valued functions, in particular that  $d\nu_A(X)=-A^{-1} X A^{-1}$ for the derivative of the function $\nu:\GL{n,\K}\to \GL{n,\K}$ defined by $\nu(A)=A^{-1}$ at $A\in \GL{n,\K}$ and $X\in T_{A}\GL{n,\K}\cong {\mathfrak g \mathfrak l} (n,\K)$.  Using (\ref{abeqn})
 for each of the four terms, we get: \beq
 d\mu_s(X_1,Y_1,\dots,X_g,Y_g)=\sum_{i=1}^g\left(f_i({X}_i)+g_i({Y}_i)\right),\label{mus} \eeq where
 we define  linear maps $f_i:{\mathfrak g \mathfrak l}(n,\K) \to {\mathfrak
   s \mathfrak l}(n,\K)$ and $g_i:{\mathfrak g \mathfrak l}(n,\K) \to
 {\mathfrak s \mathfrak l}(n,\K)$ by $$f_i(X)=\zeta_n
 [A_1,B_1]\dots[A_{i-1},B_{i-1}]\left( XA_i^{-1} - A_iB_i A_i^{-1}X
   B_i^{-1} A_i^{-1}\right)\left([A_1,B_1]\dots[A_{i-1},B_{i-1}]
 \right)^{-1}$$
 and $$g_i(X)=\zeta_n
 [A_1,B_1]\dots[A_{i-1},B_{i-1}]\left( A_iY_iB_i^{-1}A_i^{-1}-A_iB_iA_i^{-1}B_i^{-1}Y_i A_i B_i^{-1}A_i^{-1}\right)
 \left([A_1,B_1]\dots[A_{i-1},B_{i-1}] \right)^{-1}.$$

  Assume that $Z\in {\mathfrak s\mathfrak l}(n,\K)$ such that \beq\label{surjective}
  \tr(Zd\mu_s(X_1,Y_1,\dots,X_g,Y_g))=0\eeq for all $X_i$ and $Y_i$.
  By \eqref{mus} this is equivalent to 
  \bes \tr\left(Zf_i(X_i)\right)=\tr\left(Zg_i(X_i)\right)=0\ees for all $i$ and $X_i\in \GL{n,\K}$. We show by induction on $i$ that this implies that $Z$ commutes with $A_i$ and $B_i$.  Assume we have already proved
 this for $j<i$ and calculate $$0=  \tr\left(Zf_i(X_i)\right)=\tr\left(
 \left(A_i^{-1}Z-B_i^{-1} A_i^{-1}ZA_iB_iA_i^{-1}\right)X_i\right)$$ for all $X_i$,  thus $Z$ commutes with $A_iB_iA_i^{-1}$. Similarly we have $$0=  \tr\left(Zg_i(X_i)\right)=\tr\left(
 \left(B_i^{-1}A_i^{-1}ZA_i-A_iB^{-1}_i A_i^{-1}ZA_iB_iA_i^{-1}B_i^{-1}\right)X_i\right),$$ which implies that $Z$ commutes with $A_iB_iA_iB_i^{-1}A_i^{-1}$. Thus $Z$ commutes with $A_i$ and $B_i$. 
 The next lemma proves that this implies that $Z$ has to be central. 
 Because $Z$ was traceless, we also get $Z=0$ by \eqref{char}. 
     Thus there is no non-zero $Z$ such  that
  \eqref{surjective}  holds for all $X_i$ and $Y_i$.  Again because of \eqref{char} this implies that $d\mu$ is surjective at any solution $s$ of (\ref{abeqn}). Thus $\Un$ is non-singular.

\begin{lemma}\label{irreducible} 
Suppose $Z\in {\mathfrak g\mathfrak l}(n,\K)$ commutes with each of the $2g$
matrices $A_1, B_1,
\dots, A_g,B_g\in \GL{n,\K}$, which solve (\ref{abeqn}). Then $Z$ is 
central.  
\end{lemma} \begin{proof} If $Z$ is not central and $\lambda$ is an eigenvalue
then $E=\ker(Z-\lambda I_n)$ is a proper subspace of $\K^n$. 
Because $A_i$ and $B_i$ all commute with $Z$, they preserve
$E$. 
 Let
 $\tilde{A}_i=A_i|_E$ and $\tilde{B}_i=B_i|_E$.  Then restricting
 \eqref{abeqn} to $E$ we get \bes
 [\tilde{A}_1,\tilde{B}_1] \dots [\tilde{A}_g,\tilde{B}_g] = \zeta_n
 I_E. \ees As the determinant of a commutator is $1$,  the determinant of the right hand side has to be $1$. But this implies
 $\zeta_n^{\dim E}=1$, which is a contradiction as $\zeta_n$ is  a primitive $n$-th root of unity. The lemma follows. 
 \end{proof}

   This lemma also proves that if $g\in \GL{n,\K}$ is not 
  central then it acts {\em set-theoretically freely} on the
  solution space of (\ref{abeqn}).  We can also deduce the following
  more general
  \begin{corollary} The action $\bar{\sigma}$ of 
  $\PGL{n,\K}$ on $\Un$, defined in \eqref{action}, is {\em scheme-theoretically free} (\cite[Definition 0.8 (iv)]{mumford-etal}).
  \end{corollary}
  \begin{proof} The statement says that the map $\Psi:=(\bar{\sigma},p_2):\PGL{n,\K}\times \Un \to \Un\times \Un$ is a closed immersion. We prove it by an argument similar to the proof of \cite[Lemma 6.5]{reineke}.

  To prove  this consider the map $ \Phi:\Un\times \Un \to \Hom_\K({\mathfrak g \mathfrak l}(n,\K),{\mathfrak g \mathfrak l}(n,\K)^{2g}),$ defined by $$\Phi\left((A_1,B_1,\dots,A_g,B_g),(\tilde{A}_1,\tilde{B}_1,\dots,\tilde{A}_g,\tilde{B}_g)\right)(h)=\left(hA_1-\tilde{A}_1h,\dots,hB_g-\tilde{B}_gh \right).$$ We show that $\Phi(x,y)$ 
  has a non-trivial kernel if and only if $(x,y)\in \Un\times \Un$ is in the image of $\Psi$, i.e., there is an $\bar{h}\in \PGL{n,\K}$ such that $y=\bar{h}x$. The if part is clear.
  For the other direction assume that $\Phi(x,y)$ has a non-trivial
  kernel, i.e.,  $0\neq h\in {\mathfrak g\mathfrak l}(n,\K)$, such that $\Phi(x,y)(h)=0$.  Then the matrices $A_1,B_1,\dots,A_g,B_g$,
  which solve \eqref{abeqn}, 
  will leave $\ker (h)$ invariant. As in the proof of
  Lemma~\ref{irreducible} this implies that $\ker(h)$ is trivial i.e., $h$ is invertible. So indeed $\ker(\Phi(x,y))\neq 0$ implies that there exists $\bar{h}\in \PGL{n,\K}$ such that $y=\bar{h}x$.  We also get that in this case $\dim(\ker{\Phi(x,y)})=1$ as the action $\bar{\sigma}$ is set-theoretically
  free. 

Now we fix a basis for $\gl(n,\K)$, and take the closed subscheme 
  given by the vanishing of all $n^2\times n^2$ minors in the entries
  of the matrices $\Phi(x,y)\in \Hom_\K({\mathfrak g \mathfrak l}(n,\K),{\mathfrak g \mathfrak l}(n,\K)^{2g})$. This shows that the image
  of $\Psi$ is a closed subscheme $Z$ of $\Un\times \Un$. 

  Moreover on the Zariski open subscheme of $Z$ where a given
  $(n^2-1)\times (n^2-1)$ minor of $\Psi(x,y)$ is non-zero, we can solve
  algebraically for the unique $\bar{h}\in\PGL{n,\K}$ such that $y=\bar{h}x$,
  giving us locally an inverse $Z\to \Un\times \Un$ to $\Psi$; showing
  that $\Psi$ is an isomorphism onto its image. The Corollary follows.
  \end{proof}
           In particular $\Psi$ is a closed map. Consequently the action is closed so \cite[Amplification 1.3, Proposition 0.9]{mumford-etal} imply
  \begin{corollary} \label{principal} The categorical quotient $(\M_n,\pi_n)$ is a geometric quotient and $\pi_n$ in \eqref{categorical} is a $\PGL{n,\K}$-principal bundle, in particular $\pi_n$ is flat. 
  \end{corollary}   
 Because the geometrical fibres of the flat morphism $\pi_n$ are non-singular (they
 are all isomorphic to $\PGL{n,\K}$) $\pi_n$ is a smooth morphism by
 \cite[Theorem III.10.3']{mumford-red}. By 
 \cite[Corollary 17.16.3]{grothendieck2} a smooth surjective morphism  locally has an \'etale section, so \'etale-locally the principal bundle $\pi_n:\Un\to \M_n$ is trivial. As
 $\Un$ is non-singular, we get  that $\M_n$ is also non-singular.
 \end{proof}

We will see in Corollary~\ref{connected} that our varieties $\Un$ and $\M_n$ are connected. Here we can determine their dimension.

\begin{corollary} \label{dimension} For $g>0$ the dimension of (each connected component of) $\M_n$ is $d_n:=n^2(2g-2)+2$.
\end{corollary}

\begin{proof} From the previous proof we see that the dimension of
(each connected component of) $\Un$ is $$\dim(\GL{n,\K}^{2g})-\dim( \SL{n,\K})=2gn^2-(n^2-1).$$ Because $\pi_n$ is flat we have to
subtract $\dim(\PGL{n,\K})=n^2-1$ from this to get the dimension of $\M_n$ proving the claim.
\end{proof}

\begin{definition} The torus $(\K^\times)^{2g}$ acts on $\Un\subset \GL{n,\K}^{2g}$ by the following formula: \beq \label{tauaction} \tau:  (\K^\times)^{2g}\times \GL{n,\K}^{2g} & \to & \GL{n,\K}^{2g}\\ \nonumber ((\lambda_1,\dots,\lambda_{2g} ),(A_1,B_1,\dots,A_g,B_g))&\mapsto& (\lambda_1 A_1,\lambda_2B_1,\dots,\lambda_{2g-1}A_g,\lambda_{2g}B_g).\eeq This action  commutes with the action $\bar{\sigma}$.  Thus $(\K^\times)^{2g}$ acts on  $\M_n$. We call the  categorical quotient $$\tM_n:=\M_n/\!/(\K^\times)^{2g}\cong \Un/\!/ \left(\PGL{n,\K}\times (\K^\times)^{2g}\right)$$  the {twisted \em $\PGL{n,\K}$-character
variety} of the genus $g$ Riemann surface $\Sigma$. 
\end{definition}
\begin{remark} $\tM_n$ could be considered as a component of the variety of homomorphisms
of $\pi_1(\Sigma)$ into $\PGL{n,\K}$ modulo conjugation, this motivates its name. 
\end{remark}

\begin{theorem} \label{tmn} The variety $\tM_n$ is an orbifold.  
Each connected component  of $\tM_n$ has   dimension  $\tilde{d}_n=(n^2-1)(2g-2)$.  Moreover when $\K=\C$ its cohomology satisfies  $$H^*(\M_n/\C)= H^*(\tM_n/\C)\otimes H^*(\M_1/\C),$$
 and the mixed Hodge polynomial satisfies:\beq\label{htm} H(\M_n/\C;x,y,t)=H(\tM_n/\C;x,y,t)(1+xyt)^{2g}. \eeq 
\end{theorem}
\begin{proof} Let $\mu^\prime: \SL{n,\K}^{2g}\to \SL{n,\K}$  be
given by
 $$\mu^\prime(A_1,B_1,\dots,A_g,B_g):= [A_1,B_1] \dots [A_g,B_g].$$ 
 We define  \beq \Unp:=(\mu^\prime)^{-1}(\zeta_nI_n).\eeq 
 The $\PGL{n,\K}$-action \eqref{actiongln} on $\SL{n,\K}^{2g}\subset \GL{n,\K}^{2g}$ will leave the affine variety $\Unp$ invariant. We have  the  categorical quotient  \beq \label{categoricalprime} \pi^\prime_n: \Unp\to \Mp_n,\eeq defining the {\em twisted $\SL{n,\K}$-character variety} $\Mp_n$. Exactly as in the $\GL{n,\K}$ case we can
 argue that $\Unp$ and $\Mp_n$ are non-singular, $\pi^\prime_n$ is
 a $\PGL{n,\K}$-principal bundle and the components of $\Mp_n$ have dimension 
$$\dim (\SL{n,\K}^{2g})-\dim \SL{n,\K}-\dim \PGL{n,\K}= (n^2-1)(2g-2).$$

We denote by ${\bbmu}_n$  the group scheme of $n$-th roots of unity.  
$\bbmu_n^{2g}\subset (\K^\times)^{2g}$ acts on $\Unp\subset \Un\subset \GL{n,\K}^{2g}$ induced from the action \eqref{tauaction}.  It commutes with  the action $\bar{\sigma}$ in \eqref{action} and so $\bbmu_n^{2g}$ also acts on $\Mp_n$. 
Note that the map  $\SL{n,\K}\times \K^\times \to \GL{n,\K}$ given by multiplication is the categorical quotient of the action of the subgroup scheme $\bbmu_n=\{(\zeta_n^d I_n,\zeta_n^{-d}),d=1,\dots,n\}\subset \SL{n,\K}\times\K^\times $ on $\SL{n,\K}\times \K^\times$.  Therefore we can identify $\Un=(\Unp \times (\K^\times)^{2g} )/\!/\bbmu_n^{2g}$ and taking quotients we have
\beq \label{quotientm} \M_n\cong (\M_n^\prime \times (\K^\times)^{2g})/\!/\bbmu_n^{2g}\eeq

In particular
we see that the categorical quotient \beq \label{orbifold}\tM_n=\M_n/\!/(\K^\times)^{2g}\cong \left(\Mp_n\times (\K^\times)^{2g}\right)/\!/ \left( \bbmu_n^{2g} \times  (\K^\times)^{2g} \right)\cong \M^\prime_n/\!/\bbmu_n^{2g}\eeq is an orbifold
of dimension $(n^2-1)(2g-2)$.

When we take cohomologies  in \eqref{quotientm} we get:
$$H^*(\M_n/\C)=(H^*(\M_n^\prime/\C \times \M_1/\C))^{\bbmu_n^{2g}}=H^*(\M_n^\prime/\C)^{\bbmu_n^{2g}}\otimes H^*(\M_1/\C)=H^*(\tM_n/\C)\otimes
H^*(\M_1/\C),$$ by \eqref{orbifold}, the K\"unneth theorem, the fact
that $\bbmu_n^{2g}$ acts trivially on $H^*(\M_1/\C)$ and  the observation
of Grothendieck \cite{grothendieck} that the rational cohomology
of a quotient of a smooth variety by a finite group like  
\eqref{orbifold} and \eqref{quotientm} is the invariant part of the cohomology of the space.  The Theorem follows.
\end{proof}

For $g=1$ we now
determine our varieties $\M_n$ and $\tM_n$  explicitly.

\begin{lemma}
\label{g=1-lemma-1}
Let $H\subset \GL{n,\K}$ be the subgroup generated by $A,B \in
\GL{n,\K}$ satisfying
\begin{equation}
 \label{g=1-eqn}
[A,B]=\zeta_nI_n,  
\end{equation}
where $I_n$ is the identity matrix. Then the corresponding action
of $H$ on $\K^n$ is irreducible.
\end{lemma}
\begin{proof}
The proof is the same as in the proof of Lemma~\ref{irreducible}.
\end{proof}
\begin{lemma}
\label{g=1-lemma-2}
 With the notation of the previous lemma  \ref{g=1-lemma-1} we have  
\beq \label{g=1-eqn-3}
A^n=\alpha I_n, \quad B^n=\beta I_n, \qquad \alpha,\beta \in \K^\times.
\eeq
\end{lemma}
\begin{proof}
{}From \eqref{g=1-eqn} we easily deduce that
\begin{equation}
 \label{g=-1-eqn-1}
 A^jB^k=\zeta_n^{jk}B^kA^j, \qquad j,k \in \Z.
\end{equation}
In particular, $A^n$ and $B^n$ are in the center of $H$ and our claim
follows from Schur's lemma and Lemma~\ref{g=1-lemma-1}.
\end{proof}
\begin{lemma}
\label{g=1-lemma-3}
 There exists a unique solution, up to conjugation, to the equations
 \begin{equation}
   \label{g=1-eqn-2}
A^n=B^n=I_n, \qquad [A,B]=\zeta_nI_n, \qquad A,B \in \GL{n,\K},    
 \end{equation}
where $I_n$ is the identity matrix.
\end{lemma}
\begin{remark}
 The group $H$ generated by the matrices in the hypothesis of the
 lemma is a finite Heisenberg group. The lemma is a version of the
 Stone--von-Neumann theorem on the uniqueness of the Heisenberg
 representation.
\end{remark}
\begin{proof}
Let $v\in \K^n$ be an eigenvector of $A$, say
 $Av=\zeta v$ with $\zeta^n=1$. Then $ABv=\zeta_nBAv=\zeta\zeta_n Bv$
 and $Bv$ is also an eigenvector of $A$. Repeating the process we see
 that $B^kv$ is an eigenvector of $A$ for all $k\in \Z$.

 Since the action of $H$ is irreducible by Lemma \ref{g=1-lemma-1} we
 must have that $v,Bv,\ldots,B^{n-1}v$ is a basis of $\K^n$ (their
 span is clearly stable under $H$). Replacing $v$ by an appropriate
 vector $B^kv$ if necessary we may assume that $\zeta=1$. Hence in
 this basis $A$ is the diagonal matrix with entries
 $1,\zeta_n,\zeta_n^2,\ldots,\zeta_n^{n-1}$ along the diagonal and
 $B$ is the permutation matrix corresponding to the $n$-cycle
 $(12\cdots n)$. It is easy to verify that these particular matrices
 are indeed solutions to the equations \eqref{g=1-eqn-2} and we have
 shown  all  pairs of matrices satisfying \eqref{g=1-eqn-2} are
 conjugate to these proving our claim.
\end{proof}

\begin{theorem}
\label{g=1-char-var}
The orbits of the action of $\GL{n,\K}$ acting on the solutions to 
\begin{equation*} 
[A,B]=\zeta_nI_n, \qquad A,B\in \GL{n,\K}
\end{equation*}
by conjugation are in bijection with $\K^\times\times\K^\times$ via
$(A,B)\mapsto (\alpha,\beta)$ where $A^n=\alpha I_n,B^n=\beta I_n$. Consequently $\M_n\cong \K^\times\times\K^\times$ and $\tM_n$ is a point when $g=1$.
\end{theorem}

\begin{proof} Consider the action of $(\K^\times)^2$ on $\M_n$ induced by \eqref{tauaction}. Let $x_0\in\M_n$ be the point corresponding to the unique $\PGL{n,\K}$ orbit in $\Un$ of pairs of matrices $(A,B)$ solving \eqref{g=1-eqn-2}. Then for such a pair of matrices $\tau(\lambda_1,\lambda_2)(A,B)=(\lambda_1A,\lambda_2B)$ will give a solution of $\eqref{g=1-eqn}$, such that 
\eqref{g=1-eqn-3} will hold with $\alpha=\lambda_1^n$ and $\beta=\lambda_2^n$.  Because of the uniqueness of $x_0$, we see that the action of $(\K^\times)^2$ on $\M_n$ is transitive. Therefore $\tM_n$ is a point.  The stabilizer of $x_0$  is 
$\bbmu_n^2\subset (\K^\times)^2$. It follows
that $\M_n\cong (\K^\times)^2/\!/\bbmu^{2}_n\cong \K^\times\times \K^\times$, and the isomorphism $\M_n\to \K^\times\times \K^\times$ is given by the map in the theorem.  \end{proof}

\subsection{Counting solutions to equations in finite groups}
\label{count}

We collect in this section various known results about counting
solutions to equations in finite groups that we will need. These and
similar results have appeared in the literature in many places see for
example \cite{Se}, \cite{freed-quinn}, \cite{mednykh}. Interestingly,
the first application in Frobenius's \cite{frobenius} of 1896, where
he introduced characters of finite groups, were formulas of similar
type.  (Those that relate to a Riemann sphere with punctures.)

These counting formulas arise naturally, when considering Fourier transform
on finite groups. This point of view will be discussed in \cite{hausel5}, where it is shown that the counting formulas below and the one 
in \cite{hausel4} have the same origin.   

Let $G$ be a finite group. For a function
$$
f:\quad G\longrightarrow \C
$$
we define
$$
\int_Gf(x)\;dx := \frac1{|G|} \sum_{x\in G} f(x)
$$
Given a word $w \in F_n$, where $F_n=\langle X_1,\ldots,X_n\rangle$
is the free group in generators $X_1,\ldots, X_n$, and a function $f$
on $G$ as above we define
\begin{equation}
 \label{char-sum}
 \{f,w\}:=\int_{G^n}f(w(x_1,\ldots,x_n))\;dx_1\ldots dx_n,
\end{equation}
where $w(x_1,\ldots,x_n)$ is a shorthand for $\phi(w)\in G$ with
$\phi:F_n\longrightarrow G$ the homomorphism mapping each $X_i$ to
$x_i$.

\begin{lemma}
 \label{char-sum-lemma}
With the above notation we have for any $z\in G$ and $\chi$ any
irreducible character of $G$
\begin{equation}
 \label{char-sum-z}
 \int_{G^n}\chi(w(x_1,\ldots,x_n)z)\;dx_1\ldots
 dx_n=\{\chi,w\}\;\frac{\chi(z)}{\chi(1)} 
\end{equation}
\end{lemma}
\begin{proof}
Consider the linear endomorphism of the vector space $V$ of a
representation $\rho$ of $G$ with character $\chi$
$$
W:=\frac1{|G|^n}\sum_{(x_1,\ldots,x_n) \in G^n}\rho(w(x_1,\ldots,x_n)).
$$
Changing each $x_i$ in  the sum defining $W$ by $zx_iz^{-1}$ for some
$z\in G$ does not change the sum. On the other hand,
$w(zx_1z^{-1},\ldots,zx_nz^{-1})=zw(x_1,\ldots,x_n)z^{-1}$, hence
$W=\rho(z)W\rho(z)^{-1}$. In other words, $W$ is $G$-linear. 
By Schur's lemma $W$ is a scalar; taking traces we find that
$$
W=\frac{\{\chi,w\}}{\chi(1)}\;\id_V.
$$
Multiplying both sides by $\rho(z)$ on the right and taking traces again we
obtain (\ref{char-sum-z}). 
\end{proof}

\begin{proposition}
\label{sol-num-prop}
With the above notation let  $N(z)$ be the number of solutions to
\begin{equation}
 \label{sol-num}
 w(x_1\cdots x_n)z=1,\qquad (x_1,\ldots,x_n)\in G^n,
\end{equation}
then
\begin{equation}
 \label{sol-num-fmla}
 N(z)=|G|^{n-1}\sum_\chi {\{\chi,w\}}\;\chi(z)
\end{equation}
where the sum is over all irreducible characters of $G$.
\end{proposition}
\begin{proof}
Write the delta function on $G$
$$
\delta(x)=
\left\{ \begin{array}{ll}
1 & x=1\\
0 & \mbox{otherwise}
\end{array}
\right.
$$
as a linear combination of the irreducible characters of $G$
\begin{equation}
\label{delta}
\delta=\sum_\chi c_\chi \chi,
\end{equation}
where
$$
c_\chi=(\chi,\delta)=\int_G\chi(x)\delta(x)\;dx=\frac{\chi(1)}{|G|}.
$$
On the other hand,
$$
N(z)=|G|^n\int_{G^n}\delta(w(x_1,\ldots, x_n)z)\;dx_1\ldots dx_n,
$$
which  combined with (\ref{delta}) and (\ref{char-sum-z})  yields our claim. 
\end{proof}

Consider now words $w_i \in \langle X_1^{(i)},\ldots,
X_{n_i}^{(i)}\rangle$  in disjoint set of variables for $i=1,\ldots,k$
and let $w=w_1\cdots w_k \in \langle
X_1^{(1)},\ldots,X_{n_1}^{(1)},X_1^{(2)},\ldots,X_{n_2}^{(2)},\ldots\rangle$.
From  lemma \ref{char-sum-lemma} it follows by induction that
\begin{equation}
 \label{char-sum-multi}
 \{\chi,w\}=\frac{\{\chi,w_1\}\cdots \{\chi,w_k\}}{\chi(1)^{k-1}}.
\end{equation}

As an application, consider $w=[x,y]=xyx^{-1}y^{-1}$. It is not hard to
verify that
$$
\{\chi,w\}=\frac1{\chi(1)}.
$$
Indeed, consider the linear endomorphism of the vector space $V$ of a
representation $\rho$ of $G$ with character $\chi$
$$
W:=\frac1{|G|}\sum_{x \in G}\rho(xyx^{-1}).
$$
By changing variables in the sum we see that $W\rho(z)=\rho(z)W$ for
all $z\in G$. Hence, by Schur's lemma $W$ is a scalar; taking traces
we find that 
$$
W=\frac{\chi(y)}{\chi(1)}\;\id_V.
$$
Now we note that $\{\chi,w\}$ is the trace of 
$$
\frac1{|G|}\sum_{y\in G} Wy^{-1}= \frac1{\chi(1)|G|}\sum_y \chi(y)\rho(y)^{-1}
$$
and our claim follows.
We conclude from  proposition \ref{sol-num-prop} that for any\footnote{In this paper we avoid the use of the notation $\N$ as the  notion of natural numbers is different for the two authors. Instead we use the notation  $\Z_{\geq 0}$ and $\Z_{>0}$ respectively.} $g\in \Z_{\geq 0} $
\begin{equation}
 \label{sol-num-comm}
 \#\{(x_1,y_1,\ldots,x_g,y_g) \in G^{2g} \;|\; [x_1,y_1]\cdots[x_g,y_g]z=1\}
 = \sum_\chi\left(\frac{|G|}{\chi(1)}\right)^{2g-1}\chi(z).
\end{equation}
\begin{remark}
\label{hom-divis}
 For $z=1$ the quantity in \eqref{sol-num-comm} equals
 $\#\Hom(\Gamma_g,G)$ where $\Gamma_g$ is the fundamental group of a
 genus $g$ Riemann surface.  Hence we have
\begin{equation}
 \label{hom-count}
 \#\Hom(\Gamma_g,G)=|G|\sum_\chi\left(\frac{|G|}{\chi(1)}\right)^{2g-2},
\end{equation}
which, in particular, implies the remarkable fact that $|G|$
always divides $\#\Hom(\Gamma_g,G)$ for $g>0$.
\end{remark}

\subsection{Partitions}
\label{partitions}
We collect in this section some  notation and concepts on partitions
that we will need later. The main reference is Macdonald's book 
\cite{Mc}.

Let $\p_m$ be the set of all partitions $\lambda$ of a non-negative
integer $m=|\lambda|$ (where for $m=0$ we only have the zero partition
$\{0\}$) and $\p=\bigcup_m\p_m$.  We write a
partition of $n$ as $\lambda =(\lambda_{1}\geq \lambda_{2}\geq \cdots \geq \lambda
_{l}>0)$, so that $\sum \lambda_i=n$.  The
{\it Ferrers diagram} $d(\lambda)$ 
of $\lambda $ is the set of lattice points
\begin{equation}\label{e:d(lambda)}\{(i,j)\in {\mathbb Z_{\leq 0}} \times
 {\Z_{\geq 0}} : j <\lambda_{-i+1} \}.
\end{equation}
The {\it arm length} $a(z)$ and
{\it leg length} $l(z)$ of a point  $z\in d(\lambda )$  (sometimes called a {\em box}) denote the number of points
strictly to the right of $z$ and below $z$, respectively, as indicated
in this example: 

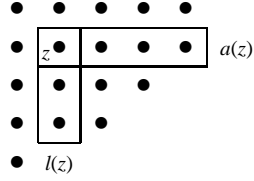
\begin{eqnarray*}\begin{array}[c]{cccccc}
\bullet &      \bullet &       \bullet &
\bullet &       \bullet  \\

\cline{2-5}
\bullet &       \multicolumn{1}{|c|}{\llap{${}_{z}$}\bullet } & \bullet &
   \bullet&     \multicolumn{1}{c|}{\bullet }& {\scriptstyle a(z)} \\
\cline{2-5}
\bullet &       \multicolumn{1}{|c|}{\bullet }& \bullet &       \bullet \\

\bullet &       \multicolumn{1}{|c|}{\bullet }& \bullet \\
\cline{2-2}
\bullet &       \hbox to 0pt{\hss $\scriptstyle l(z)$\hss }\\
\end{array}
\end{eqnarray*}
where $\lambda=(5,5,4,3,1)$, $z=(-1,1)$, $a(z)=3$ and $l(z)=2$. The {\em hook length} then is defined as \beq \label{hooklength} h(z)=l(z)+a(z)+1.\eeq

Given two partitions $\lambda, \mu\in
\p$ we define
\begin{equation}
 \label{part-pairing}
\langle \lambda,\mu \rangle  = \sum_{j\geq 1} \lambda_j'\mu_j',
\end{equation}
where $\lambda' =(\lambda'_1,\lambda'_2,\ldots)$ and $\mu'
=(\mu'_1,\mu'_2,\ldots) $ are the dual partitions.

For $\lambda =(\lambda_1,\lambda_2,\ldots) \in \p$ we define 
\begin{equation}
 \label{def-n-lambda}
 n(\lambda):=\sum_{i\geq 1}(i-1)\lambda_i
\end{equation}
then
\begin{equation}
 \label{part-pairing-1}
 \langle \lambda,\lambda \rangle = 2n(\lambda)+|\lambda|.
\end{equation}
Let the {\it hook polynomial} be
(see \cite[p. 152]{Mc})
\begin{equation}
 \label{def-dim-pol}
\tilde H_\lambda(q)=  
\prod(q^h-1),
\end{equation}
where for the product is taken for the set of boxes $d(\lambda)$ in the Ferrers diagram of $\lambda$ and we let $h=h(z)$
denote the hook length of a box $z\in d(\lambda)$ as defined in \eqref{hooklength}.

It will be convenient for us to work with Laurent polynomials in
$q^{\tfrac12}$ and scale the hook polynomial by an appropriate power
of $q$. Concretely, we let
\begin{equation}
 \label{hook-polynomial}
 \h_\lambda(q):=q^{-\tfrac12\langle \lambda,\lambda \rangle}
\prod(1-q^h),
\end{equation}

Hence with this normalization we have
\begin{equation}
 \label{hook-duality}
 \h_\lambda(q^{-1})=(-1)^{|\lambda|}\h_{\lambda'}(q)
\end{equation}
since 
\begin{equation}
 \label{sum-hook-lengths}
 \sum h= n(\lambda)+n(\lambda')+|\lambda|.
\end{equation}

For a  non-negative integer $g$ and a partition $\lambda$ we define
\begin{equation}
 \label{H-defn}
\h_\lambda(z,w):=\prod
\frac{\left(z^{2a+1}-w^{2l+1}\right)^{2g}}  
{(z^{2a+2}-w^{2l})(z^{2a}-w^{2l+2})}
\end{equation}
a rational function in $z,w$, where the product runs over the boxes in
$d(\lambda)$ with $a$ and $l$ the corresponding arm and leg length.
(Typically $g$ will be fixed and hence there is no need to indicate it
in the notation. Also the context should make clear whether
$\h_\lambda$ represents the one and two variable version.)

We note the following  easily checked properties  of $\h_\lambda$.
\begin{enumerate}
\item
\begin{equation}
 \label{H-sp}
\h_\lambda(\sqrt{q},1/\sqrt{q})=\h_\lambda(q)^{2g-2};
\end{equation}
\item 
\begin{equation}
 \label{H-symm}
 \h_\lambda(-z,-w)=\h_\lambda(z,w), \qquad 
\h_\lambda(w,z)=\h_{\lambda'}(z,w);
\end{equation}
\item
$\h_\lambda$ has a Laurent series expansion in $z$ and $w^{-1}$
\begin{equation}
 \label{H-Laurent}
 \h_\lambda = \sum_{i\geq i_0,j\geq 0} *_{i,j}\,z^{j}w^{-i} \in
 \Z[[z,w^{-1}]][w]
\end{equation}
with  $i_0=-(2g-2)\langle\lambda,\lambda\rangle$.
\end{enumerate}
(To verify the last statement, for example, write
\bes
\h_\lambda(z,w)=w^{(2g-2)\langle\lambda,\lambda\rangle} \prod
\frac{\left(1-z^{2a+1}/w^{2l+1}\right)^{2g}}  
{(1-z^{2a+2}/w^{2l})(1-z^{2a}/w^{2l+2})}  
\ees
and expand each factor of the denominator in a geometric series.)

\subsection{Formal infinite products}
\label{formal-inf-prod}
We will need the following formal manipulations of infinite products.
For a discussion for general $\lambda$-rings see \cite{getzler} whose
notation we will follow.

We first define the crucial maps $\Exp$ and $\Log$ that we need.  Let
$K:=k(x_1,\ldots,x_N)$ be the field of rational functions in the
indeterminates $x_1,\ldots,x_N$ over a ground field $k$ of
characteristic zero. In the ring $K[[T]]$ of formal power series in
another indeterminate $T$ with coefficients in $K$ we consider the
following map
\begin{eqnarray}
 \label{Exp}
 \Exp: TK[[T]] & \longrightarrow & 1+TK[[T]] \\
  V & \mapsto & \exp
\left(\sum_{r\geq 1}\frac 1 rV(x_1^r,\ldots,x_N^r,T^r)\right).
\end{eqnarray}

The map $\Exp$ has an inverse $\Log$ which we now define. Given $F\in
1+TK[[T]]$ let $U_n\in K$ be the coefficients in the expansion
$$
\log(F)=:\sum_{n\geq 1} U_n(x_1,\ldots,x_N)\frac{T^n}n.
$$
Define
\begin{equation}
 \label{U-V-general}
V_n(x_1,\ldots,x_N):=\frac 1 n \sum_{d\mid n} \mu(d)\, U_{n/d}(x_1^d,
\ldots,x_N^d),
\end{equation}
where $\mu$ is the ordinary M\"obius function, and set
$$
\Log(F):=\sum_{n\geq 1} V_n(x_1,\ldots,x_N)\,T^n.
$$
We now prove that $\Exp$ and $\Log$ are indeed inverse maps.

Let $V=\sum_{n\geq 1}V_n(x_1,\ldots,x_N)\,T^n\in TK[[T]]$ then 
\begin{eqnarray*}
\log(\Exp(V))&=&\sum_{n,r\geq 1}  V_n(x_1^r,\ldots,x_N^r)
\,\frac{T^{nr}}r\\
&=&   \sum_{n\geq 1}  \frac 1 n \sum_{d\mid n}
d\, V_d(x_1^{n/d},\ldots,x_N^{n/d})
\,T^n
\end{eqnarray*}
so that
$$
U_n(x_1,\ldots,x_N)=  \sum_{d\mid n}
d\, V_d(x_1^{n/d},\ldots,x_N^{n/d}).
$$
By M\"obius inversion this equality is equivalent to
\eqref{U-V-general}. Therefore $\Log\circ\Exp(V)=V$ and similarly
$\Exp\circ\Log(F)=F$ for $F\in 1+TK[[T]]$.

Note that   $\Exp$ and $\Log$ work the same way if we replace $K[[T]]$
by $S[[T]]$ where 
$$
S:=k[[x_1,\cdots,x_N]][x_1^{-1},\cdots,x_N^{-1}]
$$
is a Laurent series ring.

The connection with infinite products is the following one.  We
clearly have that $\Exp(V+W)=\Exp(V)\Exp(W)$ and
$$
\Exp(x^m\,T^n)=(1-x^m\,T^n)^{-1}, \qquad m=(m_1,\ldots,m_N)\in
\Z^N, n \in \N,
$$where $x^m:=x_1^{m_1}\cdots x_N^{m_N}$.
Now suppose that the coefficients in $V=\sum_{n\geq
 1}V_n(x_1,\ldots,x_N)\,T^n\in TK[[T]]$ have a Laurent expansion
\begin{equation}
\label{L-exp}
V_n(x_1,\ldots,x_N)=\sum_ma_{m,n}\,x^m, \qquad m=(m_1,\ldots,m_N)\in
\Z^N, \quad  a_{m,n}\in k
\end{equation}
in $S$.
Then formally  we may write 
\begin{equation}
 \label{inf-prod-1}
\Exp(V)=
\prod_{m,n}(1-x^m\,T^n)^{-a_{m,n}};
\end{equation}
or, in other words, for $F\in 1+TS[[T]]$ we may think of the
coefficients in $\Log(F)=\sum_{m,n}a_{m,n}\,x^m\,T^n \in TS[[T]]$ as
the exponents of a formal infinite product expansion of $F$ of the
form \eqref{inf-prod-1}.

In fact, we may actually replace $k$ by $\Z$. Let
$$
R:=\Z[[x_1,\cdots,x_N]][x_1^{-1},\cdots,x_N^{-1}].
$$
Then from \eqref{inf-prod-1} we see that $\Exp$ maps $TR[[T]]$ to
$1+TR[[T]]$.

Similarly, $\Log$ maps $1+TR[[T]]$ to $TR[[T]]$.  Indeed, we claim
that any $F\in 1+TR[[T]]$ can be written as a formal infinite
product
\begin{equation}
 \label{inf-prod}
F=\prod_{m,n}(1-x^m\,T^n)^{-a_{m,n}}, \qquad a_{m,n}\in \Z.
\end{equation}
We may in fact find the exponents $a_{m,n}$ recursively as follows.
Order the $m$'s, say lexicographically. Start with $n=1$ and let $m_0$
be the smallest $m$ such that $a_{m,1}\neq 0$. Consider
$F(1-x^{m_0}\,T)^{a_{m_0,1}} \in 1+TR[[T]]$; its coefficient of
$x^{m_0}\,T$ is zero by construction.  Repeat the process with this
series. In the limit we get a series say $F_1\in 1+TR[[T]]$ whose
coefficient of $T$ is zero.  Now set $n=2$ and start all over again
with $F_1$.  In the limit we end up with the constant series $1$ from
which we obtain an expression of the desired form \eqref{inf-prod}.
Clearly  $\Log(F)=\sum_{m,n}a_{m,n}\,x^m\,T^n$ hence, in particular,
the exponents $a_{m,n}$ are uniquely determined by $F$.

Usually given $F=\Exp(V)$ with $V\in TK[[T]]$ we have more than one
choice for what Laurent series ring to consider for the expansion
\eqref{L-exp} of the coefficients of $V$. This may result in at first
puzzlingly different infinite products for the same series $F$.

A typical example is the following. Let $V=T/(1-q)$. If we expand it
in a Laurent series in $q$ we have
$$
V=T\sum_{n\geq 0}q^n \qquad \text{in
$\Z[[q,T]]$}
$$
and hence
$$
F:=\Exp(T/(1-q))=\prod_{n\geq 0} (1-q^nT)^{-1}, \qquad \text{in
$\Z[[q,T]]$}.
$$
On the other hand $V=-Tq^{-1}/(1-q^{-1})$ and hence if we expand it in
a Laurent series in $q^{-1}$ we find
$$
V=-T\sum_{n\geq 1}q^{-n} \qquad \text{in
$\Z[[q^{-1},T]]$}
$$
and hence also
\beq \label{different}
F=\prod_{n\geq 1} (1-q^{-n}T), \qquad \text{in 
$\Z[[q^{-1},T]]$}.
\eeq
This observation becomes important, for example, when comparing
results from different sources.

\section{E-polynomial of $\M_n$}
\label{epolynomial}

\subsection{The irreducible characters of the general linear group
 over a   finite field}
\label{GLn}

Throughout this section $G_n$ will denote the group $\GL{n,\F_q}$ for a
fixed $n\in {\Z_{>0}}$ and finite field $\F_q$ of cardinality $q$. We now recall
the description of the irreducible characters of $G_n$ following 
\cite{Mc}.

Fix an algebraic closure $\overline{\F_q}$ of $\F_q$. For each $r\in {\Z_{>0}}$ let
$\F_{q^r}$ be the unique subfield of $\overline {\F_q}$ of cardinality $q^r$. Let
$\Frob_q\in \Gal(\overline{\F_q} /\F_q)$ be the Frobenius automorphism
$x\mapsto x^q$. Then $\F_{q^r}$ is the fixed field of $\Frob_q^r$. For
$r,s\in {\Z_{>0}}$ with $r|s$ we have the norm map
$\N_{s,r}: \F_{q^s} \longrightarrow \F_{q^r}$, which is surjective.

Let $\Gamma_r$ be the character group of $\F_{q^r}^\times$. Composition with
$\N_{s,r}$, when $r|s$, gives an injective map $\Gamma_r\longrightarrow
\Gamma_s$. Let
$$
\Gamma=\lim_{\rightarrow}\Gamma_r
$$
be the direct limit of the $\Gamma_r$ via these maps. The Frobenius
automorphism $\Frob_q$ acts on $\Gamma$ by $\gamma\mapsto \gamma^q$.
The fixed group of $\Frob_q^r$ is the image of $\Gamma_r$ in $\Gamma$,
which, abusing  notation, we also denote by $\Gamma_r$.

Let $\p_m(\Gamma)$ be the set of all maps $\Lambda:
\Gamma\longrightarrow \p$ which commute with $\Frob_q$ and such that
$$
|\Lambda|:=\sum_{\gamma \in \Gamma}|\Lambda(\gamma)|=m.
$$
Set $\p(\Gamma):=\bigcup_m\p_m(\Gamma)$. Given $\Lambda \in
\p_m(\Gamma)$ we let $\Lambda'\in \p_m(\Gamma)$ be the function with
values $\Lambda'(\gamma):=(\Lambda(\gamma))'$.

For $\gamma \in \Gamma$ we let $\{\gamma\}$ be its orbit in $\Gamma$
under $\Frob_q$ and $d(\gamma)$ be its degree (the size of the orbit).
Given $\Lambda \in \p_m(\Gamma)$ we let $m_{d,\lambda}$ be the
multiplicity of $(d,\lambda)$ in $\Lambda$, where $d\in \Z_{>0}$ and
$0\neq\lambda \in \p$. I.e.,
$$
m_{d,\lambda}:=\#\{\{\gamma\} \,|\, d(\gamma)=d,
\Lambda(\gamma)=\lambda\};
$$
for convenience we also set $m_{d,0}=0$ for all $d$.
We will call the collection of multiplicities $\{m_{d,\lambda}\}$ 
the {\it type} of $\Lambda$ and denote it by $\tau(\Lambda)$. We will
write
$$
|\tau|:=|\Lambda| = \sum_{d,\lambda}  m_{d,\lambda} \,d\, |\lambda|.
$$

\bigskip {\em There is a canonical bijection $\Lambda \mapsto
 \chi_\Lambda$ between $\p_n(\Gamma)$ and the irreducible characters
 of $G_n$. Under this correspondence, the dimension of the irreducible
 representation associated to $\Lambda$ is
\begin{equation}
 \label{dim-fmla}
 \chi_\Lambda(1)=\prod_{i=1}^n(q^i-1)/\prod_{\{\gamma\}}
  q_\gamma^{-n(\Lambda(\gamma)')}   
  \tilde  H_{\Lambda(\gamma)}  
(q_\gamma)  
\end{equation}
where the product is taken over orbits $\{\gamma\}$ of $\Frob_q$ in
$\Gamma$ and $q_\gamma:=q^{d(\gamma)}$.

Moreover, the value of $\chi_\Lambda$ on any central element $\alpha
I_n$ with $\alpha\in \F_q^\times$ and $I_n\in G_n$ the identity matrix is given
by
\begin{equation}
 \label{char-value}
 \chi_\Lambda(\alpha)=\Delta_\Lambda(\alpha)\chi_\Lambda(1), 
\end{equation}
where
\begin{equation}
 \label{def-delta}
 \Delta_\Lambda=\prod_{\gamma \in
 \Gamma}\gamma^{|\Lambda(\gamma)|}\in \Gamma_1.
\end{equation}
}

In particular, note that $\chi_\Lambda(1)$ only depends on the type
$\tau$ of $\Lambda$; we may hence write it as $\chi_\tau(1)$. Let
\begin{equation}
 \label{h-tau}
\h_\tau(q):=  \prod_{\{\gamma\}}
\h_{\Lambda(\gamma)}
(q_\gamma)= \prod_{d,\lambda} \h_\lambda(q^d)^{m_{d,\lambda}}
\end{equation}
where $\tau=\tau(\Lambda)$. Since
$$
|G_n|=q^{\tfrac12n(n-1)}\prod_{i=1}^n(q^i-1)
$$
we have
\begin{equation}
 \label{chi-coeff}
\frac{|G_n|}{\chi_{\tau}(1)}=(-1)^nq^{\tfrac12n^2} \h_{\tau'}(q)
\end{equation}
where $\tau':=\tau(\Lambda')$.

\begin{remark}
\label{alvis-curtis}
 With the above description the Alvis--Curtis duality \cite{alvis},
 \cite{curtis} for characters of $GL(n,\F_{q})$ is simply given by
 $\Lambda\mapsto \Lambda'$.  In particular, as polynomials in $q$
$$
q^{\frac{n(n-1)}{2}} \chi_\Lambda(1)(q^{-1})=\chi_{\Lambda'}(1)(q).
$$

\end{remark}

\subsection{Counting solutions on the general linear group}
\label{count-GLn}

We now apply the results of \S\ref{count} to $G_n=\GL{n,\F_q}$ using the
results of \S\ref{GLn}. We specialize \eqref{sol-num-comm} to this
case and where $z=\alpha I_n$ with $\alpha \in \F_q^\times$. In the resulting
sum on the right hand side we collect all irreducible characters of
the same type $\tau$ and obtain
\begin{equation}
\label{type-sum}
|G_n|\left(\frac{|G_n|}{\chi_\tau(1)}\right)^{2g-2}
\sum_{\tau(\Lambda)=\tau}\Delta_\Lambda(\alpha).
\end{equation}

Our next goal is to compute the sum in the case that $\alpha$ is a
primitive $n$-th root of unity. We will see that a tremendous
cancelation takes place and only relatively few $\Lambda$'s give 
a non-zero contribution.

Assume then that $\F_q$ contains a primitive $n$-th root of unity
$\zeta_n$ and let
\begin{equation}
 \label{def-c-tau}
 C_\tau:=\sum_{\tau(\Lambda)=\tau}\Delta_\Lambda(\zeta_n) 
\end{equation}
To simplify the notation let
\begin{equation}\label{nnq}
N_n(q):= \#\{x_1,y_1,\ldots,x_g,y_g \in \GL{n,\F_q} \;|\;
[x_1,y_1]\cdots[x_g,y_g]\zeta_n=1\}.
\end{equation}
At this point we have in combination with   \eqref{chi-coeff} 
\begin{equation}
 \label{prel-fmla}
\frac1{|G_n|}N_n(q)
 =  \sum_{|\tau|=n} C_\tau \;\left(q^{\tfrac12n^2} \h_{\tau'}(q)\right)^{2g-2}.
\end{equation}
Our next task is to compute $C_\tau$; we will find that $C_\tau$ is a
constant times $(q-1)$, independent of the choice of $\zeta_n$. In
particular, this will show that $N_n(q)/|G_n|$ is a polynomial in $q$.

\subsection{General combinatorial setup}

To compute $C_\tau$ we will use the inclusion-exclusion principle on a
certain partially ordered set. We first describe a slightly more
general setup.

Let $I:=\{1,2,\ldots,m\}$ and let $\Pi(I)$ be the poset of partitions
of $I$; it consists of all decompositions $\pi$ of $I$
into disjoint unions of non-empty subsets $I=\amalg_j I_j$ ordered by
refinement, which we denote by $\leq$. Concretely, $\pi\leq \pi'$
in $\Pi(I)$ if every subset in $\pi$ is a subset of one in $\pi'$. We
call the $I_j$'s the {\it blocks} of $\pi$.

The group $S(I)$ of permutations of $I$ acts on $\Pi(I)$ in a natural
way preserving the ordering $\leq$; for $\rho \in S(I)$ let
$\Pi(I)^\rho$ be the subposet of $\Pi(I)$ of elements fixed by $\rho$.

For $\pi \in \Pi(I)$ let $J$ be the set of its blocks and write
$I=\sqcup_{j\in J} I_j$. It will be convenient to also think of $\pi$
as the surjection $\pi: I \longrightarrow J$ that takes $i$ to $j$
where $I_j$ is the unique block containing $i$. Then the blocks $I_j$
are just the fibers of this map. For $\pi \in \Pi(I)^\rho$ the blocks
of $\pi$ are permuted by $\rho$. Denote by $\rho_\pi$ the induced
permutation in $S(J)$.

Fix a variety $X$ defined over $\F_q$ and let $(I,X):=X^I$ be the
variety of maps  $\xi: I \longrightarrow X$. We have a natural
injection $S(I) \hookrightarrow \Aut((I,X))$. For $\rho \in S(I)$ we let
$(I,X)_\rho$ be the twist of $(I,X)$ by $\rho$. Its
$\F_q$-points consists of the maps
$$
\xi: I \longrightarrow X\left(\Fqbar\right), \qquad \xi\circ\rho=\Frob_q
\circ \xi
$$
Also let $(I,X)_\rho'\subseteq (I,X)_\rho$ be the open
subset of injective maps $\xi: I \longrightarrow X$.

For $\pi \in \Pi^\rho(I)$ we let $(\pi,X)_\rho\subseteq (I,X)_\rho$ be
the closed subset of maps $\xi: I\longrightarrow X$ which are constant
on the blocks of $\pi$. (This notation is consistent with our previous
one if we think of $I$ as the partition where every block has size
$1$, the unique minimal element of $\Pi(I)^\rho$.)  There is a natural
isomorphism
\begin{eqnarray*}
 \label{isom}
\iota_\pi: (J,X)_{\rho_\pi} &\longrightarrow &(\pi,X)_\rho  \\
\xi\qquad & \mapsto &\xi\circ \pi
\end{eqnarray*}

Finally, we let $(\pi,X)_\rho'\subset (\pi,X)_\rho$ be the image of
$(J,X)_{\rho_\pi}'$ under $\iota_\pi$. Concretely, $\pi$ prescribes
some equalities on the values of $\xi: I\longrightarrow X$ (it is
constant on the blocks of $I$) and $\xi \in (\pi,X)_\rho'$ if and only
if these are the {\it only} equalities among these values. It follows
that
\begin{equation}
 \label{X-decomp}
 (I,X)_\rho=\sqcup_\pi (\pi,X)_\rho',
\end{equation}
where $\pi$ runs through the partitions in $\Pi(I)^\rho$. More
generally,
\begin{equation}
 \label{X-decomp-1}
   (\pi^*,X)_\rho=\sqcup_{\pi\leq \pi^*} (\pi,X)_\rho',
\end{equation}
for any $\pi^* \in \Pi(I)^\rho$.

Now take $X$ to be a commutative algebraic group over $\F_q$. In
particular, all the $(\pi,X)_\rho$ are subgroups of $(I,X)_\rho$ and
$\iota_\pi$ is a group isomorphism.  Fix $n \in \Z_{>0}$. Assume that there
exists  a character $\varphi: X(\F_q) \longrightarrow \mubf_n$
of exact order $n$. Let
\begin{eqnarray*}
 \label{psi-defn}
 \Phi:  (I,X)_\rho & \longrightarrow & \mubf_n \\
 \xi & \mapsto & \varphi\left(\prod_{i\in I} \xi(i)^{\eta(i)}\right),
\end{eqnarray*}
where $\eta: I \longrightarrow \Z_{>0}$ is compatible with $\rho$, i.e.,
$\eta\circ\rho=\eta$ (or, equivalently, $\eta$ is constant on the orbits
of~$\rho$) and $\sum_i\eta(i)=n$. Then $\Phi$ is a well defined character
on $(I,X)_\rho$ (the argument of $\varphi$ is in $X(\F_q)$ by the
compatibility of $\xi$ and $\eta$ with $\rho$).

For $\pi \in \Pi(I)^\rho$ and $j \in J$ define
$\eta_\pi(j):=\sum_{i\in I_j} \eta(i)$. It is easy to check that
$\eta_\pi$ is compatible with $\rho_\pi$, i.e., $\eta_\pi\circ
\rho_\pi=\eta_\pi$ and also $\sum_{j\in J}\eta_\pi(j)=n$. Let
$\Phi_\pi$ be the analogue of $\Phi$ for $(J,X)_{\rho_\pi}$
constructed using $\eta_\pi$.  Then
\begin{equation}
 \label{psi-pi}
 \Phi_\pi=\Phi\circ\iota_\pi.
\end{equation}
We will see in the next section that what we need is to compute the
following sum
\begin{equation}
 \label{c-tau}
 S'(I):=\sum_\xi \Phi(\xi),
\end{equation}
where $\xi$ runs over  $(I,X)_\rho'(\F_q)$. Thanks to
\eqref{X-decomp-1} we can calculate $S'(I)$ using the inclusion-exclusion
principle on the poset $\Pi(I)^\rho$:
\begin{equation}
 \label{incl-excl}
 S'(I)=\sum_{\pi\in \Pi(I)^\rho} \mu_\rho(\pi)S(\pi),
\end{equation}
where
\begin{equation}
 \label{S-pi}
 S(\pi):=\sum_\xi \Phi(\xi),
\end{equation}
with $\xi$ running over $(\pi,X)_\rho(\F_q)$, and where 
$\mu_\rho$ is the M\"obius function of $\Pi(I)^\rho$.

The advantage of \eqref{S-pi} over \eqref{c-tau} is that it is a
complete character sum and, hence, vanishes unless the character is
trivial.  Using \eqref{psi-pi} we get
\begin{equation}
 \label{S-pi-1}
 S(\pi)=\sum_\xi \Phi_\pi(\xi),
\end{equation}
with $\xi$ running over $(J,X)_{\rho_\pi}(\F_q)$.  We can now factor
$S(\pi)$ as a product over the orbits of $\rho_\pi$. Each factor
is a complete character sum of the form
$$
\sum_{x\in \F_{q^a}^\times} \varphi^b\circ \N_{\F_{q^a}/\F_q}(x),
$$
where $a$ and $b$ are, respectively, the size and the common value
of $\eta_\pi$ of the corresponding orbit of $\rho_\pi$. Since
$\varphi$ has exact order $n$, by assumption, the character
$\varphi^b\circ \N_{\F_{q^a}/\F_q}$ is trivial if and only if $n\mid
b$; this can only happen if $|J|=1$ because $\sum_{j\in
 J}\eta_\pi(j)=n$ and $\eta_\pi(j)>0$.

It follows that $S(\pi)=0$ unless $\pi$ is the trivial partition
$I=I$, the unique maximal element in $\Pi(I)^\rho$; in this case,
$S(\pi)=|X(\F_q)|$ since $(\pi,X)_\rho(\F_q)=X(\F_q)$. In order to
conclude the calculation we need to know the value of $\mu_\rho$ at
the maximal element of $\Pi(I)^\rho$. For simplicity denote this by
$\bar\mu_\rho$. Its value was computed by Hanlon \cite{Ha}. Abusing
notation let $\rho$ also denote the partition of $m$ determined by its
cycle structure and write it in multiplicity notation
$(1^{m_1}2^{m_2}\cdots)$ where $m_d$ is the number of cycles of size
$d$ in $\rho$. Then we have
\begin{equation}
 \label{mu-value}
 \bar\mu_\rho=\left\{ \begin{array}{ll}
\mu(d) (-d)^{m_d-1}(m_d-1)! & \qquad \rho=(d^{m_d}) \\
0 & \qquad \mbox{otherwise}
\end{array}
\right.,
\end{equation}
where $\mu$ is the ordinary M\"obius function. (To be sure,
$\rho=(d^{m_d})$ means that $\rho$ consists only of $m_d$ cycles of
size $d$ for some $d$.)  Putting this together with \eqref{incl-excl}
we finally obtain
\begin{equation}
 \label{final}
S'(I)=
 \left\{ \begin{array}{ll}
|X(\F_q)|\mu(d) (-d)^{m_d-1}(m_d-1)! & \qquad \qquad \rho=(d^{m_d}) \\
0 & \qquad \mbox{otherwise}
\end{array}
\right..
\end{equation}
Note that the value of $S'(I)$ does not actually depend on the actual
character $\varphi$.

\begin{example}
 To illustrate the previous calculation consider the simplest case
 where $\rho$ is the identity, i.e., assume the action of Frobenius
 is trivial. The situation is the following. Let $X$ be a finite
 abelian group, $\varphi: X \longrightarrow \mubf_n$ be a character
 of exact order $n$ and $(n_1,\ldots,n_m)$ be positive integers such
 that $n_1+\cdots+n_m=n$. Then \eqref{final} reduces to 
$$
\sum_{x_i\neq x_j}\varphi(x_1^{n_1}\cdots
x_m^{n_m})=|X|(-1)^{m-1}(m-1)!,
$$ 
which is not hard to prove directly.
\end{example}
\subsection{Calculation of $C_\tau$}

We now apply the general setup of the previous section to compute
$C_\tau$. We start by describing all $\Lambda \in \p_n(\Gamma)$ with a
given type $\tau$. 

For each $d\in \Z_{>0}$ and $0\neq \lambda \in \p$ let $m_{d,\lambda}$ be
the multiplicity of $(d,\lambda)$ in $\tau$. Let
$$
m_d:=\sum_{\lambda}m_{d,\lambda}, 
\qquad
m:=\sum_dm_d,
$$
then the {\it support} of $\Lambda$, i.e., those $\gamma \in
\Gamma$ with $\Lambda(\gamma)\neq 0$, has size $m$. 

Let $I:=\{1,2,\ldots,m\}$ as in the previous section and fix an
element $\rho \in S(I)$ whose cycle type has $m_d$ cycles of length
$d$ for each $d\in \Z_{>0}$.  Fix also a  map
$$
\nu: I \longrightarrow \p\setminus\{0\}
$$
which is constant on orbits of $\rho$, i.e., $\nu\circ\rho= \nu$,
and such that for any $\lambda \in \p\setminus\{0\}$ and $d\in \Z_{>0}$
there are exactly $m_{d,\lambda}$ orbits of size $d$.

Given an injective map
\begin{equation}
\label{inj-map}
\xi: I \longrightarrow \Gamma, \qquad \xi\circ\rho=\Frob_q
\circ \xi
\end{equation}
there is a uniquely determined $\Lambda \in \p(\Gamma)$ satisfying
$$
\Lambda\circ \xi = \nu.
$$
To check that it is indeed in  $\p(\Gamma)$ note that
$$
(\Lambda\circ \Frob_q) \circ \xi = \Lambda\circ \xi \circ \rho = \nu
\circ \rho =\nu=\Lambda \circ \xi
$$
hence $\Lambda\circ \Frob_q=\Lambda$.
Note also that by construction $\tau(\Lambda)=\tau$ and
$$
n:=|\Lambda|=\sum_{d,\lambda}d\,m_{d,\lambda}\,|\lambda|=\sum_{i\in
 I}|\nu(i)|. 
$$

It is clear that every $\Lambda \in \p_n(\Gamma)$ with
$\tau(\Lambda)=\tau$ arises in this manner ($\xi$ is just a labelling
of the support of $\Lambda$ and $\nu$ fixes its values) but typically
in more than one way. More precisely, the assignment $\xi \mapsto
\Lambda$ is a $z_\tau$ to $1$ map, where $z_\tau$ is the order of the
subgroup consisting of the elements of $S(I)$ which commute with
$\rho$ and preserve $\nu$. It is straightforward to check that
$$
z_\tau = \prod_{d,\lambda}d^{m_{d,\lambda}}\,m_{d,\lambda}!
=\prod_d d^{m_d}\prod_\lambda m_{d,\lambda}!
$$

Take now $X=\Gm$ in the previous section. We may (non-canonically)
identify $\Gamma$ with $X(\Fqbar)$; then the injective maps $\xi$ of
\eqref{inj-map} above correspond to the elements of $(I,X)'_\rho$.
Let $\varphi:\Gm(\F_q) \longrightarrow \mubf_n$ correspond to the
order $n$ homomorphism $\Gamma_1\longrightarrow \mubf_n$ given by
evaluation at $\zeta_n\in \F_q^\times$.  Let $\eta:I \longrightarrow
\Z_{>0}$ be defined by $\eta(i):=|\nu(i)|$. Note that
$$
\sum_i\eta(i)=n.
$$

Then
if $\xi$ corresponds to $\Lambda$ as above we have 
$$
\Phi(\xi)=\Delta_\Lambda(\zeta_n)
$$
and therefore
$$
C_\tau=\frac 1{z_\tau} S'(I).
$$
Hence by \eqref{final}
\begin{equation}
 \label{C-tau-final}
C_\tau=    \left\{ \begin{array}{ll}
(-1)^{m_d-1}(q-1)\frac{\mu(d)}d\frac{ (m_d-1)!}
{\prod_\lambda m_{d,\lambda}!}
& \qquad \rho=(d^{m_d}) \\
0 & \qquad \mbox{otherwise}
\end{array}
\right.
\end{equation}
independently of the choice of $\zeta_n$. 

\subsection{Main formula}
Let
\begin{equation}
   \label{E}
   \epol_n(q):=\frac{N_n(q)}{|\PGL{n,\F_q}|}.
\end{equation}
As we remarked at the end of \ref{count-GLn} $\epol_n$ is a polynomial
in $q$.  To see this it is enough to plug in \eqref{C-tau-final} into
\eqref{prel-fmla}.

\begin{theorem}\label{polynomialcount}
 The variety $\M_n/\C$ has polynomial count and its $E$-polynomial satisfies
 $$E(\M_n/\C;x,y)=\epol_n(xy).$$ 
\end{theorem}
\begin{proof}
 From the definition \eqref{undef} of $\Un$ it is clear that it can
 be viewed as a closed subscheme $\mathcal X$
 of $\GL{n}^{2g}$ over the ring $R:=\Z[\zeta_n,\frac{1}{n}]$. Note that we have extended the base ring in Remark~\ref{scheme}. 
Let $\varphi: R\to \C$ be an embedding, then
 $\mathcal X$ is a spreading out of $\Un/\C$. 

 For every homomorphism \beq \label{extscal} \phi: R\longrightarrow
 \F_q\eeq the image $\phi(\zeta_n)$ is a primitive $n$-th root of
 unity in $\F_q$, because the identity $$\prod_{i=1}^{n-1}(\zeta_n^i-1)=n$$ guarantees
that $1-\zeta_n^i$ is  a unit in $R$ for $i=1,\dots,n-1$, and therefore cannot be zero in the image. (This is why we have extended our ring $R$ from the one in Remark~\ref{scheme}.) Hence all of our previous considerations apply to
 compute $\#{\mathcal X}_{\phi}(\F_q)=N_n(q)$.

 On the other hand the group scheme $\PGL{n,R}$ 
 acts on $\mathcal X$ by conjugation. We define the affine
 scheme
 ${\mathcal Y}= {\rm Spec}(R[{\mathcal X}]^{\PGL{n,R}})$ over $R$.
Because $\varphi:R\to \C$ is a flat morphism  \cite[Lemma 2]{seshadri} implies that ${\mathcal Y}$ is a spreading out
of $\M_n/\C$ over $R$.

 Now take an $\F_q$-point of the scheme ${\mathcal Y}_{\phi}$,
 obtained from $\mathcal Y$ by the extensions of scalars in
 \eqref{extscal}. By \cite[Lemma 3.2]{kac0} the fiber over it in
 ${\mathcal X}_{\phi}(\F_q)$ is non-empty and an orbit of
 $\PGL{n,\F_q}$. The same argument as in the proof of
 Theorem~\ref{smooth} shows that $\PGL{n,\F_q}$ acts freely on
 ${\mathcal X}_{\phi}(\F_q)$. Consequently $$\#{\mathcal
   Y}_\phi(\F_q)=\frac{\#{\mathcal
     X}_\phi(\F_q)}{\#{\PGL{n,\F_q}}}=\frac{N_n(q)}{|\PGL{n,\F_q}|}=
 \epol_n(q)$$
 Thus $\M_n/\C$ has polynomial count. Now the theorem
 follows from Theorem~\ref{katz}.3.
\end{proof}

Let us write
$$
\epol_n(q):=\sum_k \epolc_k^n\,q^k.
$$
We also consider the following normalized version of $\epol_n$
\begin{equation}
 \label{E-bar}
 \epoln_n(q):=q^{-\tfrac12 d_n}\epol_n(q)
:=\sum_k \epolnc_k^n\,q^k
\end{equation}
a Laurent polynomial in $\Z[q,q^{-1}]$.

It will be more convenient to work with the following modified
quantity
\begin{equation}
 \label{V-defn}
V_n(q):=\frac q{(q-1)^2}\;\epoln_n(q)=\frac{q^{-(g-1)n^2}}{(q-1)^2} E_n(q)=
q^{-(g-1)n^2}\frac{N_n(q)}{(q-1)|\GL{n,\F_q}|},  
\end{equation}
(recall that $d_n=\dim(\M_n)=(2g-2)n^2+2$).
For $g>0$ this is a Laurent polynomial in $\Z[q,q^{-1}]$; for $g=0$
we have $N_1(q)=1$ and $N_n(q)=0$ for $n>1$. In this case
$V_1=q/(q-1)^2=\sum_{n\geq 1}n\,q^n$ is a power series in $\Z[[q]]$
and $V_n=0$ for $n>1$. Also $\epoln_1=\epol_1=1$ and
$\epoln_n=\epol_n=0$ for $n>1$.

By the formalism of \S \ref{formal-inf-prod} if we let 
$$
V:=\sum_{n\geq 1} V_n(q)\,T^n, 
\qquad V_n(q)=:\sum_{k\in \Z} v^n_k\;q^k, \qquad v_k^n\in \Z
$$
then
\begin{equation}
 \label{def-zeta-fctn}
\Exp(V) = \prod_{n\geq 1}\prod_{k\in \Z}(1-q^k\,T^n)^{-v^n_k}.
\end{equation}

Define
$$
\overline E :=\sum_{n\geq 1} \epoln_n(q)\,T^n,
$$ and so by \eqref{V-defn} \beq \label{E-V} \frac q{(q-1)^2} \epoln = V.\eeq

Taking $\Exp$ of \eqref{E-V} we get
\begin{equation}
 \label{epol-prod}
\Exp(V)=\prod_{j,n\geq 1}\prod_{k \in \Z} (1-q^{k+j}\,T^n)^{-j\epolnc_k^n}.
\end{equation}

The main result is the following
\begin{theorem} \label{main}
For every $g\geq 0$ we have
\begin{equation}
 \label{main-fmla}
 \sum_{\lambda\in \p} \h_\lambda(q)^{2g-2}\;T^{|\lambda|} =
\prod_{j,n\geq 1}\prod_{k \in \Z} (1-q^{k+j}T^n)^{-j\epolnc_k^n}.
\end{equation}
\end{theorem}
\begin{proof}
 Take the logarithm of the left hand side and write the resulting
 coefficient of $T^n$ as $U_n(q)/n$. Using the multinomial
 theorem we find that
 \begin{equation}
   \label{U-fmla}
\frac{U_n} n= \sum_{m_\lambda}
(-1)^{m-1}(m-1)!\prod_\lambda \frac{\h_\lambda^{(2g-2)m_\lambda}}{m_\lambda!},
\qquad m=\sum_\lambda m_\lambda
 \end{equation}
where the sum is over all $m_\lambda\in \Z_{\geq 0}$
satisfying 
\beq \label{constraintq}
\sum_\lambda m_\lambda|\lambda|=n.
\eeq

On the other hand comparing \eqref{U-fmla} with \eqref{prel-fmla} after plugging
in the value of $C_\tau$ from \eqref{C-tau-final} we obtain
\begin{equation}
\label{V-U}
V_n(q):=\frac1n\sum_{d|n}U_{n/d}(q^d)\;\mu(d).
\end{equation}
An application of the (usual) M\"obius inversion shows that 
$$
U_n(q):=\sum_{d|n}d \;V_{n/d}(q^d).
$$
We have then
$$
\sum_{n\geq 1} U_n(q)\;\frac{T^n}n=\sum_{n,r\geq 1}\frac1r 
V_n(q^r)\; T^{nr}
$$
which together with \eqref{def-zeta-fctn} and \eqref{epol-prod} imply
our claim. 
\end{proof}

The following is an immediate corollary of this result:

\begin{corollary}[Curious Poincar\'e duality]
 \begin{equation}
   \label{duality}
\epoln_n(q^{-1})=\epoln_n(q)    
 \end{equation}
\end{corollary}

\begin{proof}
 Inverting $q$ in \eqref{main-fmla} does not change the left hand
 side by  \eqref{hook-duality}. Hence looking at the right hand side
 we see that $\epolnc_{-k}^n=\epolnc_k^n$, which is equivalent to our
 claim. 
\end{proof}

\begin{remark}
 We should point out that the above duality satisfied by $\epoln_n$
 is ultimately a direct consequence of the Alvis--Curtis duality
 \eqref{alvis-curtis} for characters of $\GL{n,\F_{q}}$.
\end{remark}

\begin{corollary}\label{connected} The variety $\M_n/\C$ is connected. 
\end{corollary} 

\begin{proof}
The statement is clear for $g=0,1$ by Example~\ref{g=0-char-var} and
Theorem~\ref{g=1-char-var}. Therefore we can assume $g>1$ for the
rest of the proof. 

 Corollary~\ref{dimension} says that
  each connected component of $\M_n$ 
has dimension $d_n$. Thus the leading coefficient of $E(\M_n/\C;q)$ 
is the number of components of $\M_n$. By Theorem~\ref{polynomialcount} $E(\M_n/\C;q)=E_n(q)$, so it is enough to determine the leading coefficient of $E_n(q)$. 

For this recall the definition of $U_n(q)$ from \eqref{U-fmla}. It is a Laurent polynomial in $q$. In order to determine its lowest degree term, we
see that the lowest degree term of the summands in \eqref{U-fmla} 
are $$
n(-1)^{m-1}\frac{(m-1)!}{\prod_\lambda m_\lambda!}
q^{-(g-1)\sum_\lambda \langle\lambda,\lambda\rangle m_\lambda}.
$$
\begin{lemma} \label{leadingU} The maximum of $\sum_\lambda
\langle\lambda,\lambda\rangle\, m_\lambda$ under the constraint
\eqref{constraintq} occurs only when $m_\lambda=1$ for
$\lambda=(1^n)$ and $m_\lambda= 0$ otherwise.\end{lemma} \begin{proof}  To see this recall  \eqref{part-pairing} that 
$$
\langle\lambda,\lambda\rangle= \sum_i \lambda_i'^2.
$$
Consider a point in the simplex $\Delta: \sum_i x_i =n, x_i\geq 0$
in $\R^n$ with $\sum_\lambda m_\lambda \, m_i(\lambda')$ coordinates
equal to $i$. (Here $m_i(\lambda')$ is the multiplicity of $i$ in
$\lambda'$.) It now suffices to notice that the maximum distance to the
origin on $\Delta$ occurs at a vertex.
\end{proof}
The lemma implies that the lowest degree term of $U_n(q)$ is $nq^{-(g-1)n^2}$. Formula 
\eqref{V-U} now implies that the lowest degree term of 
$V_n(q)$ is $q^{-(g-1)n^2}$. The definition \eqref{V-defn} gives 
that the constant term of the polynomial $E_n(q)$ is $1$. 
By 
\eqref{duality} the leading term of $E_n(q)$ is $q^{d_n}$. The Corollary follows.
\end{proof}

\subsection{Special cases}

We first work out the $E$-polynomial of $\M_n$ when $n=1,2$ from our
generating function \eqref{main-fmla}. Evaluating \eqref{U-fmla} we
get that $U_1(q)=q^{-(g-1)}(1-q)^{2g-2}.$ A short calculation yields
$E_1(q)=(1-q)^{2g},$ so by
Theorem~\ref{polynomialcount} 
$$
E(\M_1;x,y)=(1-xy)^{2g},
$$
 which is consistent with \eqref{htorus}.

For $n=2$ we again evaluate \eqref{U-fmla} to get 
$$\frac{U_2(q)}{2}=-\frac{1}{2}\h_{(1)}^{2(2g-2)}(q)+\h_{(11)}^{2g-2}(q)+\h_{(2)}^{2g-2}(q)$$
substituting the hook polynomials  \eqref{hook-polynomial} we
get $$\frac{U(2)}{2}=\frac{1}{2}q^{2g-2}(1-q)^{4g-4}+q^{-{4g-4}}(1-q)^{2g-2}(1-q^2)^{2g-2}+q^{-(2g-2)}(1-q)^{2g-2}(1-q^2)^{2g-2}.$$
Using \eqref{V-U}  combined  with \eqref{E-V}, \eqref{E-bar}, and
Theorem~\ref{polynomialcount} we get 
\begin{corollary}\label{e2} The $E$-polynomial of $\M_2/\C$
  is $$E(\M_2/\C;x,y)=E_2(xy),$$ where $$E_2(q)=-\frac{1}{2}
  q^{(2g-2)}
  (1-q)^{4g-2}+(1-q)^{2g}(1-q^2)^{2g-2}+q^{2g-2}(1-q)^{2g}(1-q^2)^{2g-2}-\frac{1}{2}q^{2g-2}(1-q)^2(1-q^2)^{2g-2}.$$ 
\end{corollary}

It is also instructive to consider the special cases $g=0,1$ 
of the theorem
in detail. For $g=0$ the identity  (\ref{main-fmla}) becomes
\begin{equation}
\label{main-fmla-0}
 \sum_{\lambda\in \p} \h_\lambda(q)^{-2}\;T^{|\lambda|}=
\prod_{j\geq 1}(1-q^jT)^{-j}.
\end{equation}
This formula follows from known combinatorial identities. Indeed
$$
s_{\lambda}=\pm\;q^{n(\lambda)}H_\lambda(q)^{-1}
$$
where $s_\lambda$ is the Schur function associated to $\lambda$
evaluated at $x_i=q^{i-1}$ (see \cite[I.3 ex. 2]{Mc}  ) and
$$
H_\lambda(q)=\prod (1-q^h)
$$
is the hook polynomial.  Plugging in $x_j=Tq^j$ in the second
formula of \cite[\S I.4 ex. 2]{Mc}  yields (\ref{main-fmla-0}). This agrees
with our previous calculation: $\epolnc_k^n=1$ for $n=1,k=0$ and zero
otherwise.

For $g=1$ the identity  (\ref{main-fmla}) becomes
\begin{equation}
  \label{g=1-q-fmla} \sum_{\lambda\in \p}T^{|\lambda|}=
  \prod_{n\geq 1}\prod_{r>0}\prod_{s\geq 0}
\frac{(1-q^{r+s}\,T^n)^2}
{(1-q^{r+s-1}\,T^n)(1-q^{r+s+1}\,T^n)},
\end{equation}
which by \eqref{def-zeta-fctn} simplifies to 
$$
\sum_{n\geq 0}p(n)\;T^n = \prod_{n\geq 1}(1-T^n)^{-1}
$$
where $p(n)$ is the number of partitions of $n$ (this is an
identity of Euler).

\begin{remark}
\label{g=1-q}
We deduce that $V_n=1$ for all $n$, when $g=1$. 
Therefore $E(\M_n/\C^n;q)=(q-1)^2$
and $E(\tM_n/\C^n;q)=1$. Because $\tM_n$ is zero dimensional by Theorem~\ref{tmn}  it follows that  $\tM_n$ is a point (cf. 
Theorem~\ref{g=1-char-var}).
\end{remark}

\subsection{Euler characteristic}
\label{eulercharacteristic}
We now prove Corollary~\ref{euler-char}. 

\begin{proof}
By \eqref{htm}, the $E$-polynomial of $\tM_n/\C$ is given by 
$$E(\tM_n/\C;x,y)=\frac{E(\M_n/\C;x,y)}{(xy-1)^{2g}}.$$  By Remark~\ref{eulerchar}, Theorem~\ref{polynomialcount} and Theorem~\ref{main} the Euler characteristic
of $\tilde \M_n/\C$ equals
\beq\label{euler}
\left. \frac
 {N_n(q)}{(q-1)^{2g}|\PGL{n,\F_q}|}\right|_{q=1}.
\eeq
We should point out that the rational function in $q$ in \eqref{euler}
is actually a polynomial, the $E$-polynomial of $\tM_n/\C$.

In terms of $V_n$ we get that \eqref{euler} equals
\beq\label{eulerv}
\left. \frac
 {V_n(q)}{(q-1)^{2g-2}}\right|_{q=1}.
 \eeq
We certainly have that $(q-1)^n$ divides $\h_\lambda(q)$ for
any partition $\lambda$ of $n$. Hence, in the notation of the proof of
\eqref{main-fmla} $(q-1)^{(2g-2)n}$ divides $U_n(q)$ (note that by assumption $2g-2>0$) and it
follows from \eqref{V-U} that $(q-1)^{2g-2}$ divides $V_n(q)$.

We now see that the only contribution in \eqref{V-U} to \eqref{eulerv}
can come from the term $d=n$ and
$$
\left. \frac
 {V_n(q)}{(q-1)^{2g-2}}\right|_{q=1}=
\left. \frac
{\mu(n) U_{1}(q^n)}{n(q-1)^{2g-2}}\right|_{q=1}.
$$
But $U_{1}(q)=(q-1)^{2g-2}$ since $\h_{(1)}(q)=q-1$ for the unique
partition $(1)$ of $1$. We conclude that the Euler characteristic of $\tM_n/\C$ is 
$ \mu(n)\,n^{2g-3}$
finishing the proof.
\end{proof}

\subsection{The untwisted case}
From the above calculation we may now actually deduce the number of
solutions to the untwisted equation (see Remark \ref{hom-divis})
 \begin{equation*}
\#\Hom(\Gamma_g,\GL{n,\F_q}):= \#\{x_1,y_1,\ldots,x_g,y_g \in
\GL{n,\F_q} \;|\; [x_1,y_1]\cdots[x_g,y_g]=1\}.
\end{equation*}
Assume $g>0$ since otherwise $\Gamma_g$ is trivial.  We prove

\begin{theorem} If $\Gamma_g=\pi_1(\Sigma)$ is the fundamental
group of a closed Riemann surface of genus $g>0$, then, using the formalism  of \S\ref{formal-inf-prod}, we have
\begin{equation}
 \label{untwisted}
 \sum_{n\geq 0}
 \frac{\#\Hom(\Gamma_g,\GL{n,\F_q})}{q^{(g-1)n^2}|\GL{n,\F_q}|}\,T^n=\Exp \left((q-1)\Log \left(\sum_{\lambda\in \p} \h_\lambda(q)^{2g-2}\;T^{|\lambda|}\right)\right)
\end{equation}
\end{theorem}
\begin{remark} This gives an explicit formula for the number of representations of $\pi_1(\Sigma)$ to $\GL{n,\F_q}$. The asymptotics
for these numbers as $n$ tends to infinity has been studied
in \cite{liebeck-shalev}.
\end{remark}

\begin{proof}
One way to
express the main formula \eqref{main-fmla} is in terms of zeta
functions of colorings as in \cite{RV}, whose notation we will follow. We
consider colorings on $X=\Gm$ with values on partitions and weight
function
$$
W(\lambda):=\h_\lambda^{2g-2}(q) \in \Z[q,q^{-1}].
$$
We recognize the left hand side of \eqref{main-fmla} as
$Z_C(\bullet,q,T)$ with this setup. Hence the main formula
\eqref{main-fmla} can be written (in the notation of \S
\ref{formal-inf-prod}) as
\begin{equation}
 \label{main-fmla-Log}
 \Log\left(Z_C(\bullet,q,T)\right)=V=\sum_{n\geq 1}V_n(q)\,T^n.
\end{equation}

Similarly, by \eqref{hom-count} and \eqref{chi-coeff} we find that
\begin{equation}
 \label{zeta-color}
 Z_C(\Gm,q,T)=\sum
 \frac{\#\Hom(\Gamma_g,\GL{n,\F_q})}{q^{(g-1)n^2}|\GL{n,\F_q}|}\,T^n.
\end{equation}
In particular this implies that
$$
\#\Hom(\Gamma_g,\GL{n,\F_q})/|\GL{n,\F_q}|\in \Z[q]\ \mbox{ if  }\  g>0
$$
since $W(\lambda)$ is a Laurent polynomial.  This is consistent
with the observation at the end of Remark \ref{hom-divis} that for
$g>0$ $|G|$ always divides $\#\Hom(\Gamma_g,G)$.

{}Using formula \cite[(24)]{RV}  we deduce from \eqref{main-fmla-Log} that
$$
\Log\left(Z_C(\Gm,q,T)\right)=(q-1)\sum_{n\geq 1}V_n(q)\,T^n.
$$
If we take $\Exp$ of both sides we get
$$Z_C(\Gm,q,T)=\Exp((q-1)V).$$
This proves \eqref{untwisted} assuming
that $\F_q$ contains a primitive $n$-th root of unity. However, as we
pointed out, the coefficients on the left hand side of
\eqref{untwisted} are Laurent polynomials in $q$ so the statement is
true for all $q$. The Theorem follows.
\end{proof}

\section{Mixed Hodge polynomial of $\M_n$}
\label{mixedhodgen}
\subsection{Cohomology of $\M_n$}
\label{cohomologymn}

In this section we take $\K=\C$. According to Theorem~\ref{tmn}  $H^*(\M_n)=H^*(\tM_n)\otimes H^*(\M_1),$ where $\M_1\cong (\C^\times)^{2g}$, and the factor $H^*(\M_1)$ is generated by $2g$ degree $1$ classes $\epsilon_j\in H^1((\C^\times)^{2g})$ for $j=1,\dots,2g$. By slight abuse of notation we use the same notation for the  corresponding
classes in $\epsilon_j\in H^1(\M_n)$  for $j=1,\dots,2g$.

To get more interesting cohomology classes on $\M_n$, we construct
cohomology classes in $H^*(\tM_n)\cong H^*(\Mp_n)^{\bbmu_n^{2g}}$. 
We  construct a differentiable  
principal bundle over $\M_n^\prime \times \Sigma$ by 
following \cite{jeffrey}.  
Let $\bar{G} =
\PGL{n,\C}$.  Any $\rho \in \Unp$ induces a
well-defined homomorphism $\pi_1(\Sigma) \to \bar{G}$.  Let
$\tilde{\Sigma}$ be the universal cover of $\Sigma$, which is acted on by
$\pi_1(\Sigma)$ via deck transformations.  There is then a free action of
$\pi_1(\Sigma) \times \GL{n,\C}$ on $\bar{G} \times \Unp \times \tilde{\Sigma}$ given by
$$(p,g) \cdot (h,\rho,x) = (\overline{g} \rho(p) h, \overline{g} \rho
\overline{g}^{-1}, p\cdot x),$$
where $\overline{g}$ denotes the image
of $g$ in $\bar{G}$.  This action commutes with the action of $\bbmu^{2g}$ on $\Unp$. The quotient is the desired ($\bbmu^{2g}$-equivariant) principal
$\bar{G}$-bundle on $\Mp_n$, which we denote by $\U$.  Like any principal $\bar{G}$-bundle, it
has characteristic classes $\bar{c}_2(\U), \dots, \bar{c}_r(\U)$, where
$\bar{c}_i(\U) \in H^{2i}(\Mp_n \times \Sigma)^{\bbmu_n^{2g}}$.  In terms of formal Chern roots
$\xi_k$, $\bar{c}_i$ can be described as the $i$th elementary
symmetric polynomial in the $\xi_k - \zeta$, where $\zeta$ is the
average of the $\xi_k$.  In particular $\bar{c}_1=0$.

Now let $\si \in H^2(\Sigma)$ be the fundamental
cohomology class, and let $e_1, \dots, e_{2g}$ be a standard symplectic basis
of $H^1(\Sigma)$ .  
In terms of these, each of the characteristic classes 
has a K\"unneth decomposition
\beq \label{defineuniversal}\bar{c}_i(\U) = \al_i \si + \be_i + \sum_{j=1}^{2g} \psi_{i,j} e_j, \eeq
defining classes $\al_i \in H^{2i-2}(\Mp_n)^{\bbmu_n^{2g}}\subset H^{2i-2}(\M_n)$, $\be_i \in H^{2i}(\Mp_n)^{\bbmu_n^{2g}}\subset H^{2i}(\M_n)$, and
$\psi_{i,j} \in H^{2i-1}(\Mp_n)^{\bbmu_n^{2g}}\subset H^{2i-1}(\M_n)$ for $i=2,\dots,n$.  In \cite{markman} Markman proves that
\begin{theorem} The classes $\epsilon_j$; $\al_i$, $\psi_{i,j}$ and $\beta_i$ generate $H^*(\M_n)$. 
\end{theorem}

\begin{construction} \label{suspensionconstruction} For what follows we need the following construction. Let $f:Y\to X$, $x\in X$ and $F=f^{-1}(x)$. Then we have the following commutative diagram 
$$
\begin{array}{ccccccc} H^i(Y)&\stackrel{i_F^*}{\longrightarrow}& H^i(F)&\stackrel{d}{\longrightarrow}&H^{i+1}(Y,F)&\stackrel{i^*_{Y}}{\longrightarrow}&H^{i+1}(Y) \\ 
&&&& \mbox{\scriptsize $f^*$}\uparrow&{\nwarrow}^{q^*}& \mbox{\scriptsize $f^*$} \uparrow
\\
&&&&H^{i+1}(X,x)&\stackrel{i^*_{X}}{\cong}& H^{i+1}(X)
\end{array}
$$
Here the first row is the cohomology long exact sequence of the pair $(Y,F)$. 
The second row is the cohomology long exact sequence
of the pair $(X,x)$,  the map $i_X^*$ is an isomorphism for the map $H^i(X)\to H^i(x)$ is always surjective. Finally 
$q^*=f^*(i_X^*)^{-1}:H^{i+1}(X)\to H^{i+1}(Y,F)$. By the commutativity of the
diagram $q^*$ induces a map $\ker(f^*)\to \ker(i_Y^*) \cong 
\im(d)\cong \coker(i_F^*)$. We denote the resulting map $$\sigma^f:\ker(f^*)\to \coker(i_F^*).$$ 

We can also give the map $\sigma$ in terms of cochains. Let  $x\in C^{i+1}(X)$ be
a cocycle. Then $f^*(x)$ 
will be a cocyle in $C^{i+1}(Y)$ vanishing on $F$. if $[x]\in \ker(f^*)$ then $f^*(x)$ is exact. Let $\tilde{y} \in C^i(Y)$ be a cochain such that  $d\tilde{y}=f^*(x).$ Let $y=i_F^*(\tilde{y})$ . Then $dy=i_F^*(f^*(x))=0$, so $y$ is a cocycle. We can define $\sigma([x])=[y]$. 
\end{construction}

\begin{example}
When $f$ is a fibration the map $\sigma$ is called the
{\em suspension map} (see \cite[\S 8.2.2 page 298]{mccleary}). When the fibration is the path fibration of the space
$X$ then $Y=PX$ is the based path space and so contractible, while $F$ can homotopically be identified with
the based loop space $\Omega X$. In this case the suspension map $\sigma:
H^{i+1}(X)\to H^i(\Omega X)$  can be identified with the map 
$\sigma=p_*ev^*$, where $ev:S^1\times \Omega X\to X$ is the evaluation map and $p:S^1\times \Omega X\to \Omega X$ is the projection.  A particular case of the path fibration is the universal bundle: $\pi:EG\stackrel{G}{\longrightarrow} BG$ for a connected $G$ complex linear group $G$. Here $G\sim \Omega BG$ and $EG\sim PBG$. The
suspension map then is a map $\sigma^\pi:H^{i+1}(BG)\to H^i(G).$ 
\end{example}
\begin{remark} We will also need an equivariant version of this construction. 
If we assume that $G$ is a topological group, which
acts on $(X,x)$ and $Y$ in a way so that $f$ is equivariant, then
we have the same diagram and construction above in equivariant cohomology. This way we get the equivariant map $$\sigma^f_G
:\ker_G(\pi^*)\to \coker_G(i_F^*).$$ 
In particular,  $G$ acts on itself by conjugation and consequently on $BG$ and $EG$ making the fibration $\pi$ equivariant. We will then
have the equivariant suspension map \beq \label{universalsigma} \sigma^\pi_G:H^{i+1}_{G} (BG)\to H_G^i(G).\eeq

\end{remark}

\begin{lemma} When $X,Y$ are complex algebraic varieties and $f$ is
algebraic then $\ker(f^*)$, $\coker(i_F^*)$ have natural mixed Hodge structures, and $\sigma$ preserves it. Additionally if a complex linear
group $G$ acts on $(X,x)$ and $Y$ so that $f$ is equivariant then $\ker_G(f^*)$, $\coker_G(i_F^*)$ have natural mixed Hodge structures, and $\sigma_G$ preserves it.
\end{lemma} 
\begin{proof} Deligne in  \cite[Example 8.3.8]{De-Hodge III} constructs a mixed Hodge structure
on relative cohomology, and shows in \cite[Proposition 8.3.9 ]{De-Hodge III} that all maps in  the cohomology long exact 
sequence of a pair preserve mixed Hodge structure. The first statement follows.  

For the second statement Deligne constructs in  \cite[Theorem 9.1.1]{De-Hodge III} a mixed Hodge structure on $H^*(BG)$, by considering a model for $BG$ as a simplicial scheme. Similarly
one can construct the Borel construction $X\times_G EG$ as a simplicial scheme, which will give a mixed Hodge structure on $H^*_G(X)=H^*( X\times_G EG)$ see e.g. \cite{franz-weber}. Then we see
that all maps in the equivariant cohomology sequence of a pair preserve mixed Hodge structures. In turn we get that  $\sigma_G$ too preserves mixed Hodge structures.  \end{proof}

\begin{definition} We say that a  cohomology class $\gamma\in H^i(X)=H^i(X,\Q)$ or $\gamma\in H^i(X;\R)$ has {\em homogenous weight} $k$ if  its complexification  satisfies $\gamma^\C=\gamma\otimes 1\in {W_{2k}}H^i(X)^\C\cap F^{k}H^i(X;\C)$. \end{definition}

\begin{remark}\label{weightremark} 
Note that if  $\gamma\in H^i(X)$ and $\gamma^\C\in {W_{l}}^\C \cap F^m$  with  $2m>l$  then ${Gr^{W^\C}_{l}} \cap F^m \cap \overline{F^{m}}=0$ by \eqref{conjugates} and $\overline{\gamma^\C}=\gamma^\C$ imply that $\gamma^\C\in W_{l-1}^\C$. By induction we have $\gamma=0$. 
Thus we get that
\beq\label{weightlemma} \mbox{ if $\gamma\in H^i(X)$ has homogenous weight $k$ and $\gamma^\C\in F^{k+1}$ or $\gamma\in W_{2k-1}$  then $\gamma=0$.}\eeq
\noindent In particular, a non-zero cohomology class cannot have different homogenous weights.  Moreover as the cup-product preserves mixed Hodge structures  by Theorem~\ref{mixedhodgeproperties}.\ref{kunneth} we have that if $\gamma_1$ has homogenous weight $l_1$ and $\gamma_2$ has homogenous weight $l_2$ then $\gamma_1\cup \gamma_2$ has homogenous weight $l_1+l_2$. In particular we see that if the cohomology $H^*(X)$ of an algebraic variety is generated by classes with homogenous weight then the MHS on $H^*(X)$ will be of type $(p,p)$ i.e. \beq\label{hodgetate}W_l H^*(X)^\C\cap F^mH^*(X;\C) =0, \mbox{ when }2m>l.\eeq Namely if $0\ne x=\sum_k a_k+ib_k\in W_l H^*(X)^\C\cap F^mH^*(X;\C)$ 
with $a_k,b_k\in H^*(X;\R)$ homogenous of weight $k$ then we can consider $k_{min}:=\min_k\{a_k+ib_k\neq 0\}$ and $k_{max}=\max_k\{a_k+ib_k\neq 0\}$ and get $m\leq k_{min}\leq k_{max}\leq l/2$ from 
\eqref{weightlemma}.

Finally for a complex algebraic map $f:X\to Y$ the map
$f^*:H^*(Y)\to H^*(X)$  preserves mixed Hodge structures, we have that if $\alpha\in H^*(Y)$ has homogenous weight $l$ so does $f^*(\alpha)$.

\end{remark}

Now we determine the weights of the universal generators. 
First we know from \eqref{htorus} that the
homogenous weight  of $\epsilon_j$ is $1$. To determine the weight  of the rest of the universal classes we will use Jeffrey's \cite{jeffrey} group cohomology description of them as interpreted in \cite{bott-etal,meinrenken,racaniere}. 

We note that \cite{jeffrey,bott-etal,meinrenken,racaniere} work
with the compact groups ${\rm SU}(n)$, however the arguments are correct with complex groups too. Another way to see that Jeffrey's formulas \eqref{jeffrey1} , \eqref{jeffrey2} and \eqref{jeffrey3} for the universal classes are valid for $G:=\SL{n,\C}$  is to note that Lemma~\ref{bound} below implies that the natural inclusion map of the twisted 
${\rm SU}(n)$-character variety  into the twisted $\SL{n,\C}$-character variety $\M^\prime_n$ induces an isomorphism on ($\bbmu_n^{2g}$-invariant) cohomology below degree $2(g-1)(n-1)+2$. Now $2(g-1)(n-1)+2$  is larger than the degree of any universal class, except possibly of $\beta_n$ (which has degree $2n$), when $g=2$. However Jeffrey's formula for $\beta_n$ is trivially correct
for the complex character varieties as we will see below. Another
difference in our application of \cite{jeffrey, bott-etal,meinrenken,racaniere} is that 
we work on the level of cohomology instead 
of differential forms or cochains, but our cohomological interpretation of 
\cite{jeffrey,bott-etal,meinrenken,racaniere} is straightforward using the last paragraph in Construction~\ref{suspensionconstruction}. 

The easiest is to determine the 
weight of the $\beta_k$. By their construction $\beta_k=c_k(\U |_{\Mp_n\times \{p\}})$ are the Chern classes of the differentiable 
$\PGL{n,\C}$-bundle 
$\U$ constructed above, restricted to $\Mp_n\times \{p\}$, where $p$ is a point on 
$\Sigma$. It is straightforward to identify  $\U |_{\Mp_n\times \{p\}}$ with the $\bar{G}$-bundle $\pi^\prime_n:\Unp\to \Mp_n$ in \eqref{categoricalprime}, thus \beq\label{jeffrey1} \beta_k=c_k(\Unp).\eeq Now $\pi_n'$  is an algebraic principal bundle, 
therefore its Chern classes are pulled back from $H^*(B\bar{G})$ by a complex algebraic map. 
It now follows from  \cite[Theorem 9.1.1]{De-Hodge III}
that the homogenous weight  of $\beta_k$ 
is indeed $k$.

We next determine the weight of the $\psi_{k;j}$. Let 
$\bar{c}_k\in H_G^{2k}(B{\bar G})$ be the $k$-th equivariant Chern class of the $\bar{G}$-equivariant bundle $\pi:E\bar{G}\to B\bar{G}$. 
Clearly $H_G^*(B\bar{G})\cong H^*(B(G\ltimes_\phi\bar{G}))$, where $\phi:G\to \Aut(\bar{G})$ is given by conjugation. By  \cite[Theorem 9.1.1]{De-Hodge III} $\bar{c}_k$ has homogeneous weight $k$. Using the map \eqref{universalsigma} we construct the class $\eta_G^k=\sigma_G(\bar{c}_k)\in H^{2k-1}_{G}(G)$.   It follows that $\eta^k_G$ has  homogenous weight  $k$. Let $p_j:G^{2g}\to G$ be the projection to the $j$th factor, which is equivariant with respect to the conjugation action of $G$. 
Thus $p_j^*(\eta^k_G)\in H^{2k-1}_G(G^{2g})$ has homogenous weight $k$. If $i$ denotes the $G$-equivariant embedding of $\Unp$ into $G^{2g}$, then we have that $i^*p_j^*(\eta^k_G)\in H^{2k-1}_G(p)$ has homogenous weight $k$. Now \cite[Theorem 3.2]{racaniere} implies that \beq\label{jeffrey2} \psi_{k,j}= i^*p_j^*(\eta^k_G) \in H^{2k-1}_G(\Unp)\cong H^{2k-1}(\Mp_n).\eeq Thus $\psi_{k,j}$ has
homogenous weight $k$ as claimed. 

To calculate the weight  of $\alpha_k$ we recall Construction~\ref{suspensionconstruction}. for the case $Y=G^{2g}$, $X=G$, $\pi=\mu^\prime$
and $x=\zeta_n I_n\in G$. Then $F=\Unp$. This gives a map $\sigma_G:\ker_G(\pi^*)\to \coker_G(i_{F}^*)$. Now it follows from \cite[Lemma 2.4]{racaniere} that
$$\eta_G^k\in \ker_G(\pi^*)\subset H^{2k-1}_G(G)$$ and also by \cite[Theorem 3.2]{racaniere} that \beq \label{jeffrey3} \sigma^{\pi}_G(\eta_G^k)=p(\alpha_k),\eeq where $$p:H_G^{2k-2}(\Unp)\to \coker^{2k-2}_G(i_{F}^*)$$ denotes the projection. We know that the homogenous weight of  $\eta_G^k$ is $k$  and
so $p(\alpha_k)$ has homogenous weight $k$.  By the previous paragraph $\im(i_{F}^*)\subset
H_G^*(\Unp)\cong H^*(\Mp_n)$ is exactly the subring generated by the $\psi_{k,j}$
and $\beta_k$ for $k=2,\dots,n$ and $j=1,\dots,2g$. This shows in particular that $p$ is an isomorphism when $k=2$.  Thus 
the homogenous weight of $\alpha_2\in H^2(\Mp_n)$ is $2$.

We summarize our findings in the following

\begin{proposition} \label{universalweights} The cohomology classes $\epsilon_j$ have homogenous weight $1$, while $\psi_{k;j},\beta_k$ have homogenous weight $k$. Finally $\alpha_2$ has homogenous weight $2$ and $p(\alpha_k)\in \coker^{2k-2}_G(i^*_F)$ have homogenous weight $k$.
\end{proposition}

\begin{remark} It is most probably true that $\alpha_k$ has homogenous weight $k$ even for $k>2$, the result for $p(\alpha_k)$ however will suffice for our purposes. Here we show that $p(\alpha_k)\ne 0$.  By the previous paragraph $\im(i_{F}^*)\subset
H_G^*(\Unp)\cong H^*(\Mp_n)$ is exactly the subring generated by the $\psi_{k,j}$
and $\beta_k$ for $k=2,\dots,n$ and $j=1,\dots,2g$. \noindent Because the degree of $\alpha_k$ is $2k-2\leq 2n-2\leq2(g-1)(n-1)$ Lemma~\ref{bound} below implies that $\alpha_k\notin \im(i_F^*)$ i.e. $p(\alpha_k)\ne 0$.

\end{remark}

\begin{corollary} \label{purecorollary} The pure part 
$PH^*(\M_n)=\oplus_k W_k H^k(\M_n)$ is generated by
the classes $\beta_i\in H^{2i}(\M_n)$ for $i=2,\dots,n$. 
\end{corollary} 
\begin{proof} The previous Proposition shows that among the $\psi_{i;j}$ and $\beta_i$ only the classes $\beta_i\in H^{2i}(\M_n)$ are pure classes,
i.e., have pure homogenous weight $i$. This shows that the pure part of  the subring $\im(i^*F)\subset H^*_G(\Unp)\cong H^*(\M^\prime_n)$ they generate is generated by the $\beta_i$ classes.  Moreover  the $\bbmu_n^{2g}$-invariant part of $\coker_G(i^*_F)$ is generated by the classes $p(\alpha_i)$ none of
which has pure homogenous weight. Thus the pure part of $H^*(\M_n^\prime)^{\bbmu_n^{2g}}\cong H^*(\tM_n)$ is generated by the classes $\beta_i$.  By Theorem~\ref{tmn} the result follows.
\end{proof}
\begin{corollary} \label{hodgetriv} The cohomology of $\M_n$ is of type $(p,p)$ , i.e., 
$h^{p,q;j}(\M_n)=0$ unless $p=q$. In particular $H(\M_n;x,y,t)$ is a polynomial in $xy$ and $t$. In the notation of \eqref{hqt}
$$H(\M_n;x,y,t)=H(\M_n;xy,t).$$
\end{corollary}
\begin{proof} By Remark~\ref{weightremark} and Proposition~\ref{universalweights} we know that both the $\bbmu_n^{2g}$-invariant part of $\coker_G(i^*_F)$, which is generated by the classes $p(\alpha_i)$ and the subring $\im(i^*_F)\subset H^*_G(\Unp)^{\bbmu_n^{2g}}\cong H^*(\M_n^\prime)^{\bbmu_n^{2g}}$ generated by the $\beta_i$ and $\psi_{i;j}$ have MHS of type $(p,p)$ in other words \eqref{hodgetate} holds. Thus $H^*(\M_n^\prime)^{\bbmu_n^{2g}}\cong H^*(\tilde{\M}_n)$ has MHS of type $(p,p)$. By Theorem~\ref{tmn} so does $H^*(\M)$.
\end{proof} 
\begin{lemma}\label{bound} There are no relations among the universal generators in the  cohomology of $H^*(\Mp_n)$ until degree $2(g-1)(n-1)+2$. \end{lemma} \begin{proof} This follows from the same statement for the twisted ${\rm SU}(n)$ character variety, which in turn follows from \cite[(7.16)]{atiyah-bott}.\end{proof} 
\subsection{Main Conjecture}
Recall the definition of the $\h_\lambda$ from \eqref{H-defn} and
its properties thereafter.

Let $U_n(z,w)$ be defined by
$$
\log\left(\sum_\lambda
\h_\lambda(z,w)\,T^{|\lambda|}\right)=
\sum_{n\geq 0}U_n(z,w)\frac{T^n}n.
$$
As in \eqref{U-fmla} we find that
\beq \label{Uzw-fmla}
\frac{U_n(z,w)} n= \sum_{m_\lambda}
(-1)^{m-1}(m-1)!\prod_\lambda
\frac{\h_\lambda(z,w)^{m_\lambda}}{m_\lambda!}, 
\qquad m=\sum_\lambda m_\lambda
\eeq
where the sum is over all $m_\lambda\in \Z_{\geq 0}$
satisfying 
\begin{equation}
 \label{constraint}
\sum_\lambda m_\lambda|\lambda|=n.  
\end{equation}
Expanding $U_n(z,w)$ in Laurent series in $z,w^{-1}$ as in
\eqref{H-Laurent} we see that the leading term in $w^{-1}$ of the 
summand is
$$
n(-1)^{m-1}\frac{(m-1)!}{\prod_\lambda m_\lambda!}
w^{(2g-2)\sum_\lambda \langle\lambda,\lambda\rangle m_\lambda}.
$$ From Lemma~\ref{leadingU} it follows  that the leading term of $U_n$ in $w^{-1}$ is
$nw^{(2g-2)n^2}$.

Let 
\beq \label{Vzw-Uzw}
V_n(z,w):=\frac1n \sum_{d\mid n}\mu(d)\,U_{n/d}(z^d,w^d).
\eeq
By the formalism explained in \S \ref{formal-inf-prod}  we know that
\begin{equation}
 \label{V-exp}
\sum_\lambda
\h_\lambda(z,w)\,T^{|\lambda|}=
\exp\left(\sum_{k,n\geq 1} V_n(z^k,w^k)\,\frac{T^{nk}}k\right).
\end{equation}
From our previous calculation we deduce that the leading term in
$w^{-1}$ of $V_n$ is $w^{(2g-2)n^2}$.

Let also
\begin{equation}
 \label{H-defn-1}
\hpoln_n(z,w):=(z^2-1)(1-w^2)V_n(z,w).  
\end{equation}
Both $V_n$ and $\hpoln_n$ are rational functions of $z$ and
$w$.  We should remark that by \eqref{H-sp} we have
\begin{equation*}
V_n(\sqrt{q},1/\sqrt{q})=V_n(q)
\end{equation*}
and therefore
\begin{equation}
 \label{H-sp1}
\hpoln_n(\sqrt{q},1/\sqrt{q})=\epoln_n(q).
\end{equation}

From \eqref{H-symm} we deduce that
\begin{equation*}
 \label{V-symm}
V_n(w,z)=V_n(z,w), \qquad V_n(-z,-w)=V_n(z,w)  
\end{equation*}
and
\begin{equation}
 \label{h-symm}
\hpoln_n(w,z)=\hpoln_n(z,w), \qquad
\hpoln_n(-z,-w)=\hpoln_n(z,w).
\end{equation}
We expand $V_n$ and $\hpoln$ as Laurent series in $z,1/w$
$$
V_n(z,w)=\sum_{i\geq i_0,j\geq 0}v_{i,j}^n \, z^jw^{-i},
\qquad
\hpoln_n(z,w)=\sum_{i\geq -d_n,j\geq 0 }\hpolnc_{i,j}^n \, z^jw^{-i}
$$
where $i_0=-(2g-2)n^2$ and hence $i_0-2=-d_n=-\dim(\M_n)$. Our
calculation of leading terms implies that $\hpolnc^n_{-d_n,j}=1$ for
$j=0$ and is $0$ otherwise.

In terms of these coefficients we can write our generating function as
the infinite products
\begin{eqnarray}
\nonumber
\sum_\lambda
\h_\lambda(z,w)\,T^{|\lambda|}&=&
\prod_{n\geq 1}\prod_{i\geq i_0,j\geq
 0}(1-z^j w^{-i}\,T^n)^{-v_{i,j}^n}\\
 &=& 
\prod_{n\geq 1}\prod_{r>0,s\geq 0}\prod_{i\geq -d_n,j\geq 0}
(1-z^{2s+j}w^{-(2r+i)}\,T^n)^{-\hpolnc_{i,j}^n}
\label{zw-main-fmla}
\end{eqnarray}
Our main conjecture is the following
\begin{conjecture}
\label{main-conj}
 \begin{equation}
   \label{conj}
H(\M_n;q,t)=(t \sqrt q)^{d_n}\hpoln_n\left( \sqrt 
 q, \frac {-1}{t \sqrt q }, \right).
 \end{equation}
\end{conjecture}

\begin{remark}
 In view of \eqref{H-sp1} and \eqref{h-symm}  Conjecture~\ref{main-conj} is true
 specialized to $t=-1$ as it reduces to \eqref{main-fmla}.\end{remark}

Because of the second identity in \eqref{h-symm} and because $d_n$ is even by Corollary~\ref{dimension}, we have that 
the RHS of \eqref{conj} is actually a rational function in $q$.
The geometric Conjecture~\ref{main-conj} implies the following combinatorial
conjectures
 \begin{conjecture}
\label{comb-conj}
   \begin{enumerate}
   \item 
$\hpoln_n(z,w)$ is a polynomial in $z,w$.
\item
The coefficients $(-1)^j\hpolnc_{i,j}^n$ of $\hpoln_n(z,-w)$ are
non-negative integers. 
   \end{enumerate}
\end{conjecture}

\label{deducecpd} In light of \eqref{h-symm}, our main Conjecture~\ref{main-conj} implies the following.

\begin{conjecture}[Curious Poincar\'e Duality] \label{cpd} $$H\left(\M_n;\frac{1}{qt^2},t\right)=(qt)^{-{d_n}}H(\M_n;q,t)$$
\end{conjecture}

\begin{remark} When $t=-1$, this formula specializes to the known
Corollary~\ref{duality}.
\end{remark}

\begin{remark} \label{toppure} On the level of mixed Hodge numbers this conjecture
is equivalent to \beq \label{cpdnumbers} h^{p,p;k}(\M_n)=h^{d_n-p,d_n-p;d_n+k-2p}(\M_n).\eeq
Because $\M_n$ is non-singular, $h^{p,p;k}(\M_n)=0$ for $2p<k$. Dually \eqref{cpdnumbers} implies that $h^{p,p;k}(\M_n)=0$ for $k>d_n$.  
The vanishing of the cohomology of $\M_n$ above middle dimension can be deduced from the fact that $\M_n$ is diffeomorphic to the
space of twisted flat $\GL{n,\C}$-connections on the Riemann surface $\Sigma$, which is a Stein manifold with its natural hyperk\"ahler
metric \cite{hitchin}. 

In particular \eqref{cpdnumbers} implies that the 
 pure mixed Hodge numbers $h^{p,p;2p}(\M_n)$ should be curious Poincar\'e dual to $h^{d_n-p,d_n-p;d_n}(\M_n)$, i.e., to the mixed Hodge numbers of the middle (top non-vanishing) cohomology of $\M_n$.
\end{remark}

Finally we have a geometric conjecture which would imply the 
above curious Poincar\'e duality. Define the Lefschetz map $L:H^i(\tM_n)\to H^{i+2}(\tM_n)$ by
$x\mapsto \alpha\cup x$, where $\alpha=\alpha_2$ is the universal class in $H^2(\tM_n)$ defined in \eqref{defineuniversal}. As it respects mixed Hodge structures by Theorem~\ref{mixedhodgeproperties}.\ref{kunneth}  and $\alpha$ has homogenous weight $2$ by Proposition~\ref{universalweights} it defines a map on the graded pieces of the homogenous weight filtration 
$L: Gr^W_{l} H^i(\tM_n)\to Gr^W_{l+4}H^{i+2}(\tM_n)$. 
\begin{conjecture}[Curious Hard Lefschetz] \label{curioushardconj} Recall that $\tilde{d}_n=\dim(\tM_n)=(n^2-1)(2g-2)$. Then $$L^{l}:Gr^W_{{\tilde{d}}_n-2l} H^{i-l}(\tM_n)\to Gr^W_{{\tilde{d}}_n+2l} H^{i+l}(\tM_n)$$ is an isomorphism.\end{conjecture}

\begin{remark} \label{curioushardremark} Here we prove a consequence of this conjecture. As $\tM_n$ is an orbifold  and the non-trivial
weights in the weight filtration on $H^*(\tM_n)$ are even by Proposition~\ref{universalweights}, we have that
for $0<k\leq \tilde{d}_n/2$ 
\bes Gr^W H^{\tilde{d}_n/2-k} (\tM_n) 
&=& \bigoplus^{[\tilde{d}_n/4-k/2]}_{j=0}Gr_{d_n-2k-2j}^W H^{\tilde{d}_n/2-k}(\tM_n)
.\ees Conjecture~\ref{curioushardconj} says that the map $$L^{k+j}:Gr_{d_n-2k-2j}^W H^{\tilde{d}_n/2-k}(\tM_n)\to Gr_{d_n+2k+2j}^W H^{\tilde{d}_n/2+k+j}(\tM_n)$$ is an isomorphism.  This implies 
that $$L^k:Gr_{d_n-2k-2j}^W H^{\tilde{d}_n/2-k}(\tM_n)\to Gr_{d_n+2k-2j}^W H^{\tilde{d}_n/2+k}(\tM_n)$$  is injective. Thus Conjecture~\ref{curioushardconj} implies that the map $$L^k: H^{\tilde{d}_n/2-k}(\tM_n)\to H^{\tilde{d}_n/2+k}(\tM_n)$$ is an injection. This statement follows from \cite[Corollary 4.3]{hausel2} (cf. also \cite[Remark 4.4]{hausel2}) when applied to the moduli space of Higgs bundles diffeomorphic to $\tM_n$. 
\end{remark}

\subsection{Special cases of the main conjecture}
\label{g01}
First we verify the cases of $n=1,2$ of Conjecture~\ref{main-conj}.
{} From \eqref{Vzw-Uzw} and
\eqref{Uzw-fmla} $$V_1(z,w)=U_1(z,w)=\frac{(z-w)^{2g}}{(z^2-1)(1-w^2)}.$$
By \eqref{H-defn-1} $$\hpoln_1(z,w)=(z-w)^{2g}.$$ Thus
Conjecture~\ref{main-conj}
implies $$H_1(\M_1,q,t)=(t\sqrt{q})^{2g}\left(\sqrt{q}+\frac{1}{t\sqrt{q}}\right)^{2g}=(1+tq)^{2g},$$
which checks with \eqref{htorus}.

From \eqref{Uzw-fmla}  we have
$$
\frac{U_2(z,w)}{2}=-\frac{1}{2}\h_{(1)}^{2(2g-2)}(z,w)+\h_{(11)}^{2g-2}(z,w)+\h_{(2)}^{2g-2}(z,w).
$$
Combining \eqref{Vzw-Uzw}, \eqref{H-defn} and  \eqref{H-defn-1} 
$$
\hpoln_2(z,w)=-\frac{1}{2}\frac{(z-w)^{4g}}{(z^2-1)(1-w^2)}+\frac{(z^3-w)^{2g}(z-w)^{2g}}
{(z^4-1)(z^2-w^2)}+\frac{(z-w^3)^{2g}
  (z-w)^{2g}}{(z^2-w^2)(1-w^4)}- 
\frac{1}{2}\frac{(z^2-w^2)^{2g}}{(z^2+1)(1+w^2)}.
 $$ 
Substituting $z=\sqrt q$ and $w=\frac{-1}{t\sqrt q}$  we see that
Theorem~\ref{mhp2}, proved in \S \ref{hm2}, is equivalent to
Conjecture~\ref{main-conj} for $n=2$. 

Next we consider the special cases of $g=0,1$. For $g=0$ we have 
$\M_n$ is
a point for $n=1$ and is empty otherwise. Hence
$$
H(\M_n;q,t)=
\left\{
\begin{array}{ll}
1 & \qquad n=1 \\
0 & \qquad \mbox{otherwise}
\end{array}
\right.
$$
and according to the conjecture \eqref{main-conj} we find
$$
\hpolnc_{i,j}^n=
\left\{
\begin{array}{ll}
1 & \qquad n=1,i=j=0 \\
0 & \qquad \mbox{otherwise}
\end{array}
\right.
$$
hence, after replacing $z^2$ by $z$ and $w^2$ by $w$, we should
have
\begin{equation}
 \label{g=0}
\sum_\lambda \frac 1{\prod(z^{a+1}-w^{l})(z^{a}-w^{l+1})} \,T^{|\lambda|} =
\prod_{r>0,s\geq 0} (1-z^{s}w^{-r}\,T)^{-1}.\end{equation}
In fact we can prove this identity.
\begin{theorem}\label{g=0iden}
 The identity \eqref{g=0} is true.
\end{theorem}
\begin{proof}
 We know from \cite[Thm 3.10 (f)]{garsia-haiman} that 
 \begin{equation}
   \label{GH}
\sum_{|\lambda|=n} \frac 1{\prod(w^l-z^{a+1})(z^a-w^{l+1})} = 
\sum_{|\lambda|=n}  \frac{z^{n(\lambda')}w^{n(\lambda)}}
{\prod (1-z^h)(1-w^h)},
 \end{equation}
where $h=a+l+1$ is the hook length. On the other hand we know \cite[
I.3 ex. 2]{Mc} that
$$
s_\lambda(1,x,x^2,\ldots)=\frac{x^{n(\lambda)}}{\prod(1-x^h)}
$$
where $s_\lambda$ is the Schur function and hence
$$
s_\lambda(1,1/x,1/x^2,\ldots)=\frac{(-x)^{|\lambda|}
x^{n(\lambda')}}{\prod(1-x^h)}.
$$
Summing over all $n$ we then find
$$
\sum_{\lambda} \frac 1{\prod(w^l-z^{a+1})(z^a-w^{l+1})} T^{|\lambda|}=
\sum_{\lambda}
s_\lambda(1,z,z^2,\ldots)s_\lambda(T/w,T/w^2,T/w^3,\ldots) 
$$
and by Cauchy's formula \cite[I  (4.3)]{Mc} this equals the right hand
side of \eqref{g=0}.
\end{proof}

Now let us consider the case $g=1$. We have that $\M_n\simeq
\C^\times\times \C^\times$ for all $n$ (see Theorem~\ref{g=1-char-var}). Hence
$$
H(\M_n;q,t)=(1+qt)^2
$$
and according to Conjecture~\ref{main-conj} we  should have
$$
\hpoln_n(z,w)=(z- w)^2, \qquad n\in \Z_{>0}.
$$
Consequently Conjecture~\ref{main-conj} implies 
\begin{conjecture}
 The following identity holds
 \begin{equation}
   \label{g=1-conj}
\sum_\lambda
\prod
\frac{\left(z^{2a+1}-w^{2l+1}\right)^2}
{(z^{2a+2}-w^{2l})(z^{2a}-w^{2l+2})}\,T^{|\lambda|}=
\prod_{n\geq 1}\prod_{r>0}\prod_{s\geq 0}
\frac{(1-z^{2s+1}w^{-2r+1}\,T^n)^2}
{(1-z^{2s}w^{-2r+2}\,T^n)(1-z^{2s+2}w^{-2r}\,T^n)}
 \end{equation}
\end{conjecture}
\begin{remark}
 The conjecture is a purely combinatorial one. The specialization
 $z=\sqrt q, w=1/\sqrt q$ is essentially Euler's identity which we
 already encountered in \eqref{g=1-q-fmla}. We also prove below
in Remark~\ref{g=1-hua} the specialization $z=0 , w=\sqrt{q}$.

 We have checked \eqref{g=1-conj} numerically up to the $T^6$ terms.
 For this it is more convenient to write in its additive form
 \eqref{V-exp}
$$
\sum_\lambda \prod\frac{\left(z^{2a+1}-w^{2l+1}\right)^2}
{(z^{2a+2}-w^{2l})(z^{2a}-w^{2l+2})}\,T^{|\lambda|}=
\exp\left( \sum_{k\geq 1}
 \frac{(z^k-w^k)^2}{(z^{2k}-1)(1-w^{2k})(1-T^k)}\frac{T^k}k\right) 
$$
and check that the coefficient of $T^n$ on both sides (as rational
functions in $z,w$) agree.
\end{remark}

\subsection{Purity conjecture}
\label{purepart}
\begin{theorem}\label{purity}  Let $A_n(q)$ be the number of absolutely indecomposable $g$-tuples of $n$ by $n$  matrices over the 
finite field $\F_q$ modulo conjugation. Then 
 \begin{equation}
   \label{hua}  
 \hpoln_n(0,\sqrt q)= A_n(q)    
 \end{equation}
\end{theorem}
 \begin{proof}
It is immediate to verify that
\begin{equation}
\label{pure-h}
 \h_\lambda(0,\sqrt q) = \frac{q^{(g-1)\langle\lambda, \lambda
     \rangle}}{b_\lambda(1/q)},
\end{equation}
where $b_\lambda(q)=\prod_{i\geq 1}(1-q)\cdots(1-q^{m_i})$  with $m_i$
is the multiplicity of $i$ in $\lambda$. 

It follows that the left hand side of \eqref{zw-main-fmla} for $z=0, \
w=\!\!\sqrt q$ equals the left hand side of Hua's formula
\cite[Theorem 4.9]{hua} for the $S_{\! g}$ quiver. On the right hand
side we get
$$
\prod_{n\geq 1}\prod_{r>0,i\geq -d_n}
(1-q^{-(r+i)}\,T^n)^{-\hpolnc_{2i,0}^n}
$$
(note that $\hpolnc_{i,0}=0$ for $i$ odd thanks to \eqref{H-symm}).
By the formalism of \eqref{different} we may rewrite this as 
$$
\prod_{n\geq 1}\prod_{r\geq 0,i\leq d_n}
(1-q^{r+i}\,T^n)^{\hpolnc_{-2i,0}^n}.
$$ 
Comparing with the right hand side of Hua's formula we deduce that
$\hpolnc_{-2i,0}=t^n_i$ proving our claim.
 \end{proof}
  \begin{remark} Combining Theorem~\ref{purity} and Conjecture~\ref{main-conj} is what we call the {\em purity conjecture}: the pure part
 of the mixed Hodge polynomial of the character variety $\M_n$ is
 the reverse of the $A$-polynomial of the quiver $S_{\! g}$ (a vertex with $g$    loops)  with dimension vector $n$.
   By a result of Kac \cite{kac} $A_n(q)$, and therefore also
   $\hpoln_n(0,\sqrt q)$, is a polynomial in $q$, which is implied
   by part (1) of Conjecture \ref{comb-conj}. Then Part (2) of Conjecture \ref{comb-conj}
   implies non-negativity of the coefficients of $A_n(q)$, which is
   conjecture 2 of Kac \cite{kac} for the $S_{\! g}$ quiver with dimension $n$ at the vertex. Since $S_{\! g}$ with this
   dimension vector is divisible (for $n>1$) the conjecture is still
   open (the indivisible case was proved in \cite{crawley-boevey-etal} ). To summarize this discussion we can claim: 
    the purity conjecture
   implies Kac's  \cite[Conjecture 2]{kac} for the quiver $S_{\! g}$. 
   In \cite{hausel5} a detailed discussion, motivation and the origin  for this and more general purity conjectures will be given. 
    \end{remark}
 \begin{remark}
   \label{g=1-hua}
   For $g=0$ Theorem~\ref{purity} implies that Hua's formula
     \cite[(5.1)]{hua} is the specialization of \eqref{g=0} at $z=0$. 
   On the other
   hand, for $g=1$ the theorem shows that Hua's formula \cite[(5.2)]{hua} is the
   specialization $z=0, w=\sqrt q$ of our conjecture
   \eqref{g=1-conj}.
 \end{remark}
\begin{proposition}\label{nomidpure}
For all $n,g>0$ we have that  $q^{(g-1)n+1}$ divides
 $A_n(q)$  and 
$$
\left. \frac{A_n(q)}{q^{(g-1)n+1}}\right|_{q=0} = 1
$$
\end{proposition}
\begin{proof}
 This is a consequence of Hua's formula but for convenience we will
 express the result in our notation using \eqref{hua}.  For a
 partition $\lambda$ of $n>0$ we have from  \eqref{pure-h}
$$
\h_\lambda(0,\sqrt q)=\frac{q^{(g-1)\langle\lambda, \lambda
   \rangle+m}}{\prod_{i\geq 1}
     (q-1)(q^2-1) \cdots (q^{m_i}-1)}
$$
where 
$$
m:=\sum_{i\geq 1} \binom{m_i+1} 2
$$
with $m_i=m_i(\lambda)$ the multiplicity of $i$ in $\lambda$.
Among all partitions $\lambda$ of $n$ the exponent
$(g-1)\langle\lambda, \lambda \rangle+m$ of $q$ takes its minimum
value $(g-1)n+1$ only for $\lambda=(n)$. In particular, $q^{(g-1)n+1}$ divides
$\h_\lambda(0,\sqrt q)$ and, moreover, 
$$
\left. \frac{\h_{(n)}(0,\sqrt q)}{q^{(g-1)n+1}}\right|_{q=0}=-1.  
$$
After some calculation we find that
$$
\left. \frac{V_n(0,\sqrt
 q)}{q^{(g-1)n+1}}\right|_{q=0} = -1,
$$
which combined with \eqref{H-defn-1} proves our claim.
\end{proof}
\begin{remark} One consequence of Proposition~\ref{nomidpure}
is that the purity conjecture or more generally our main 
Conjecture~\ref{main-conj} implies, that the largest non-trivial
degree of $PH^*(\M_n)$ is $2(g-1)n(n-1)$. 
Interestingly, \cite[Theorem 7 and Proposition 9]{earl-kirwan} proves
the same about the "pure ring" of the twisted $U(n)$-character variety $\Nn^d$, i.e., the subring generated by the classes $\beta_k$. This and the known situation
for $n=2$ (see the next section) indicates that the "pure ring" of $\Nn^d$ maybe isomorphic with $PH^*(\M_n)$.  An interesting consequence of this would be  that the "pure ring" of $\Nn^d$  is  independent of $d$, unlike the whole  
cohomology $H^*(\Nn^d)$, which does depend on $d$.  Finally, combining the reasoning above with the purity conjecture suggests that the "pure ring"
of $\Nn^d$ could also be used for 
a cohomological interpretation of the $A$ -polynomials $A_n(q)$, 
implying \cite[Conjecture 2]{kac} for the $S_{\!g}$ quiver.
\end{remark}

\begin{remark} \label{an1} If we combine the purity conjecture with 
Remark~\ref{toppure}, we get that the middle cohomology of $\M_n$ should have dimension $A_n(1)$. We  list below the formulas
for the value of $A_n(1)$ for $n=2,3$ and $4$ as a polynomial in $\chi=2g-2$ obtained by computer calculations:

\bes A_2(1)&=&\frac{1}{2}\chi+1\\
A_3(1)&= &\frac{1}{2}\chi^2 + \frac{3}{2}\chi+ 1\\
A_4(1)&=&\frac{2}{3}\chi^3 + \frac{5}{2}\chi^2 + \frac{17}{6}\chi + 1\\
\ees 
It is known that the middle cohomology of $\M_2$ has dimension
$g=\frac{1}{2}(2g-2)+1$ by \cite{hitchin} and that the middle Betti
number of $\M_3$ is $2g^2-g=\frac{1}{2}(2g-2)^2 + \frac{3}{2}(2g-2)+
1$ dimensional \cite{gothen}. For $n\geq 4$ the middle Betti number of
$\M_n$ is not known. However one can say something about the leading
coefficient of $A_n(1)$ as a polynomial in $\chi$.  In the above
formulas it is $\frac{n^{n-3}}{(n-1)!}$ (a proof of this fact will
appear elsewhere). One can also guess the
leading coefficient of $\dim H^{d_n}(\M_n)$ as a function of
$\chi=2g-2$. The dimension $\dim H^{d_n}(\M_n)$ is exactly the number
of fixed point components of the natural circle-action on the
corresponding moduli space of Higgs bundles. These fixed point
components are not well understood in general, but one class of fixed
point components, the so-called type $(1,1,\dots,1)$ is well
understood (see the proof of \cite[Proposition 10.1]{HT4}). Their
number turns out to be a degree $n-1$ polynomial in $\chi$ with
leading coefficient $\frac{n^{n-3}}{(n-1)!}$, which is the volume of a
certain skew hypercube, given by inequalities dictated by the
stability condition for Higgs bundles of type $(1,1,\dots,1)$, which
appear in the proof of \cite[Proposition 10.1]{HT4}.  As the rest of
the fixed point components are expected to be counted by a polynomial
in $\chi$ of degree less then $n-1$, the quantity
$\frac{n^{n-3}}{(n-1)!}$ should be the leading coefficient of $\dim
H^{d_n}(\M_n)$, in agreement with the prediction coming from the above
conjecture.

\end{remark}

\subsection{Intersection form}

Another consequence of Conjecture~\ref{main-conj}  and Proposition~\ref{nomidpure} is that the $\tilde{d}_n$ (=middle) dimensional cohomology of $\tM_n$ does not have pure part. 
This implies the following
\begin{corollary}\label{l2} Conjecture~\ref{main} implies that the intersection form on $H_c^{\tilde{d}_n}(\tM_n)$
is trivial. Equivalently the forgetful map $H^*_c(\tM_n)\to H^*(\tM_n)$
is $0$.\end{corollary}
\begin{proof} Conjecture~\ref{main} and  Proposition~\ref{nomidpure}
implies that there is no pure part in $H^{\tdn}(\tM_n)$, consequently all the non-trivial weights in the weight filtration are $>\tdn$. Now
Poincar\'e duality \eqref{poincareduality} implies that $H^{\tdn}_c(\tM_n)$ has no pure part either; consequently all the non-trivial weights in the weight filtration $<\tdn$. However Theorem~\ref{mixedhodgec}.\ref{forget} shows that the map
$H^{\tdn}_c(\tM_n)\to H^{\tdn}(\tM)_n$ preserves the 
weight filtration.  This proves that the map has to be $0$.
\end{proof}

\section{Mixed Hodge polynomial of $\M_2$}
\label{mixedhodge2}
\subsection{Cohomology ring of $\M_2$}

Here we compute $H(\M_2;q,t)$ by using the explicit description of the ring $H^*(\M_2)$ given in \cite{hausel-thaddeus-generators, hausel-thaddeus-relations}. 
According to \cite{hausel-thaddeus-generators}  the cohomology ring $H^*(\M_2)$ is generated by 
classes $\epsilon_i\in H^1(\M_2)$,
$\psi_i\in H^3(\M_2)$ for $i=1,\dots,2g$ and $\alpha\in H^2(\M_2)$ and $\beta\in H^4(\M_2)$. In the notation of \eqref{defineuniversal} $\alpha=\alpha_2$, $\psi_j=\psi_{2,j}$ and $\beta=\beta_2$. The paper 
\cite{hausel-thaddeus-relations} then proceeds by determining the relations in these universal generators.
The result is as follows.  

Let $\Gamma$ be the group $\rm{Sp}(2g,\Z)$. Let
$\La^k (\psi)$ be the $k$th exterior power of the standard
representation of $\Gamma$, with basis $\psi_1, \dots, \psi_{2g}$.  Define the
{\em primitive part} $\La^k_0(\psi)$ to be the kernel of the natural
map $\La^k(\psi) \to \La^{2g+2-k}(\psi)$ given by the wedge product with
$\ga^{g+1-k}$, where $\ga=2\sum_{i=1}^g \psi_i \psi_{i+g}$ .  The primitive part is complementary to $\ga \La^{k-2}(\psi)
\subset \La^k(\psi)$, and is an irreducible representation of $\Ga$: this is
well-known for $\rm{Sp}(2g, \C)$, and so remains true for the Zariski
dense subgroup $\Ga$. Consequently, \beq \label{dim} \dim (\La^k_0(\psi))=\bino{2g}{k}-\bino{2g}{k-2}\eeq
For any $g, n \geq 0$, let $I^g_n$ be the ideal within the polynomial
ring $\Q[\al,\be,\ga]$ generated by $\ga^{g+1}$ and the polynomials
\beq
\label{mmm}
\rho^{n,g}_{r,s,t}
= \sum_{i=0}^{\min(r,s,g-t)} \, (c-i)! \,
\frac{\al^{r-i}}{(r-i)!} \,
\frac{\be^{s-i}}{(s-i)!} \,\frac{(2 \ga)^{t+i}}{i!},
\eeq
where $c = r+3s+2t-2g+2-n$,
for all $r,s,t \geq 0$ such that
\beq
\label{ppp} t\leq g, \phantom{xxx} 
r+3s+3t > 3g-3+n \mbox{\phantom{xxx} and \phantom{xxx}}
r+2s+2t \geq 2g-2+n.
\eeq

The following is then the main result of \cite{hausel-thaddeus-relations}.

\begin{theorem}
\label{ee}
As a $\Ga$-algebra,
$$H^*(\M_2) = \La(\epsilon)\otimes \left(\bigoplus_{k=0}^g
\La^k_0(\psi)\otimes\Q[\al, \be, \ga]/I^{g-k}_{k}\right).$$
\end{theorem}
\begin{proof} There is a slight difference in the classes $\rho^{n,g}_{r,s,t}$
in \eqref{mmm} and the classes $\rho^c_{r,s,t}$ in \cite{hausel-thaddeus-relations}. In \cite{hausel-thaddeus-relations} the sum for $i$ is between $0$ and  $\min(r,s,c)$. 
First of all $c$ is unnecessary in the $\min$ because $s\leq c$ by
the third inequality in \eqref{ppp}. Second difference is that in 
\eqref{mmm}
we have the sum going from $0$ to $\min(r,s,g-t)$. So the relations
are slightly different, here any monomial which is divisible by 
$\gamma^{g+1}$ is left out. However as $\gamma^{g+1}\in I^g_n$ the two sets of polynomials generate the same ideal. 
\end{proof}

\begin{remark}\label{monomial}
We note that the $I^{g-k}_k$ has the following additive basis: take
all classes $\rho^{n,g}_{r,s,t}$ satisfying \eqref{ppp} and monomials
of the form $\alpha^r\beta^s\gamma^t$ with $t>g$. It is an additive
basis because their leading terms in the lexicographical ordering
additively generate an ideal. 

For the calculation of the mixed Hodge polynomial of $\M_2$ we only need to know 
that a monomial basis for $\Q[\al, \be, \ga]/I^{g-k}_{k}$ is given by $\al^r\be^s\ga^t$, 
for \beq \label{cond}  
{0\leq r\, , 0\leq s\,, 0\leq t \leq g'}  \mbox{\phantom{xx} and \phantom{xx}}
(\, \, r+3s+3t \leq 3g'-3+k \mbox{\phantom{xx} or \phantom{xx}}
r+2s+2t <2g'-2+k\, \, ), 
\eeq
where $g'=g-k$.
\end{remark} 

\subsection{Calculation of the Mixed Hodge polynomial}
\label{hm2}

We introduce the notation $S^{g'}_k$ for the set of triples $(r,s,t)$ of non-negative integers
satisfying  \eqref{cond}. To simplify notation we will use $g$ instead of $g^\prime$ below. 

\begin{lemma} \beq\label{lemma} \sum_{(r,s,t)\in S^{g}_k} a^rb^sc^t&=& {\frac {1-{c}^{g+1}}{ \left( 1-a \right)  \left( 1-b \right)  \left( 1
-c \right) }}-\frac{{a}^{k-2} {b}^{g} \left( 1-{\frac {{c}^{g+1}}{{b}^{g+1}}}
\right) }{ \left( 1-a \right) \left( 1-{\frac {c}{b}}
\right) \left( 1-{\frac {b}{{a}^{2}}} \right) }- \frac{\left( {b}^{g+\left[(k+1)/{2}\right]-1}+a{b}^{g+\left[{k}/{2}\right]-1}
\right)  \left( 1-{\frac {{c}^{g+1}}{{b}^{g+1}}} \right) }{ 
\left( 1-b \right)  \left( 1-{\frac {c}{b}} \right) 
\left( 1-{\frac {{a}^{2}}{b}} \right) }- \nonumber
\\ && \frac{{a}^{3\,g+k-2} \left( 
1-{\frac {{c}^{g}}{{a}^{3\,g}}} \right)}{  \left( 1-a \right) 
\left( 1-{\frac {c}{{a}^{3}}} \right)  \left( 1-{\frac {b}{{a}^{
3}}} \right) } +\frac{{a}^{k-2} {b}^{g} \left( 1-{\frac {{c}^{g}}{{b}^{g}}}
\right) }{\left( 1-a \right) \left( 1-{\frac {c}{b}}
\right)  \left( 1-{\frac {b}{{a}^{3}}} \right) }
\eeq
\end{lemma}

\begin{proof} 
Fix $g$. It is clear that $S_k^{g}\subset S_{k+1}^{g}$. Furthermore
we can separate $S_{k+1}^{g}\setminus S_k^{g}=R^k_1\coprod R^k_2$ 
into the following two disjoint sets:
$$R^k_1:=\{ (r,s,t)\in \Z_{\geq 0}^{3}\, |\, r+3s+3t=3g-3+k+1 \mbox{  and  } r+2s+2t > 2g-2 +k \mbox{ and } t\leq g \}$$ 
$$R^k_2:=\{ (r,s,t)\in \Z_{\geq 0}^{3}\, |\, r+3s+3t\geq 3g-3+k+1 \mbox{  and  } r+2s+2t = 2g-2 +k \mbox{ and } t\leq g \}.$$ 
We can calculate \beq \sum_{(r,s,t)\in R^k_1} a^rb^sc^t&=&\sum_{t=0}^{g-1}\sum_{s=0}^{g-1-t} a^{3g-3+k+1} (b/a^3)^{s}(c/a^3)^t=\sum_{t=0}^{g-1}
a^{3g-3+k+1}\frac{1-(b/a^3)^{g-t}}{1-b/a^3}  (c/a^3)^t\nonumber \\ &=&\frac{a^{3g-3+k+1}(1-(c/a^3)^{g})}{(1-c/a^3)(1-b/a^3)}-\frac{a^{k-2}b^{g}(1-(c/b)^{g})}{(1-b/a^3)(1-c/b)}\label{s1}\eeq
and
\beq \sum_{(r,s,t)\in R^k_2} a^rb^sc^t  &=& \sum_{t=0}^g\sum_{s=g-t}^{g-1+[k/2]-t} a^{2g-2+k}(b/a^2)^s(c/a^2)^t=\sum_{t=0}^g a^{2g-2+k} (c/a^2)^t
(b/a^2)^{g-t} \frac{1-(b/a^2)^{[k/2]}}{1-(b/a^2)} \nonumber \\  &=& \frac{a^{k-2}b^g(1-(c/b)^{g+1})\left(1-(b/a^2)^{[k/2]}\right)}{(1-c/b)(1-b/a^2)},\label{s2}\eeq
thus \beq \label{kul} \sum_{(r,s,t)\in S_{k+1}^{g}\setminus S_k^{g}}a^rb^sc^t=
\sum_{(r,s,t) \in R^k_1} a^rb^sc^t +\sum_{(r,s,t) \in R^k_2} a^rb^sc^t\eeq
As $$\cup_{k'=k}^\infty S^{g}_{k'} = \{(r,s,t)\in \Z_{\geq 0}^3 | t\leq g\} $$ we can deduce that 
\bes {\frac {1-{c}^{g+1}}{ \left( 1-a \right)  \left( 1-b \right)  \left( 1
-c \right) }}&=& \sum_{(r,s,t)\in \cup_{k'=k}^\infty S^{g}_{k'}} a^rb^sc^t= \sum_{(r,s,t)\in S^{g}_k} a^rb^sc^t + \sum_{k'=k}^{\infty}\sum_{(r,s,t)\in S_{k'+1}^{g}\setminus S_{k'}^{g}}a^rb^sc^t.
\ees
Using \eqref{kul} and summing up \eqref{s1} and \eqref{s2} proves the 
Lemma.  \end{proof}

We can now prove Theorem~\ref{mhp2}.

\begin{proof} By Proposition~\ref{universalweights} we know that the classes $\alpha$,
$\psi_k$ and $\beta$ have homogenous weight $2$. 
Thus $\gamma$ has homogenous weight $4$. 
As the cup product is compatible with mixed Hodge structures
by  Theorem~\ref{mixedhodgeproperties}.\ref{preservescup} the homogenous weights of a monomial in the universal generators will be
the sum of the homogenous weights of the factors (see Remark~\ref{weightremark}). Thus a method
to calculate the mixed Hodge polynomial of $\M_2$ is to take
the monomial basis of $\Q[\alpha,\beta,\gamma]/I^{g-k}_k$ 
given in \eqref{cond} evaluate the homogenous weights
of the individual monomials and sum this up over all monomials.

First we have
\begin{lemma}\label{klemma}\begin{multline*} \sum_{(r,s,t)\in S^{g}_k} (q^2t^2)^r(q^2t^4)^s(q^4t^6)^t= \frac {q^{2\,g-2}{t}^{4\,g-4-2\,k} \left( 1-{q}^{4\,g-4\,k+4}{t}^{2\,g-2\,k+2} \right) }{ \left( 1-{q}^{4}{t}^{2} \right)  \left( {q}^{2}
-1 \right)  \left( {q}^{2}{t}^{2}-1 \right) } +\frac {1-{q}^{4\,g-4\,
k+4}{t}^{6\,g-6\,k+6}}{ \left( 1-{q}^{4}{t}^{6} \right)  \left( {q}^{2
}{t}^{2}-1 \right)  \left( {q}^{2}{t}^{4}-1 \right) }
\\ -\frac{1}{2}{\frac {{q}^{2\,g-2-k}{t}^{4\,g-4-2\,k} 
\left( 1-{q}^{2\,g-2\,k+2}{t}^{2\,g-2\,k
+2} \right) }{ \left( 1-{q}^{2}{t}^{2} \right)  \left( q-1 \right) 
\left( q{t}^{2}-1 \right) }}-\frac{1}{2}\,{\frac {\,(-q)^{2\,g-2-k}
{t}^{4\,g-4-2\,k} \left( 1-{q}^{2\,g-2\,k+2}{t}^{2\,g-2\,k+2} \right) 
}{ \left( 1-{q}^{2}{t}^{2} \right)  \left( q+1 \right)  \left( q{t}^{2
}+1 \right) }}
\end{multline*}
\end{lemma}
\begin{proof} Substitute $a=q^2t^2$, $b=q^2t^4$, $c=q^4t^6$ in 
(\ref{lemma}). To prove that the resulting rational function is the same
as the one above, one can multiply over with the denominators and get an identical expression. 
\end{proof}

We can now use the description of the cohomology ring of $\M_2$ in Theorem~\ref{ee} to get the mixed Hodge polynomial $H(\M_2;q,t)$.
We have 
$$\frac{H(\M_2;q,t)}{(1+qt)^{2g}}=\sum^g_{k=0} \left(\bino{2g}{k}-\bino{2g}{k-2}\right)(q^2t^3)^k
\sum_{(r,s,t)\in S^{g-k}_k} (q^2t^2)^r(q^2t^4)^s(q^4t^6)^t.$$
Writing in Lemma~\ref{klemma} and summing it up with $k$ yields
exactly Theorem~\ref{mhp2}.
\end{proof}
\begin{remark} The above proof of Theorem~\ref{mhp2} follows
closely the geometry behind the proof of Theorem~\ref{ee} in \cite{hausel-thaddeus-relations}. There certain spaces ${\mathcal H}_k$ of rank $2$ Higgs bundles with a pole of order at most $k$ are introduced. It is 
shown there that ${\mathcal H}_0\cong \M_2$ are diffeomorphic, they form a tower: ${\mathcal H}_k\subset {\mathcal H}_{k+1}$ and the direct limit ${\mathcal H}_\infty :=\cup_{k=0} {\mathcal H}_k$ is homotopically equivalent
with the classifying space of a certain gauge group. The cohomology 
ring $H^*({\mathcal H}_k)$ is also generated by the same classes $\epsilon_i,\alpha,\psi_i,\beta$. One can show by the description of their cohomology ring
in \cite{hausel-thaddeus-relations} that there exists an abstract weight filtration
on $H^*({\mathcal H}_k)$ by setting the weight of the universal generators $\alpha,\psi_i,\beta$ be $4$, and the weight of $\epsilon_i$ to be $2$ . Lemma~\ref{lemma} then can be considered as calculating the
natural two-variable polynomial associated to this abstract weight filtration on the $\Gamma$-invariant part of $H^*({\mathcal H}_k)$. Similarly to the calculation above, one can obtain the following formula for the two-variable polynomial associated
to this filtration on the whole cohomology $H^*({\mathcal H}_k)$:

\begin{multline*} \frac{(q^2t^3+1)^{2g}(qt+1)^{2g}}{(q^2t^2-1)(q^2t^4-1)}+
\frac{q^{2g-2}t^{4g-4+2k}(q^2t+1)^{2g}(qt+1)^{2g}}{(q^2-1)(q^2t^2-1)}-\\
-\frac{1}{2}\frac{q^{2g-2+k}t^{4g-4+2k}(qt+1)^{2g}(qt+1)^{2g}}{(qt^2-1)(q-1)}-\frac{1}{2}
\frac{ (-q)^{2g-2+k}t^{4g-4+2k}(qt-1)^{2g}(qt+1)^{2g}}{(q+1)(qt^2+1)}\end{multline*}

This polynomial has some remarkable properties. First we see that in the $k\rightarrow \infty$ limit only the first term survives, which gives
the two-variable rational function associated to this abstract
filtration on the cohomology $H^*({\mathcal H}_\infty)$ of the classifying space
of the gauge group, which is known to be a free anticommutative algebra on the universal generators $\epsilon_i,\alpha,\psi_i,\beta$. 

Second, the polynomial satisfies a curious Poincar\'e duality, when replacing $q$ by $1/qt^2$. Thus when we set $t=-1$ in the above polynomial we have a palindromic polynomial in $q$. It has degree 
$8g-6+2k$. We may therefore expect that there is a character variety
version $\M_2^k$ of the Higgs moduli spaces ${\mathcal H}_k$ so that the abstract
weight filtration we put on $H^*({\mathcal H}_k)$ is the actual weight filtration
coming from the mixed Hodge structure on $H^*(\M_2^k)$. However if 
this was the case the $E$-polynomial of $\M^k_2$ would be of degree $8g-6+2k$, and therefore $\M_2^k$ would have dimension 
$8g-6+2k$. The dimension of ${\mathcal H}_k$ is $8g-6+3k$. Therefore what we
could expect is perhaps a deformation retract of ${\mathcal H}_k$ being diffeomorphic to a certain character variety $\M_2^k$ of dimension $8g-6+2k$ with the above mixed Hodge polynomial.

\end{remark}

\begin{remark} \label{provecurious2} We can now deduce Corollary~\ref{curious2} by combining Theorem~\ref{mhp2} and Remark~\ref{deducecpd}. 
\end{remark} 
\subsection{Curious Hard Lefschetz}
\label{curioushard}
Define the Lefschetz map $L:H^i(\tM_2)\to H^{i+2}(\tM_2)$ by
$x\mapsto \alpha\cup x$, where $\alpha=\alpha_2$ is the universal class in $H^2(\tM_2)$ defined in \eqref{defineuniversal}. As it respects mixed Hodge structures and $\alpha$ has homogenous weight $2$ by Proposition~\ref{universalweights} it defines a map on the graded pieces of the weight filtration 
$L: Gr^W_{l} H^i(\tM_2)\to Gr^W_{l+4}H^{i+2}(\tM_2)$. We now
prove Theorem~\ref{lefschetz}.
\begin{proof}
We start with a few lemmas. Let us call the last (in the ordering of the 
sum in \eqref{mmm}) monomial  $\alpha^{r_0}\beta^{s_0}\gamma^{t_0}$ 
appearing in $\rho^{n,g}_{r,s,t}$ its {\em tail}. Clearly $r_0=r-\min(r,s,g-t)$, 
$s_0=s-\min(r,s,g-t)$ and $t_0=t-\min(r,s,g-t)$. Thus if a monomial $\alpha^{r_0}\beta^{s_0}\gamma^{t_0}$ is the tail of the polynomial
$\rho^{n,g}_{r,s,t}$ then \beq \label{tailcondition} r_0=0 \mbox{ or } s_0=0 \mbox{ or } t_0=g.\eeq 
Let us denote by $T^g_n$ the set of triples $(r,s,t)\in \Z_{\geq 0}^3$ 
which satisfy \eqref{ppp}.
\begin{lemma}\label{det} Let $r_0,s_0,t_0\in \Z_{\geq 0}$, $t_0\leq g$ and satisfying
\eqref{tailcondition}. Denote by $d$ the number of polynomials
$\rho^{n,g}_{r,s,t}$ with  $(r,s,t)\in T^g_n$ and tail $\alpha^{r_0}\beta^{s_0}\gamma^{t_0}$. It  
satisfies \beq\label{dineq} d= \min\left[ t_0+1,\max(r_0+3s_0+4t_0-(3g-3+n),0),\max(r_0+2s_0+3t_0-(2g-3+n),0)\right].\eeq The $d$ times $d$  matrix $A=(a_{ij})_{i,j=0}^{d-1}$  is non-singular, where
$$a_{ij}=\left\{ \begin{array}{cc}   0 & i+j>t_0 \\  \frac{\left(r_0+3s_0+3t_0-2g+2-n+i-j\right)! 2^{t_0-i}}{(r_0+i)!(s_0+i)!(t_0-i-j)!} & i+j\leq t_0 \end{array}\right.,$$ which is the
coefficient of $\alpha^{r_0+i}\beta^{s_0+i}\gamma^{t_0-i}$
in $\rho^{n,g}_{r_0+t_0-j,s_0+t_0-j,j}$. \end{lemma}
\begin{proof} To prove the first statement we need to count
the number of $0\leq i\leq t_0$ such that $(r_0+i,s_0+i,t_0-i)\in T^n_g$ 
consequently  they satisfy $$r_0+3s_0+3t_0+i>3g-3+n,$$ thus
$$t_0\geq i > 3g-3+n-(r_0+3s_0+3t_0).$$ Similarly, we have $$r_0+2s_0+2t_0+i>2g-3+n,$$ which yields $$t_0\geq i > 2g-3+n-(r_0+2s_0+2t_0).$$This proves (\ref{dineq}).

For the second statement consider the matrix $B=(b_{ij})_{i,j=0}^{d-1}$
with $b_{ij}=\frac{(r_0+3s_0+3t_0-2g+2-n+i-j)!}{(t_0-i-j)!},$  if $t_0-i-j\geq 0$, and $b_{ij}=0$ otherwise. As the matrix
$B$ is obtained from $A$ by multiplying the rows and columns by
non-zero constants it is enough to show that $B$ is non-singular. 
Introduce the notation $(a)_j=a(a+1)\dots (a+j-1)$ and $(a)_0=1$. 
Now we can write \bes b_{ij}&=&\frac{(r_0+3s_0+3t_0-2g+2-n+i)!(t_0-i-j+1)_j}{(r_0+3s_0+3t_0-2g+2-n+i-j+1)_j(t_0-i)!}\\ &=&\frac{(r_0+3s_0+3t_0-2g+2-n+i)!(t_0-i-j+1)_j}{(-1)^j(-r_0-3s_0-3t_0+2g-2+n-i)_j(t_0-i)!} ,\ees which is valid for any $i$ and $j$ as $(t_0-i-j+1)_j=0$ if and only if $t_0-i-j< 0$ (note that $$t_0-i\geq t_0-d+1\geq 0$$  and $$r_0+3s_0+3t_0-2g+2-n+i-j+1\geq r_0+3s_0+3t_0-2g+2-n+i-d+2>0$$ by \eqref{dineq}). Now multiplying the rows and columns of $B$ by non-zero constants we get the matrix $C= (c_{ij})_{i,j=0}^{d-1}$
with $$c_{ij}=\frac{(\alpha_i-\beta_j)_j}{(\alpha_i)_j},$$ where $\alpha_i=-r_0-3s_0-3t_0+2g-2+n-i$ and $\beta_j=-r_0-3s_0-4t_0+2g-3+n+j.$
The determinant of a matrix like  $C$ was calculated in \cite[Lemma 19]{gessel-viennot}.
Their formula gives $$|C|=\prod_{i=0}^{d-1}\frac{(\beta_i)_i}{(\alpha_i)_{d-1}}\prod_{0\leq i<j<d} (i-j).$$ Because $$\beta_i+i-1\leq -r_0-3s_0-4t_0+2g-3+n+2d-2-1\leq -r_0-3s_0-4t_0+3g-5+n+d<0$$  
by \eqref{dineq} we get $|C|\neq 0$ and consequently $|A|\neq 0$.
This completes the proof. 
\end{proof}

\begin{lemma} \label{cases} Let $(r_0,s_0,t_0)$ satisfy (\ref{tailcondition}) and let
\beq \label{weight} w=6g-6+2n-(2r_0+2s_0+4t_0).\eeq The number of monomials of the
form $\alpha^{r_0+i}\beta^{s_0+i}\gamma^{t_0-i}$ for which 
$       0\leq r_0+ i<3g-3+n-w$ and $0\leq i \leq t_0$ is at least $d$.
\end{lemma}

\begin{proof}
We distinguish three cases depending on 
which of the cases of (\ref{tailcondition}) is satisfied.

First case is when $r_0=0$.  The number of monomials of the
form $\alpha^{i}\beta^{s_0+i}\gamma^{t_0-i}$ for which 
$       0\leq i<3g-3-n-w$ and $0\leq i \leq t_0$ is clearly \bes \min(t_0+1,3g-3+n-w)&=&\min(t_0+1,3g-3+n-w)\\&=& \min(t_0+1,2r_0+2s_0+4t_0-(3g-3+n))\\&\geq & \min(t_0+1,r_0+2s_0+3t_0-(2g-3+n))\\ &\geq& d \ees
because of  (\ref{weight}), $0\leq r_0$, $t_0\leq g$ and
(\ref{dineq}).

Second case is when $s_0=0$. The number of monomials of the
form $\alpha^{r_0+i}\beta^{i}\gamma^{t_0-i}$ for which 
$       r_0+i<3g-3-n-w$ and $0\leq i \leq t_0$ is clearly \bes \min(t_0+1,3g-3+n-w-r_0)&=&\min(t_0+1,3g-3+n-w-(r_0-s_0))\\&=& \min(t_0+1,2r_0+2s_0+4t_0-(3g-3+n)-(r_0-s_0))\\&=& \min(t_0+1,r+3s_0+4t_0-(3g-3+n))\\ &\geq& d \ees
because of  (\ref{weight}) and
(\ref{dineq}). 

Finally the third case is when $t_0=g$. Now the number of monomials
of the form $\alpha^{r_0+i}\beta^{s_0+i}\gamma^{g-i}$ for which
$       r_0+i<3g-3-n-w$ and $0\leq i \leq t_0$ is clearly
\bes \min(t_0+1,3g-3+n-w-r_0)&=&\min(t_0+1,3g-3+n-w-(r_0+t_0-g))\\
 &=&\min(t_0+1,2r_0+2s_0+4t_0-(3g-3+n)-(r_0+t_0-g))\\ &=&\min(t_0+1,r_0+2s_0+3t_0-(2g-3+n))\\ &\geq& d
\ees
\end{proof}

We say that $x\in \Q[\alpha,\beta,\gamma]$ has homogenous  weight $w=w(x)$ if it is
homogeneous of degree $w$ when $\alpha,\beta,\gamma$ are assigned degrees $2,2$ and $4$ respectively. Note that all the classes
$\rho^{n,g}_{r,s,t}$ have homogenous weight $2r+2s+4t$.

\begin{lemma}\label{lemmalefschetz}  Let $x\in \Q[\alpha,\beta,\gamma]$ have homogenous weight $w<3g-3+n$. Then  $x\alpha^{3g-3+n-w} \in I^g_n$ implies
$x\in I^g_n$. 
\end{lemma}
\begin{proof} By Remark~\ref{monomial} we  can write \beq \label{assume}x\alpha^{3g-3+n-w(x)}=\gamma^{g+1}y+\sum_{(r,s,t) \in T^g_n(w)  } \lambda_{r,s,t} \rho^{n,g}_{r,s,t}, \eeq
where $y\in \Q[\alpha,\beta,\gamma]$, $\lambda_{r,s,t}\in \Q$  and $T^g_n(w)$ is the set of non-negative triples $(r,s,t)$
satisfying (\ref{ppp}) and $w(\rho^{n,g}_{r,s,t})=2r+2s+4t=6g-6+2n-w=w(x\alpha^{3g-3+n-w} )$. 

We show that all $\lambda_{r,s,t}=0$. Take $(r,s,t)\in T^g_n(w)$.
 Let $\alpha^{r_0}\beta^{s_0}\gamma^{t_0}$ be the 
tail of $\rho^{n,g}_{r,s,t}$. In particular $w(\rho^{n,g}_{r,s,t})=2r_0+2s_0+4t_0$. According to Lemma~\ref{det} the number of relations of $(r',s',t')\in T^g_n(w)$ with tail $\alpha^{r_0}\beta^{s_0}\gamma^{t_0}$ is
$d$ given by (\ref{dineq}). On the other hand Lemma~\ref{cases} 
implies that the number of monomials $\alpha^{r_0+i}\beta^{s_0+i}\gamma^{t_0-i}$ such that $    0\leq r_0+ i<3g-3+n-w$ and 
$0\leq i \leq t_0$ is at least $d$. But these monomials do not 
appear on the LHS of (\ref{assume}) because all the terms there are divisible by $\alpha^{3g-3+n-w}$. They are only contained in relations with tail $\alpha^{r_0}\beta^{s_0}\gamma^{t_0}$ therefore
Lemma~\ref{det} implies $\lambda_{r,s,t}=0$.

Consequently $x\alpha^{3g-3+n-w(x)}=\gamma^{g+1}y$, which 
implies $x$ is divisible by $\gamma^{g+1}$, thus $x\in I^g_n$. Lemma~\ref{lemmalefschetz}
follows.
\end{proof}

 We can now prove Theorem~\ref{lefschetz}. By Corollary~\ref{curious2} (proved in Remark~\ref{provecurious2}) we know  that $$\dim(Gr^W_{6g-6-2l} H^{i-l}(\tM_2))=h^{3g-3-l,3g-3-l;i-l}(\tM_2)=h^{3g-3+l,3g-3+l;i+l}(\tM_2)=\dim(Gr^W_{6g-6+2l} H^{i+l}(\tM_2))$$ thus it is enough to 
show that $L^{l}$ in the Theorem above is an injection. 

Now by Theorem~\ref{ee} any element  $z\in Gr^W_{6g-6-2l} H^{i-l}(\tM_2)$
can be represented as $z=\sum_{k=0}^g y_k[x_k]$, where $y_k\in\Lambda^k_0(\psi)$ and $[x_k]\in \Q[\alpha,\beta,\gamma]/I^{g-k}_k$,
with a representative $x_k\in \Q[\alpha,\beta,\gamma]$ of homogenous weight. 
As  $w(y_k)=2k$ we have $w(x_k)= 3g-3-l-2k=3(g-k)-3+k-l$ or equivalently $l=3(g-k)-3+k-w(x_k)$. Assume 
now that $z\alpha^{l}=0$. This would imply that $x_k\alpha^{l}\in I^{g-k}_k$ for all $k$. By the previous Lemma~\ref{lemma} $x_k\in I^{g-k}_k$ and so $z=0$. Theorem~\ref{lefschetz} follows.
\end{proof} 

Corollary~\ref{purecorollary} implies that the  pure part of $H^*(\tM_2)$ is $g$ dimensional,
due to the Newstead relation $\beta^g=0$ proved in \cite{hausel-thaddeus-relations}.   This combined with Theorem~\ref{mhp2} and Theorem~\ref{purity} proves Theorem~\ref{pp2}.   

We know from \cite{hitchin} that the middle cohomology $H^{6g-6}(\tM_2)$ is  also $g$-dimensional. The curious Hard Lefschetz map then gives a natural isomorphism between the associated graded of the weight filtration on the vector spaces $PH^*(\tM_2)$ and $H^{6g-6}(\tM_2)$ (cf. Remark~\ref{toppure} and Remark~\ref{an1}).

\subsection{Intersection form}
 Theorem~\ref{mhp2} and Proposition~\ref{nomidpure} implies
that the middle cohomology of $\tM_2$ does not have pure part and
as explained in Corollary~\ref{l2}  we have

\begin{corollary} \label{corollarym2intersect}
The intersection form on $H^{6g-6}_c(\tM_2)$ is trivial,
i.e., there are no "topological $L^2$ harmonic forms" on $\tM_2$. 
\end{corollary}

\
\section{Appendix by Nicholas M. Katz: E-polynomials, zeta-equivalence, and polynomial-count varieties}
\stepcounter{subsection}
Given a noetherian ring $R$, we denote by $(\Sch/R)$ the category of separated
$R$-schemes of finite type, morphisms being the $R$-morphisms. We
denote by $K_{0}(\Sch/R)$ its Grothendieck group. By definition,
$K_{0}(\Sch/R)$
is the quotient of the
free abelian group on elements $[X]$, one for each separated $R$-scheme
of finite type, by the subgroup generated by all the relation elements
$$[X]-[Y],{\rm\ whenever\ }X^{red}\cong Y^{red},$$
and
$$[X] -[X\setminus Z] -[Z],{\rm\  whenever\ }Z\subset X{\rm\ is\ a\ closed\ subscheme}. $$

It follows easily that if $X$ is a finite union of locally closed
subschemes $Z_{i}$, then in $K_{0}(\Sch/R)$ we have the
inclusion-exclusion relation
$$[X] = \sum_{i}[Z_{i}] - \sum_{i<j}[Z_{i}\cap Z_{j}] +\ldots.$$

For any ring homomorphism $R \longrightarrow R'$ of noetherian rings,
the ``extension of scalars'' morphism from $(\Sch/R)$ to $(\Sch/R')$
which sends $X/R$ to $X\otimes_{R}R'/R'$, extends to a group
homomorphism from $K_{0}(\Sch/R)$ to $K_{0}(\Sch/R')$.

Suppose $A$ is an abelian group, and $\rho$ is an ``additive
function'' from $(\Sch/R)$ to $A$, i.e., a rule which assigns to
each $X\in (\Sch/R)$ an element $\rho(X) \in A$, such that $\rho(X)$
depends only on the isomorphism class of $X^{red}$, and such that
whenever $Z\subset X$ is a closed subscheme, we have
$$\rho(X)=\rho(X-Z)+\rho(Z).$$
Then $\rho$ extends uniquely to a group homomorphism from
$K_{0}(\Sch/R)$
to $A$, by defining $\rho(\sum_{i} [X_{i}])=\sum_{i}\rho(X_{i})$.

When $R$=$\C$, we have the following simple lemma, which we record now
for later use.
\begin{lemma}\label{K-diff}
Every element of $K_{0}(\Sch/\C)$ is of
   the form $[S]-[T]$, with $S$ and $T$ both projective smooth
   (but not necessarily connected) $\C$-schemes. 
\end{lemma}
\begin{proof}
   To show this, we argue as
   follows. It is enough to show that for
   any separated $\C$-scheme of finite type $X$, $[X]$ is of this type.
   For then $-[X]=[T]-[S]$, and  
   $$[S_{1}]-[T_{1}] +[S_{2}]-[T_{2}] = [S_{1}\sqcup
   S_{2}]-[T_{1}\sqcup 
   T_{2}],$$
   and the disjoint union of two projective smooth schemes is again
   one. [Indeed, if we embed each in a large projective space, 
   say $S_{i} \subset \bfP^{N_{i}}$ and pick a point $a_{i}\in \bfP^{N_{i}}\setminus
   S_{i}$, then $S_{1}\times a_{2}$ and $a_{1}\times S_{2}$ are
   disjoint in $\bfP^{N_{1}}\times\bfP^{N_{2}}$.]

   We first remark that for any $X$ as above, $[X]$ is of the
   form $[V]-[W]$ with $V$ and $W$ affine. This follows from
   inclusion-exclusion by taking a finite covering of $X$ by affine
   open sets, and noting that the disjoint union of two affine
   schemes of finite type is again an affine scheme of finite type.
   So it suffices to prove our 
   claim for affine $X$. Embedding $X$ as a closed subscheme of some
   affine space $\A^{N}$ and using the relation
   $$[X] = [\A^{N}]-[\A^{N} \setminus X],$$
   it now suffices to prove our claim for smooth quasiaffine $X$.
   By resolution, we can find a projective smooth compactification $Z$
of $X$, such that $Z\setminus X$ is a union of smooth divisors
$D_{i}$ in $Z$ with normal crossings. Then by inclusion-exclusion we
have
$$[X] = [Z] - \sum_i [D_i] 
       +\sum_{i,j} [D_i \cap D_j] + \ldots.$$
In this expression, each summand on the right hand side is 
projective and smooth. Taking for $S$ the disjoint union of the summands with a
plus sign and for $T$ the disjoint union of the summands with a
minus sign, we get the desired expression of our $[X]$ as $[S]-[T]$,
with $S$ and $T$ both projective and smooth. 
\end{proof}

Now take for $R$ a finite field $\F_{q}$. For each integer $n \ge 1$,
the function on
$(\Sch/\F_{q})$ given by $X\mapsto \#X(\F_{q^{n}})$ is visibly an
additive function from $(\Sch/\F_{q})$ to $\Z$. Its extension to $K_{0}(\Sch/\F_{q})$ will be
denoted
$$\gamma \mapsto \#\gamma(\F_{q^{n}}).$$
We can also put all these  functions together, to form the zeta function. Recall 
that the zeta function
$Z({X/\F_{q}},t)$ of 
$X/\F_{q}$ is the power series (in fact it is a rational function)
defined by $$Z(X/\F_{q},t)=exp\left(\sum_{n\ge 1}
\#X(\F_{q^{n}})t^{n}/n\right).$$ Then $X\mapsto
Z(X/\F_{q},t)$ is an additive function with values in
the multiplicative group $\Q(t)^{\times}$.
We denote by
$$\gamma \mapsto \Zeta(\gamma/\F_{q},t)$$
its extension to $K_{0}(\Sch/\F_{q})$.
We say that an element $\gamma \in K_{0}(\Sch/\F_{q})$ is {\it
zeta-trivial} if $\Zeta(\gamma/\F_{q},t)=1$, i.e., if
$\#\gamma(\F_{q^{n}})=0$ for all $n \ge 1$. We say that two elements
of $K_{0}((\Sch/\F_{q})$ are {\it zeta-equivalent}  if they have the same 
zeta functions, i.e., if their difference 
is zeta-trivial.

We say that an element $\gamma \in K_{0}(\Sch/\F_{q})$
is {\em polynomial-count} (or has {\it polynomial count}) if there exists a (necessarily
unique) polynomial $P_{\gamma/\F_{q}}(t) =\sum_i a_it^i \in \C[t]$ 
such that for every
finite extension $\F_{q^{n}}/\F_{q}$, we have
$$\#\gamma(\F_{q^{n}})=P_{\gamma/\F_{q}}(q^{n}).$$
If $\gamma/\F_{q}$ has polynomial count, its counting
polynomial $P_{\gamma /\F_{q}}(t)$ lies in $\Z[t]$. 
(To see this, we argue as follows. On the one hand, from the series
definition of the zeta function, and the polynomial formula for the
number of rational points, we have
$$(td/dt)log(Z(\gamma/\F_{q},t))=\sum_{i}
a_{i}q^{i}t/(1-q^{i}t).$$
As the zeta function is a rational function, say
$\prod_{i}(1-\alpha_{i}t)/\prod_{j}(1-\beta_{j}t)$ in lowest
terms, we first see by comparing logarithmic derivatives that each of 
its zeroes and poles is a non-negative power of $1/q$.
Thus for some integers $b_{n}$, the zeta function is of the form
$\prod_{n\ge 0}(1 - q^{n}t)^{-b_{n}}.$ Again comparing logarithmic
derivatives, we see that we have $a_{n}=b_{n}$ for each $n$.)

Equivalently, an element  $\gamma \in K_{0}(\Sch/\F_{q})$
is polynomial-count if it is zeta-equivalent to a $\Z$-linear combination of
classes of affine spaces $[\A^{i}]$, or, equivalently, to a 
$\Z$-linear combination of
classes of projective spaces $[\bfP^{i}]$ (since
$[\A^{i}]=[\bfP^{i}]-[\bfP^{i-1}]$, with the convention that
$\bfP^{-1}$ is the empty scheme).
If $\gamma/\F_{q}$ is polynomial-count, then so is its extension of
scalars from $\F_{q}$ to any finite extension field, with the 
{\it same} counting polynomial. (But an element $\gamma/\F_{q}$ which is not
polynomial-count can become polynomial-count after extension of scalars, e.g.,
a nonsplit torus over $\F_{q}$, or, even more simply, the zero locus 
of a square-free polynomial $f(z)\in \F_{q}[z]$ which does not
factor completely over $\F_{q}$.)

Now let $R$ be a ring which is finitely generated as a $\Z$-algebra. We say that an element
$\gamma \in K_{0}(\Sch/R)$ is zeta-trivial if, for
every finite field $k$, and for every ring homomorphism $\phi : R
\longrightarrow k$, the element $\gamma_{\phi,k}/k$ in $K_{0}(\Sch/\F_{q})$
deduced from $\gamma$ by
extension of scalars is zeta-trivial. And we say that two elements are
{\it zeta-equivalent} if their difference is zeta-trivial.

We say that an element
$\gamma \in K_{0}(\Sch/R)$ is {\it strongly polynomial-count} with (necessarily
unique) counting polynomial $P_{\gamma/R}(t) \in \Z[t]$ if,
for
every finite field $k$, and for every ring homomorphism $\phi : R
\longrightarrow k$, the element $\gamma_{\phi,k}/k$ in $K_{0}(\Sch/\F_{q})$
deduced from $\gamma$ by
extension of scalars is polynomial-count with counting polynomial
$P_{\gamma/R}(t)$.

We say that an element
$\gamma \in K_{0}(\Sch/R)$ is {\it fibrewise polynomial-count} if, 
for every ring homomorphism $\phi : R
\longrightarrow k$, the element $\gamma_{\phi,k}/k$ in $K_{0}(\Sch/\F_{q})$
deduced from $\gamma$ by
extension of scalars is polynomial-count (but we allow its counting
polynomial to vary with the choice of $(k,\phi)$).

All of these notions, zeta-triviality, zeta equivalence, being
strongly or fibrewise polynomial-count, are stable by extension of
scalars of finitely generated rings.

We now pass to the complex numbers $\C$. Given an element
$\gamma \in K_{0}(\Sch/\C)$, by a ``spreading out'' of $\gamma/\C$, we
mean an element $\gamma_{R} \in K_{0}(\Sch/R)$, $R$ a subring 
of $\C$ which is finitely generated as a $\Z$-algebra, which gives
back $\gamma/\C$ after extension of scalars from $R$ to $\C$.
It is standard that such spreadings out exist, and that given two
spreadings out $\gamma_{R} \in K_{0}(\Sch/R)$ and $\gamma_{R'} \in
K_{0}(\Sch/R')$, then over some larger finitely generated ring $R''$
containing both $R$ and $R'$, the two spreadings out will agree in
$K_{0}(\Sch/R'')$.

We say that an element $\gamma \in K_{0}(\Sch/\C)$ is zeta-trivial if
it admits a spreading out $\gamma_{R} \in K_{0}(\Sch/R)$ which is
zeta-trivial. One sees easily, by taking spreadings out to a common
$R$, that the zeta-trivial elements form a subgroup of $K_{0}(\Sch/\C)$. 

We say that two elements are zeta-equivalent if their 
difference is zeta-trivial. We say that an element is strongly
polynomial-count, with counting polynomial $P_{\gamma}(t) \in \Z[t]$,
(respectively fibrewise polynomial-count)
if it admits a spreading out which has this property. 

Given $X/\C$ a separated scheme of finite type, its {\it E-polynomial}
$E(X; x,y) \in \Z[x,y]$ is defined as follows. The compact cohomology
groups $H_{c}^{i}(X^{an},\Q)$ carry Deligne's mixed Hodge
structure, cf. \cite{De-Hodge II} and  \cite[8.3.8]{De-Hodge III}, and one defines
$$E(X;x,y)=\sum_{p,q}e_{p,q}x^{p}y^{q},$$
where the coefficients $e_{p,q}$ are the virtual Hodge
numbers, defined in terms of the pure Hodge structures which are the associated gradeds 
for the weight filtration on the compact cohomology as follows:
$$e_{p,q} :=
\sum_{i}(-1)^{i}h^{p,q}(gr^{p+q}_{W}(H_{c}^{i}(X^{an},\C))).$$
Notice that the value of $E(X;x,y)$ at the point $(1,1)$ is just the
(compact, or ordinary, they are equal, by \cite{Lau}) Euler
characteristic of $X$.
One knows that the formation of the E-polynomial is additive (because
the excision long exact sequence is an exact sequence in the abelian 
category of mixed Hodge structures, cf. \cite[8.3.9]{De-Hodge III}). 
So we can speak of the
E-polynomial $E(\gamma;x,y)$ attached to an element $\gamma \in K_{0}(\Sch/\C)$.

\begin{theorem}We have the following results. \label{katz}
   \begin{itemize}
\item[(1)]  If $\gamma \in K_{0}(\Sch/\C)$ is zeta-trivial, then 
$$E(\gamma;x,y) = 0.$$
\item[(2)]  If $\gamma_{1}\in K_{0}(\Sch/\C)$ and 
$\gamma_{2}\in K_{0}(\Sch/\C)$ are zeta-equivalent, then 
$$E(\gamma_{1};x,y) = E(\gamma_{2};x,y).$$
In particular, if $X$ and $Y$ in $(\Sch/\C)$ are zeta-equivalent, then
$$E(X;x,y) = E(Y;x,y).$$
\item[(3)]  If $\gamma \in K_{0}(\Sch/\C)$ is strongly
polynomial-count, with counting polynomial $P_{\gamma}(t) \in \Z[t]$, 
then
$$E(\gamma;x,y) = P_{\gamma}(xy).$$
In particular, if $X \in (\Sch/\C)$ is strongly polynomial-count, with
counting polynomial $P_{X}(t) \in \Z[t]$, 
then
$$E(X;x,y) = P_{X}(xy).$$
\end{itemize}
\end{theorem}

\begin{proof}Assertion (2) is an immediate consequence of (1), by the 
   additivity of the E-polynomial. Statement (3) results from (2) as 
   follows. If $\gamma \in K_{0}(\Sch/\C)$ is strongly
polynomial-count, with counting polynomial
$P_{\gamma}(t)=\sum_{i}a_{i}t^{i} \in \Z[t]$, then by definition $\gamma$ is
zeta-equivalent to  $\sum_{i}a_{i}[\A^{i}]\in  K_{0}(\Sch/\C)$.
So we are reduced to noting that $E(\A^{i};x,y) =x^{i}y^{i}$,
which one sees by writing $[\A^{i}]=[\bfP^{i}]-[\bfP^{i-1}]$ and using
the basic standard fact that $E(\bfP^{i};x,y) =\sum_{0\le j \le i}x^{j}y^{j}$.
So it remains only to prove assertion (1) of the theorem. By
lemma~\ref{K-diff}, every element $\gamma \in K_{0}(\Sch/\C)$ 
is of the form $[X]-[Y]$, with $X$ and $Y$ are projective smooth $\C$-schemes.
So assertion (1)
results from the following theorem, which is proven, but not quite
stated, in \cite{Wang}. (What Wang proves is that ``K-equivalent''
projective smooth connected $\C$-schemes have the same Hodge numbers, 
through the intermediary of using motivic integration to show that K-equivalent
projective smooth connected $\C$-schemes are zeta-equivalent.)
   \end{proof}

\begin{theorem}Suppose $X$ and $Y$ are projective smooth $\C$-schemes
   which are zeta-equivalent. Then $$E(X;x,y)=E(Y;x,y).$$
   \end{theorem}

\begin{proof}   
Pick spreadings out
$\X/R$ and $\mathcal{Y}/R$ over a common $R$ which are 
zeta-equivalent. At the expense of 
inverting some nonzero element in $R$, we may further assume that
both $\X/R$ and $\mathcal{Y}/R$ are projective and smooth, and that $R$ is
smooth over $\Z$. We denote the structural morphisms of $\X/R$ and $\mathcal{Y}/R$
by
$$ f: \X \longrightarrow \Spec(R), g:\mathcal{Y}\longrightarrow
\Spec(R).$$

One knows \cite[5.9.3]{Ka-RLS} that, for any
finitely generated
subring $R\subset \C$, there exists an integer $N\ge 1$ such that
for all primes
$\ell$ which are prime to N, there exists a finite extension
$E/\Q_{\ell}$, with ring of integers $\Oz$ and
an injective ring homomorphism from $R$ to $\Oz$. Fix one
such prime number $\ell$, which we choose larger than
both $\dim(X)$ and $\dim(Y)$, and one such inclusion of $R$ into $\Oz$.

Over $\Spec(R[1/\ell])$, the $\Q_{\ell}$-sheaves $R^{i}f_{\star}\Q_{\ell}$ and
$R^{i}g_{\star}\Q_{\ell}$ are lisse, and pure of weight $i$
\cite[3.3.9]{De-Weil II}. By the 
Lefschetz Trace Formula and proper base change, for each finite field 
$k$, and for each $k$-valued point $\phi$ of $\Spec(R[1/\ell])$, we
have
$$\Zeta(\X_{k,\phi}/k,t) =
\prod_{i}\det(1-t\Frob_{k,\phi}|R^{i}f_{\star}\Q_{\ell})^{(-1)^{i+1}}$$
and
$$\Zeta(\mathcal{Y}_{k,\phi}/k,t)=
\prod_{i}\det(1-t\Frob_{k,\phi}|R^{i}g_{\star}\Q_{\ell})^{(-1)^{i+1}}$$
By the assumed zeta-equivalence, we have, for
each finite field 
$k$, and for each $k$-valued point $\phi$ of $\Spec(R[1/\ell])$, the
equality of rational functions
$$\Zeta(\X_{k,\phi}/k,t)=\Zeta(\mathcal{Y}_{k,\phi}/k,t).$$
Separating the reciprocal zeroes and poles by absolute value, we
infer by purity that for every $i$, we have
$$\det(1-t\Frob_{k,\phi}|R^{i}f_{\star}\Q_{\ell}) =
\det(1-t\Frob_{k,\phi}|R^{i}g_{\star}\Q_{\ell}).$$

Therefore by Chebotarev the virtual semisimple representations of $\pi_{1}(\Spec(R[1/\ell]))$ 
given by $(R^{i}f_{\star}\Q_{\ell})^{ss}$ and $(R^{i}g_{\star}\Q_{\ell})^{ss}$
are equal:
$$(R^{i}f_{\star}\Q_{\ell})^{ss}\cong (R^{i}g_{\star}\Q_{\ell})^{ss}.$$

Now make use of the inclusion of $R$ into
$\Oz$, which maps $R[1/\ell]$ to $E$. 
The pullbacks $X_{\Oz}$ and $\mathcal{Y}_{\Oz}$ of $\X/R$ and $\mathcal{Y}/R$ to $\Oz$
are proper and smooth over $\Oz$. Thus their generic fibres, $X_{E}$ 
and $\mathcal{Y}_{E}$ are projective and smooth over $E$, of dimension
strictly less than $\ell$, and they have {\it good reduction}. Via the
chosen map from $\Spec(E)$ to $\Spec(R[1/\ell])$, we may pull back 
the representations $R^{i}f_{\star}\Q_{\ell})$ and
$R^{i}g_{\star}\Q_{\ell})$ 
of $\pi_{1}(\Spec(R[1/\ell]))$ to
$\pi_{1}(\Spec(E))$, the galois group $Gal_{E}:= Gal(E^{sep}/E)$.
Their pullbacks are the etale cohomology groups $H^{i}(X_{E^{sep}},
\Q_{\ell})$ and $H^{i}(\mathcal{Y}_{E^{sep}},
\Q_{\ell})$ respectively, viewed as representations of $Gal_{E}$. These
representations of $Gal_{E}$ need not be semisimple, but 
their semisimplifications are isomorphic:
$$H^{i}(X_{E^{sep}},\Q_{\ell})^{ss}\cong H^{i}(\mathcal{Y}_{E^{sep}},
\Q_{\ell})^{ss}.$$
By a fundamental result of Fontaine-Messing \cite[Theorems A and
B]{F-M} (which applies in the case
of good reduction, $E/\Q_{\ell}$ unramified, and dimension less than $\ell$) 
and Faltings
\cite[4.1]{Fal} (which
treats the general case, of a projective smooth generic fibre), 
we know that $H^{i}(X_{E^{sep}},
\Q_{\ell})$ and $H^{i}(\mathcal{Y}_{E^{sep}},\Q_{\ell})$
are
Hodge-Tate representations of $Gal_{E}$, with Hodge-Tate numbers exactly the Hodge 
numbers of
the complex projective smooth varieties $X$ and $Y$ respectively
(i.e., the dimension of the $Gal_{E}$-invariants in $H^{a}(X_{E^{sep}},
\Q_{\ell})(b)\otimes\C_{\ell}$ under the semilinear action of $Gal_{E}$ 
is the Hodge number
$H^{b,a-b}(X)$, and similarly for $Y$). By an elementary argument of
Wang \cite[5.1]{Wang}, 
the semisimplification of a Hodge-Tate
representation is also Hodge-Tate, with the same Hodge-Tate numbers.
So the theorem of
Fontaine-Messing and Faltings tells us that for all $i$, $H^{i}(X)$ and
$H^{i}(Y)$ have the same Hodge numbers. This is
precisely the required statement, that
$E(X;x,y)=E(Y;x,y)$.
\end{proof}

The reader may wonder why we introduced the notion of being fibrewise 
polynomial-count, for an element $\gamma \in K_{0}(\Sch/\C)$. In fact, 
this notion is entirely superfluous, as shown by the following Theorem.

\begin{theorem}Suppose $\gamma \in K_{0}(\Sch/\C)$ is fibrewise
   polynomial-count. Then it is strongly polynomial-count.\label{strongly}
   \end{theorem}
\begin{proof}Write $\gamma$ as $[X]-[Y]$, with $X$ and $Y$ projective smooth 
   $\C$-schemes. Repeat the first paragraph of the proof of the
   previous theorem. Extending $R$ if necessary, we may assume that
   the element $[\X/R]-[\mathcal{Y}/R]\in K_{0}(\Sch/R)$ is fibrewise 
   polynomial-count. So for each finite field $k$ and each ring
   homomorphism $\phi: R \longrightarrow k$, there exists a
   polynomial $P_{k,\phi}=\sum_{n}a_{n,k,\phi}t^{n} \in \Z[t]$ such
   that  
$$\Zeta(\X_{k,\phi}/k,t)/\Zeta(\mathcal{Y}_{k,\phi}/k,t)=\prod_{n}(1-(\#k)^{n}t)^{-a_{n,k,\phi}}.$$
 Writing the cohomological expressions of the zeta functions and
 using purity, we see that the coefficient $a_{n,k,\phi}$ is just the
 difference of the $2n$'th $\ell$-adic Betti numbers of
 $\X_{k,\phi}\otimes {\overline k}$
 and $\mathcal{Y}_{k,\phi}\otimes {\overline k}$, which is in turn the
 difference of the ranks of the two lisse sheaves 
 $R^{2n}f_{\star}\Q_{\ell}$ and $R^{2n}g_{\star}\Q_{\ell}$.
 This last difference is independent of the particular choice of $(k,
 \phi)$.
   \end{proof}

\end{document}